\documentclass[leqno]{article}
\usepackage{latexsym,amsmath,amsfonts,mathrsfs,amssymb,amsthm}
\usepackage{yhmath}
\usepackage[usenames]{color}
\usepackage{graphicx}
\def\R{\mathbb R}
\def\N{\mathbb N}
\def\C{\mathbb C}
\def\Z{\mathbb Z}
\def\s{\sharp}
\def\a{\alpha}
\def\b{\beta}
\def\e{\epsilon}
\def\d{\delta}
\def\t{\tilde}
\def\pa{\partial}
\def\Y{\mathbb Y}
\def\T{\mathbb T}
\def\P{\mathbb P}
\def\F{{\cal F} }
\def\oF{\overline{\cal F} }
\def\oB{\overline B}
\def\H{{\cal H} }

\def\mC{\mathfrak C}
\def\mA{\mathfrak A}
\def\mc{\mathfrak c}
\def\mU{\mathfrak U}
\def\mS{\mathfrak S}
\def\ms{\mathfrak s}
\def\cU{{\cal U}}

\def\mI{\mathfrak I}

\def\be{\begin{equation}}
\def\ee{\end{equation}}
\def\bs{\backslash}
\def\wg{\wedge}

\def\qed{\hfill$\Box$\bigskip}
\def\nd{\noindent \textbf{Proof. }}

\def\tb{\textcolor{blue}}
\numberwithin{equation}{section}
\newtheorem{thm}{Theorem}[section]
\newtheorem{lem}[thm]{Lemma}
\newtheorem{pro}[thm]{Proposition}
\newtheorem{defn}[thm]{Definition}
\newtheorem{cor}[thm]{Corollary}
\newtheorem{rem}[thm]{Remark}
\begin{document}
\bigskip

\begin{center}{\Large \textbf{Minimality for unions of 2-dimensional minimal cones with non-isolated singularities}}
\end{center}

\bigskip

\centerline{\large Xiangyu Liang}

%

\vskip 1cm

\centerline {\large\textbf{Abstract.}}

In this article we prove that for a large class of 2-dimensional minimal cones (including almost all 2-dimensional minimal cones that we know), the almost orthogonal union of any two of them is still a minimal cone. Comparing to existing results for minimality of almost orthogonal union of planes \cite{2p,2ptopo}, here we are dealing with unions of cones with non isolated singularities, which results in a series of essential difficulties, and new ideas are required. 

The proof in this article can be generalized to other types of minimalities, e.g. topological minimality, Reifenberg minimality, etc..

\bigskip

\textbf{AMS classification.} 28A75, 49Q20, 49K21, 49K99

\bigskip

\textbf{Key words.} Plateau's problem, Classification of Singularities, Minimal sets, Non isolated singularities, Almost orthogonal unions, Regularity, Uniqueness, Hausdorff measure

\tableofcontents

\section{Introduction}

This paper is the core of a series of papers (\cite{stablePYT, uniquePYT, stableYXY}), which aims at finding new families of singularities for 2-dimensional Almgren minimal sets, by taking unions of known singularities in transversal directions. This is motivated by the attempt to classify all possible singularities for 2-dimensional Almgren minimal sets in $\R^n$.

 In this paper we prove that for a large class of 2-dimensional minimal cones (including almost all 2-dimensional minimal cones that we know), the almost orthogonal union of any two of them is still a minimal cone. 

The notion of Almgren minimality was introduced by Almgren \cite{Al76} to modernize Plateau's problem, which aims at understanding physical objects that minimize the area while spanning a given boundary. The study of regularity and existence for these sets is one of the centers of interest in geometric measure theory. 

Among all other notions of minimality for Plateau's problem (such as currents, varifolds, minimal surfaces, see works of De Giorgi, Douglas, Federer, Fleming, Reifenberg, etc.), Almgren's notion of minimal sets gives a very good description of the local behavior for soap films, in the sense that the list of all possible singularities (''tangent objects'') for 2-dimensional Almgren minimal sets in 3-dimensional ambient spaces coincides exactly with what people can observe in soap film experiments.

So let us say more about Almgren's minimal sets. 

Briefly, given any open set $U\subset\R^n$ ($U$ is the domain in which our minimal set lives), a relatively closed set $E$ in $U$ is $d-$dimensional Almgren minimal in $U$ if its $d-$dimensional Hausdorff measure could not be decreased by any Lipschitz deformation with compact support in $U$. (See Section \tb{2.1} for the precise definition.) In other words, $E$ is Almgren minimal if it minimizes Hausdorff measure among all its Lipschitz deformation with compact support in $U$. (For soap films, the dimension $d$ of the set is 2, and the ambient dimension $n=3$).

\bigskip

Note that we are working in the setting of sets, hence the point of view here is very different from those of minimal surfaces and mass minimizing currents under certain boundary conditions. Here we see no algebraic or differential structure from the definition, and no equational description, hence traditional tools on functional analysis and PDE do not work directly. On the other hand since we are talking about absolute minima, not just stationary points for some area functional, hence we are expecting better regularity. Comparing to the large number of results in the theory of mass minimizing currents, or classical minimal surfaces, in our case, very few results of regularity are known. (On the other hand, there have been a number of progress on the existence for minimizers when we work with sets, see \cite{Fv, Fang16, DGM17, FS17} etc. )

\medskip

The first regularity results for Almgren minimal sets have been given by Frederick Almgren \cite{Al76} (rectifiability, Ahlfors regularity in arbitrary dimension), then generalized by Guy David and Stephen Semmes \cite{DS00} (uniform rectifiability, big pieces of Lipschitz graphs), Guy David \cite{GD03} (minimality of the limit of a sequence of minimizers). 

Since Almgren minimal sets are rectifiable and Ahlfors regular, they admit a tangent plane at almost every point. But our main interest is to study those points where there is no tangent plane. These points are called singular points. 

\medskip

A first finer description of the interior regularity for Almgren minimal sets is due to Jean Taylor, who gave in \cite {Ta} an essential regularity theorem for 2-dimensional minimal sets in 3-dimensional ambient spaces: if $E$ is a minimal set of dimension 2 
in an open set of $\R^3$ , then  every point $x$ of $E$ has a neighborhood where $E$ can be parametrized (modulo a negligible set) through a $C^1$ diffeomorphism by an Almgren minimal cone (that is, an Almgren minimal set which is also a cone. We will call them minimal cones in the rest part of the paper). It might be worth mentioning that the list of 2 dimensional minimal cones in $\R^3$ coincides exactly with all possible singularities that we can observe in soap films.

In \cite{DJT}, Guy David generalized Jean Taylor's theorem to 2-dimensional Almgren minimal sets in $\R^n$, but with a local bi-H\"older parametrization, that is, every point $x$ of $E$ has 
a neighborhood where $E$ can be parametrized through a bi-H\"older diffeomorphism by a minimal cone $C$ (but the minimal cone might not be unique). In addition, in \cite{DEpi}, David also proved that, if this minimal cone $C$ satisfies an epiperimetric condition, we will have the $C^1$ parametrization (called $C^1$ regularity). In particular, the tangent cone of $E$ at the point $x\in E$ exists and is a minimal cone, and the blow-up limit (see Definition \tb{2.11}) of $E$ at $x$ is unique.

Thus, the study of singular points is transformed into the classification of singularities, i.e., into looking for a list of minimal cones. Besides, getting such a list would also help deciding locally what kind (i.e. $C^1$ or bi-H\"older) of perimetrization by a minimal cone can we get.

\bigskip

In $\R^3$, the list of 2-dimensional minimal cones has been given by several mathematicians a century ago. (See for example \cite{La} or \cite{He}). They are, modulo isometry: a plane (which we also call a $\P$ set), a $\Y$ set (the union of 3 half planes that meet along a straight line where they make angles of 120 degrees. This straight line is called its spine), and a $\T$ set (the cone over the 1-skeleton of a regular tetrahedron centered at the origin). See the pictures below. 

\centerline{\includegraphics[width=0.16\textwidth]{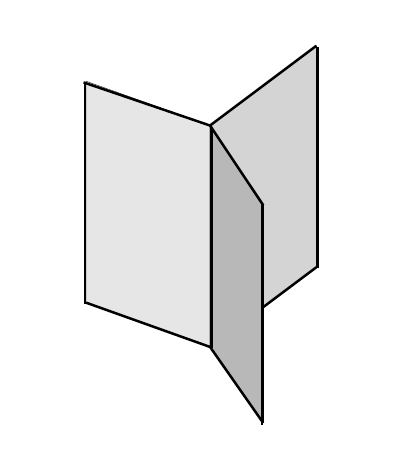} \hskip 2cm\includegraphics[width=0.2\textwidth]{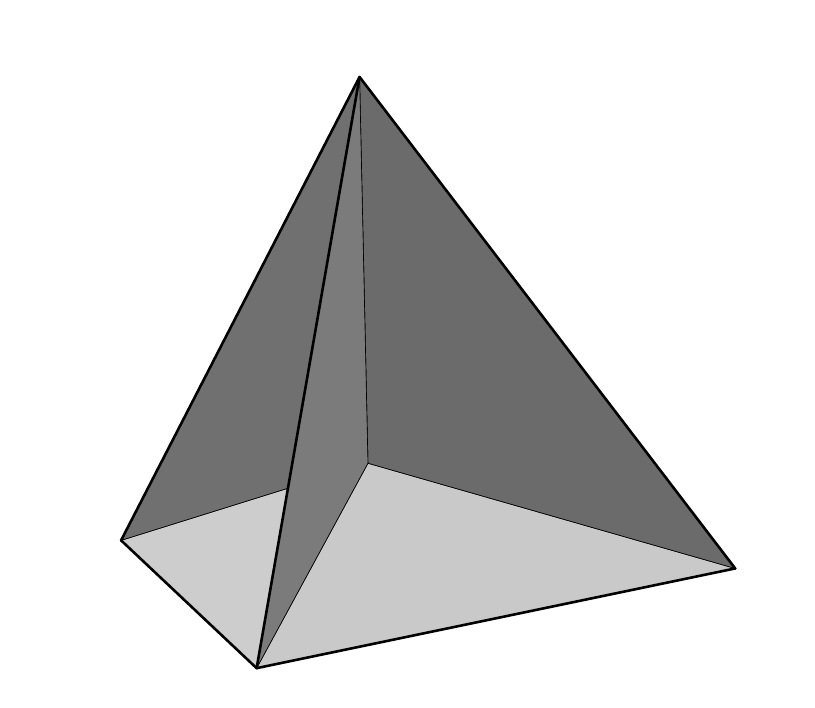}}
\nopagebreak[4]
\hskip 4cm a $\Y$ set\hskip 3.9cm  a $\T$ set

In higher dimensions, even in dimension 4, the list of minimal cones is still very far from clear. Except for the three minimal cones that already exist in $\R^3$, the only 2-dimensional minimal cones that were known before this paper is the product of two 1-dimensional $\Y$ sets \cite{YXY}, and unions of several almost orthogonal planes (cf.\cite{2p} and \cite{2ptopo}). The later gives continuous one-parameter families of minimal cones, and this gives already an interesting phenomenon that does not occur in $\R^3$----in $\R^3$, no small perturbation of any minimal cones ever preserves the minimality, and moreover, each minimal cone admits a different topology from the others. 

Note that in case of unions of two orthogonal planes, the minimality comes easily from an algebraic lemma (cf. for example Lemma 5.2 of \tb{\cite{Mo84}}). But in case of almost orthogonal planes, the minimality is far from intuitive, and the proof of minimality is already quite envolved.

In this article we aims at proving the minimality for unions of almost orthogonal 2-dimensional minimal cones, and let us first state the main theorem:

\begin{thm}\label{main}[minimality of the union of two almost orthogonal 2-dimensional minimal cones]Let $K_0^1\subset \R^{n_1}$ and $K_0^2\subset \R^{n_2}$ be two 2-dimensional Almgren unique minimal cones. We suppose in addition that both cones are sliding stable. Then there exists $0<\theta_0<\frac\pi2$, such that if $K^1$ and $K^2$ are two copies (centered at 0) of $K_0^1$ and $K_0^2$ in $\R^{n_1+n_2}$, with $\a_{K^1,K^2}\ge\theta_0$, 
then the union $K^1\cup K^2$ is a minimal cone.
\end{thm}

\begin{rem}
$1^\circ$
The angle $\a_{K^1,K^2}$ between two minimal cones $K^1$ and $K^2$ is defined to measure the orthogonality between the two sets. When $\a=0$ the two sets are orthogonal, and almost orthogonal when $\a$ is small. See Definition \tb{2.40};

$2^\circ$ The Almgren uniqueness says roughly that a minimal set is the only one that attains the minimal measure in the class of its deformations and their limits. See the beginning of \tb{Section 3} for the precise definition.

$3^\circ$ The condition of sliding stability means that the set is of minimal measure among deformations if we allow the boundary of the set to move a little bit along the boundary of the domain. See Definition \tb{2.29}.
\end{rem}
 
Comparing to existing results for minimality of almost orthogonal union of planes \cite{2p,2ptopo}, here we are dealing with unions of cones with non isolated singularities, which means that we cannot decompose things into regular parts and turn ourselves into treating problems with non singularities (for which many well known tools in minimal graphs and differential geometry can easily apply). As a result, a series of essential difficulties occurs throughout the proof, such as the lack of good local parametrization, the lack of translation invariance, the lack of good 1-Lipschitz projection, the lack of Almgren uniqueness, etc. Hence new ideas are required. We list here the main ingredient to overcome some of the above difficulties:

$1^\circ$ The proof of the uniqueness theorem for the orthogonal union:

Like in \cite{2p}, to prove uniqueness, we mainly use the fact that given any minimizer other than the orthogonal union, the tangent plane is analytic almost everywhere. If the two components of the union are planes, things turn to graphs, and complex analysis can essentially do almost everything. But here we treat components with singularities, hence we have to carefully control the behavior of singular sets, which are of codimension 1, for which the property of ''analytic tangent plane almost everywhere'' says nothing. We have to analyse the manner under which the regular parts can meet, to deal with this codimension 1 set.

$2^\circ$ Flatten boundary near minimal cones:

To estimate the measure near the non isolated singular sets ($\Y$ sets), again because there is no 1-Lipschitz projection, we construct a particular regions $\cU$ (in stead of simply balls) for each 2-dimensional minimal cone (see Section \tb{5}), to insure that the variation of the measure of the minimal cone while sliding along the boundary of the region is zero (cf. \cite{stablePYT} Theorem 3.1). In this particularly defined region, we can also prove the sliding stability for all known 2-dimensional minimal cones (\cite{stablePYT,stableYXY}), which is a necessary condition in our main theorem. 

$3^\circ$ Fine estimates of measure based on topological invariant:

In contrary to minimal sets near a plane, we cannot use minimal graph to immitate the behavior of a minimal set close to a $\Y$ set, 
the fine estimate of measure for such sets are mainly done under the help of topological invariant, like connectness and homology to control the topology near the sets (see Section \tb{6}), and then to control the measure (see Section \tb{7}).

$4^\circ$ Reduction from a nearly conic minimal set to any of its small part:

When we want to decompose a minimal set near the union of $K^1$ and $K^2$, due to the non translation invariance for minimal cones other than planes, it may occur that at smaller and smaller scales, our set goes gradually near the union of translated $K^1$ and $K^2$ with different centers. The whole Section \tb{8} will be devoted to solve this. Briefly, we introduce a reduction procedure, and reduce this case step by step to the case of the union of two planes. As a key step for this reduction, we construct retractions from a minimal cone to any arbitrary small part, and further we do the same for any minimal set that is very near a cone.

Another problem issued from the non translation invariance is the measure stability for an arbitrary 2-dimensinoal minimal cone. We prove it in \cite{stablePYT}.

$5^\circ$ We have also changed the whole workspace: 

In the proof of minimality of almost orthogonal union of planes, we carried on a series of processes that pass to the limit. By definition, a minimal cone only minimizes measure among all its deformations, but not among limits of its deformations. But one can see in the proof of \cite{2p,2ptopo}, that we still work in the frame of deformations, because if the minimal cone is a plane, taking limits is not a real issue. In fact, some manipulations on orthogonal projections (which are 1-Lipschitz) to planes are enough to show that a plane minimizes the measure among a much larger class of sets, including limits of its deformations, and even if we allow the deformations to move the boundary orthogonally with respect to the plane. 

For an arbitrary 2-dimensional minimal cone, when there is no 1-Lipschitz neighorhood retraction, we are forced to pay our attention to the much larger class of sets-- the class of all limits of deformations. We mainly use homology property to analyse the local structure for limits; and to further estimate the measure, we prove an upper-semi-continuity result (cf. Theorem \tb{4.1} of \cite{uniquePYT}). This plays an important role in the whole proof.

There are still a number of other problems during the analysis of local structure and measure estimates around singularities, as we will see throughout the proof. 

\bigskip 

After all the above efforts, there are still two issues that we are not able to cover: The Almgren uniqueness, and the sliding stability for the two minimal cones $K^1$ and $K^2$. That is why we add these two assumptions in our main theorem. We do not know whether these are necessary condition for the minimality of almost orthogonal union of $K^1$ and $K^2$. But on the other hand, so far for all of the 2-dimensional minimal cones we know, they all satisfy these two conditions. See \cite{uniquePYT}, \cite{stablePYT} and \cite{stableYXY}.

We also remark that for minimal cones with calibrations or paired calibrations, or satisfy some topological condition, they automatically satisfy the sliding stability. Sets with calibrations and paired calibrations are very likely to satisfy the Almgren uniqueness as well. See proof in \cite{stablePYT} and \cite{stableYXY}.

\bigskip

The whole article before the last section will be devoted to the proof of Theorem \ref{main}. In the last section we give some perspectives. We prove that the almost orthogonal union of two sliding stable Almgren unique 2-dimensional minimal cones is still sliding stable and Almgren unique, so that we can continue take the union of these newly found minimal cones with others to generate new families of Almgren minimal cones. We also announce a theorem on the topological minimality of the almost orthogonal union of two topological minimal cones, the proof of which is simpler, because with the topology condition, things pass to the limit directly.

\textbf{Acknowledgement:} This work was supported by China's Recruitement Program of Global Experts, School of Mathematics and Systems Science, Beihang University. 

\section{Definitions and preliminaries}

\subsection{Some useful notation}

$[a,b]$ is the line segment with end points $a$ and $b$;

$\overrightarrow{ab}$ is the vector $b-a$;

$R_{ab}$ denotes the half line issued from the point $a$ and passing through $b$;

$B(x,r)$ is the open ball with radius $r$ and centered on $x$;

$\overline B(x,r)$ is the closed ball with radius $r$ and center $x$;

$\H^d$ is the Hausdorff measure of dimension $d$ ;

$d_H(E,F)=\max\{\sup\{d(y,F):y\in E\},\sup\{d(y,E):y\in F\}\}$ is the Hausdorff distance between two sets $E$ and $F$. 

For any subset $K\subset \R^n$, the local Hausdorff distance in $K$ $d_K$ between two sets $E,F$ is defined as $d_K(E,F)=\max\{\sup\{d(y,F):y\in E\cap K\},\sup\{d(y,E):y\in F\cap K\}\}$;

For any open subset $U\subset \R^n$, let $\{E_n\}_n$, $F$ be closed sets in $U$, we say that $F$ is the Hausdorff limit of $\{E_n\}_n$, if for any compact subset $K\subset U$, $\lim_n d_K(E_n,F)=0$;

$d_{x,r}$ : the relative distance with respect to the ball $B(x,r)$, is defined by
$$ d_{x,r}(E,F)=\frac 1r\max\{\sup\{d(y,F):y\in E\cap B(x,r)\},\sup\{d(y,E):y\in F\cap B(x,r)\}\}.$$

For any (affine) subspace $Q$ of $\R^n$, and $x\in Q$, $r>0$, $B_Q(x,r)$ stands for $B(x,r)\cap Q$, the open ball in $Q$.

For any subset $E$ of $\R^n$ and any $r>0$, we call $B(E,r):=\{x\in \R^n: dist (x,E)<r\}$ the $r$ neighborhood of $E$.

If $E$ is a $d$-rectifiable set, denote by $T_xE$ the tangent plane (if it exists and is unique) of $E$ at $x$.
$\H^d$ is the Hausdorff measure of dimension $d$ ;

For any $d\le n$, any abelien group $G$, and any subset $E\subset \R^n$, $H_d(E,G)$ denote the $d$-th singular homological group of $E$ with coefficient in the group $G$.

\subsection{Basic definitions and notations about minimal sets}

In the next definitions, fix integers $0<d<n$. We first give a general definition for minimal sets. Briefly, a minimal set is a closed set which minimizes the Hausdorff measure among a certain class of competitors. Different choices of classes of competitors give different kinds of minimal sets.

\begin{defn}[Minimal sets]Let $0<d<n$ be integers. Let $U\subset \R^n$ be an open set. A relative closed set $E\subset U$ is said to be minimal of dimension $d$ in $U$ with respect to the competitor class $\mathscr F$ (which contains $E$) if 
\be \H^d(E\cap B)<\infty\mbox{ for every compact ball }B\subset U,\ee
and
\be \H^d(E\bs F)\le \H^d(F\bs E)\ee
for any competitor $F\in\mathscr F$.
\end{defn}

\begin{defn}[Almgren competitor (Al competitor for short)] Let $E$ be relatively closed in an open subset $U$ of $\R^n$. An Almgren competitor for $E$ is an relatively closed set $F\subset U$ that can be written as $F=\varphi_1(E)$, where $\varphi_t:U\to U,t\in [0,1]$ is a family of continuous mappings such that 
\be \varphi_0(x)=x\mbox{ for }x\in U;\ee
\be\mbox{ the mapping }(t,x)\to\varphi_t(x)\mbox{ of }[0,1]\times U\mbox{ to }U\mbox{ is continuous;}\ee
\be\varphi_1\mbox{ is Lipschitz,}\ee
  and if we set $W_t=\{x\in U\ ;\ \varphi_t(x)\ne x\}$ and $\widehat W=\bigcup_{t\in[0.1]}[W_t\cup\varphi_t(W_t)]$,
then
\be \widehat W\mbox{ is relatively compact in }U.\ee
 
Such a $\varphi_1$ is called a deformation in $U$, and $F$ is also called a deformation of $E$ in $U$.
\end{defn}

For future convenience, we also have the following more general definition:

\begin{defn}Let $U\subset \R^n$ be an open set, and let $E\subset \R^n$ be a closed set (not necessarily contained in $U$). We say that $E$ is minimal in $U$, if $E\cap U$ is minimal in $U$. A closed set $F\subset \R^n$ is called a  deformation of $E$ in $U$, if $F=(E\bs U)\cup \varphi_1(E\cap U)$, where $\varphi_1$ is a deformation in $U$.
\end{defn}

We denote by $\F(E,U)$ the class of all deformations of $E$ in $U$. In the article, we need to use Hausdorff limits in $\F(E,U)$. However, if we regard elements of $\F(E,U)$ as sets in $\R^n$, and take the Hausdorff limit, the limit may have positive measure on $\partial U\bs E$. In other words, sets in $\F(E,U)$ may converge to the boundary. We do not like this. Hence we let $\oF(E,U)$ be the class of Hausdorff limits of sequnces in $\F(E,U)$ that do not converge to the boundary. That is: we set
\be \begin{split}\oF(E,U)&=\{F\subset\bar U: \exists \{E_n\}_n\subset\F(E,U)\mbox{ such that }d_K(E_n,F)\to 0\\
&\mbox{ for all compact set }K\subset \R^n, 
\tb{(2.1)}\mbox{ holds for }F\mbox{, and }\H^d(F\cap\partial U\bs E)=0\}.\end{split}\ee

It is easy to see that both classes are stable with respect to Lipschitz deformations in $U$. 

\begin{defn}[Almgren minimal sets]
Let $0<d<n$ be integers, $U$ be an open set of $\R^n$. An Almgren-minimal set $E$ in $U$ is a minimal set defined in Definition \tb{2.1} while taking the competitor class $\mathscr F$ to be the class of all Almgren competitors for $E$.\end{defn}

For the need of our future argument, we also have the following definition:

\begin{defn}Let $0<d<n$ be integers, $U$ be an open set of $\R^n$. A closed set $E\subset \R^n$ is said to be Almgren minimal in $U$, if $E\cap U$ is Almgren minimal in $U$.
\end{defn}

%
%
%
%

\begin{defn}[Topological competitors] Let $G$ be an abelian group. Let $E$ be a closed set in an open domain $U$ of $\R^n$. We say that a closed set $F$ is a $G$-topological competitor of dimension $d$ ($d<n$) of $E$ in $U$, if there exists a convex set $B$ such that $\bar B\subset U$ such that

1) $F\bs B=E\bs B$;

2) For all Euclidean $n-d-1$-sphere $S\subset U\bs(B\cup E)$, if $S$ represents a non-zero element in the singular homology group $H_{n-d-1}(U\bs E; G)$, then it is also non-zero in $H_{n-d-1}(U\bs F;G)$.
We also say that $F$ is a $G$-topological competitor of $E$ in $B$.

When $G=\Z$, we usually omit $\Z$, and say directly topological competitor.
\end{defn}

And Definition \tb{2.1} gives the definition of $G$-topological minimizers in a domain $U$ when we take the competitor class to be the classe of $G$-topological competitors of $E$.

The simplest example of a topological minimal set is a $d-$dimensional plane in $\R^n$.  

\begin{pro}[cf.\cite{topo} Proposition 3.7 and Corollary 3.17]  Let $U$ be a domain in $\R^n$.

$1^\circ$ Let $E\subset U$ be relatively closed. Then for any $d<n$, every Almgren competitor of $E$ in $U$ is a $G$-topological competitor of $E$ of dimension $d$ in $U$.

$2^\circ$ All $G$-topological minimal sets in $U$ are Almgren minimal in $U$.
\end{pro}

\begin{defn}[reduced set] Let $U\subset \R^n$ be an open set. For every closed subset $E$ of $U$, denote by
\be E^*=\{x\in E\ ;\ \H^d(E\cap B(x,r))>0\mbox{ for all }r>0\}\ee
 the closed support (in $U$) of the restriction of $\H^d$ to $E$. We say that $E$ is reduced if $E=E^*$.
\end{defn}

It is easy to see that
\be \H^d(E\bs E^*)=0.\ee
In fact we can cover $E\bs E^*$ by countably many balls $B_j$ such that $\H^d(E\cap B_j)=0$.

\begin{rem}
 It is not hard to see that if $E$ is Almgren minimal (resp. $G$-topological minimal), then $E^*$ is also Almgren minimal (resp. $G$-topological minimal). As a result it is enough to study reduced minimal sets. An advantage of reduced minimal sets is, they are locally Ahlfors regular (cf. Proposition 4.1 in \cite{DS00}). Hence any approximate tangent plane of them is a true tangent plane. Since minimal sets are rectifiable (cf. \cite{DS00} Theorem 2.11 for example), reduced minimal sets admit true tangent $d$-planes almost everywhere.
\end{rem}

\begin{rem}The notion of (Almgren or $G$-topological) minimal sets does not depend much on the ambient dimension. One can easily check that $E\subset U$ is $d-$dimensional Almgren minimal in $U\subset \R^n$ if and only if $E$ is Almgren minimal in $U\times\R^m\subset\R^{m+n}$, for any integer $m$. The case of $G$-topological minimality is proved in \cite{topo} Proposition 3.18.\end{rem}

\subsection{Regularity results for minimal sets}

We now begin to give regularity results for minimal sets. They are in fact regularity results for Almgren minimal sets, but they also hold for all $G$-topological minimizers, after Proposition \tb{2.7}. By Remark \tb{2.9}, from now on all minimal sets are supposed to be reduced. 

\begin{defn}[blow-up limit] Let $U\subset\R^n$ be an open set, let $E$ be a relatively closed set in $U$, and let $x\in E$. Denote by $E(r,x)=r^{-1}(E-x)$. A set $C$ is said to be a blow-up limit of $E$ at $x$ if there exists a sequence of numbers $r_n$, with $\lim_{n\to \infty}r_n=0$, such that the sequence of sets $E(r_n,x)$ converges to $C$ for the local Hausdorff distance in any compact set of $\R^n$.
\end{defn}

\begin{rem}
 $1^\circ$ A set $E$ might have more than one blow-up limit at a point $x$. However it is not known yet whether this can happen to minimal sets. 
 
 When a set $E$ admits a unique blow-up limit at a point $x\in E$, denote this blow-up limit by $C_xE$.
 
 $2^\circ$ Let $Q\subset \R^n$ be any subpace, denote by $\pi_Q$ the orthogonal projection from $\R^n$ to $Q$. Then it is easy to see that if $E\subset \R^n$, $x\in E$, and $C$ is any blow-up limit of $E$ at $x$, then $\pi_Q(C)$ is contained in a blow-up limit of $\pi_Q(E)$ at $\pi_Q(x)$. 
\end{rem}

\begin{pro}[c.f. \cite{DJT} Proposition 7.31]Let $E$ be a reduced Almgren minimal set in an open set $U$ of $\R^n$, and let $x\in E$. Then every blow-up limit of $E$ at $x$ is a reduced Almgren minimal cone $F$ centred at the origin, and $\H^d(F\cap B(0,1))=\theta(x):=\lim_{r\to 0} r^{-d}\H^d(E\cap B(x,r)).$\end{pro}

An Almgren minimal cone is just a cone which is also Almgren minimal. We will call them minimal cones throughout this paper, since we will not talk about any other type of minimal cones. 

\begin{rem}$1^\circ$ The existence of the density $\theta(x)$ is due to the monotonicity of the density function $\theta(x,r):=r^{-d}\H^d(E\cap B(x,r))$ at any given point $x$ for minimal sets. See for example \cite{DJT} Proposition 5.16.

$2^\circ$ After the above proposition, the set $\Theta(n,d)$ of all possible densities for points in a $d$-dimension minimal set in $\R^n$ coincides with the set of all possible densities for $d$-dimensional minimal cones in $\R^n$. When $d=2$, this is a very small set. For example, we know that $\pi$ is the density for a plane, $\frac 32\pi$ is the density for a $\Y$ set, and for any $n$, and any other type of 2-dimensional minimal cone in $\R^n$, its density should be no less than some $d_T=d_T(n)>\frac 32\pi$, by \cite{DJT} Lemma 14.12.

$3^\circ$ Obviously, a cone in $\R^n$ is minimal if and only if it is minimal in the unit ball, if and only if it is minimal in any open subset containing the origin.

$4^\circ$ For future convenience, we also give the following notation: let $U\subset \R^n$ be a convex domain containing the origin. A set $C\subset U$ is called a cone in $U$, if it is the intersection of a cone with $U$.
\end{rem}

We now state some regularity results on 2-dimensional Almgren minimal sets. 

\begin{defn}[bi-H\"older ball for closed sets] Let $E$ be a closed set of Hausdorff dimension 2 in $\R^n$. We say that $B(0,1)$ is a bi-H\"older ball for $E$, with constant $\tau\in(0,1)$, if we can find a 2-dimensional minimal cone $Z$ in $\R^n$ centered at 0, and $f:B(0,2)\to\R^n$ with the following properties:

$1^\circ$ $f(0)=0$ and $|f(x)-x|\le\tau$ for $x\in B(0,2);$

$2^\circ$ $(1-\tau)|x-y|^{1+\tau}\le|f(x)-f(y)|\le(1+\tau)|x-y|^{1-\tau}$ for $x,y\in B(0,2)$;

$3^\circ$ $B(0,2-\tau)\subset f(B(0,2))$;

$4^\circ$ $E\cap B(0,2-\tau)\subset f(Z\cap B(0,2))\subset E.$  

We also say that B(0,1) is of type $Z$ for $E$.

We say that $B(x,r)$ is a bi-H\"older ball for $E$ of type $Z$ (with the same parameters) when $B(0,1)$ is a bi-H\"older ball of type $Z$ for $r^{-1}(E-x)$.
\end{defn}

\begin{thm}[Bi-H\"older regularity for 2-dimensional Almgren minimal sets, c.f.\cite{DJT} Thm 16.1]\label{holder} Let $U$ be an open set in $\R^n$ and $E$ a reduced Almgren minimal set in $U$. Then for each $x_0\in E$ and every choice of $\tau\in(0,1)$, there is an $r_0>0$ and a minimal cone $Z$ such that $B(x_0,r_0)$ is a bi-H\"older ball of type $Z$ for $E$, with constant $\tau$. Moreover, $Z$ is a blow-up limit of $E$ at $x$.
\end{thm}

\begin{defn}[point of type $Z$] 

$1^\circ$ In the above theorem, we say that $x_0$ is a point of type $Z$ (or $Z$ point for short) of the minimal set $E$. The set of all points of type $Z$ in $E$ is denoted by $E_Z$. 

$2^\circ$ In particular, we denote by $E_P$ the set of regular points of $E$ and $E_Y$ the set of $\Y$ points of $E$. Any 2-dimensional minimal cone other than planes and $\Y$ sets are called $\T$ type cone, and any point which admits a $\T$ type cone as a blow-up is called a $\T$ type point. Set $E_T=E\bs (E_Y\cup E_P)$ the set of all $\T$ type points of $E$.  Set $E_S:=E\bs E_P$ the set of all singular points in $E$.
\end{defn}

\begin{rem} Again, since we might have more than one blow-up limit for a minimal set $E$ at a point $x_0\in E$, the point $x_0$ might be of more than one type (but all the blow-up limits at a point are bi-H\"older equivalent). However, if one of the blow-up limits of $E$ at $x_0$ admits the``full-length'' property (see Remark \ref{ful}), then in fact $E$ admits a unique blow-up limit at the point $x_0$. Moreover, we have the following $C^{1,\a}$ regularity around the point $x_0$. In particular, the blow-up limit of $E$ at $x_0$ is in fact a tangent cone of $E$ at $x_0$.
\end{rem}

\begin{thm}[$C^{1,\a}-$regularity for 2-dimensional minimal sets, c.f. \cite{DEpi} Thm 1.15]\label{c1} Let $E$ be a 2-dimensional reduced minimal set in the open set $U\subset\R^n$. Let $x\in E$ be given. Suppose in addition that some blow-up limit of $E$ at $x$ is a full length minimal cone (see Remark \ref{ful}). Then there is a unique blow-up limit $X$ of $E$ at $x$, and $x+X$ is tangent to $E$ at $x$. In addition, there is a radius $r_0>0$ such that, for $0<r<r_0$, there is a $C^{1,\a}$ diffeomorphism (for some $\a>0$) $\Phi:B(0,2r)\to\Phi(B(0,2r))$, such that $\Phi(0)=x$ and $|\Phi(y)-x-y|\le 10^{-2}r$ for $y\in B(0,2r)$, and $E\cap B(x,r)=\Phi(X)\cap B(x,r).$ 

We can also ask that $D\Phi(0)=Id$. We call $B(x,r)$ a $C^1$ ball for $E$ of type $X$.
\end{thm}

\begin{rem}[full length, union of two full length cones $X_1\cup X_2$]\label{ful}We are not going to give the precise definition of the full length property. Instead, we just give some information here, which is enough for the proofs in this paper.

$1^\circ$ The three types of 2-dimensional minimal cones in $\R^3$, i.e. the planes, the $\Y$ sets, and the $\T$ sets, all verify the full-length property (c.f., \cite{DEpi} Lemmas 14.4, 14.6 and 14.27). Hence all 2-dimensional minimal sets $E$ in an open set $U\subset\R^3$ admits the local $C^{1,\a}$ regularity at every point $x\in E$. But this was known from \cite{Ta}.

$2^\circ$ (c.f., \cite{DEpi} Remark 14.40) Let $n>3$. Note that the planes, the $\Y$ sets and the $\T$ sets are also minimal cones in $\R^n$. Denote by $\mathfrak C$ the set of all planes, $\Y$ sets and $\T$ sets in $\R^n$. Let $X=\cup_{1\le i\le n}X_i\in \R^n$ be a minimal cone, where $X_i\in \mathfrak{C}, 1\le i\le n$, and for any $i\ne j$, $X_i\cap X_j=\{0\}$. Then $X$ also verifies the full-length property. 
\end{rem}

\begin{thm}[Structure of 2-dimensional minimal cones in $\R^n$, cf. \cite{DJT} Proposition 14.1] Let $K$ be a reduced 2-dimensional minimal cone in $\R^n$, and let $X=K\cap \partial B(0,1)$. Then $X$ is a finite union of great circles and arcs of great circles $C_j,j\in J$. The $C_j$ can only meet at their endpoints, and each endpoint is a common endpoint of exactly three $C_j$, which meet with $120^\circ$ angles. In addition, the length of each $C_j$ is at least $\eta_0$, where $\eta_0>0$ depends only on the ambient dimension $n$.
\end{thm}

An immediate corollary of the above theorem is the following:

\begin{cor}
$1^\circ$ If $C$ is a minimal cone of dimension 2, then for the set of regular points $C_P$ of $C$, each of its connected components is a sector. 

$2^\circ$ Let $E$ be a 2-dimensional minimal set in $U\subset \R^n$. Then $\bar E_Y=E_S$.

$3^\circ$ The set $E_S\bs E_Y$ is isolated. \end{cor}

As a consequence of the $C^1$ regularity for regular points and $\Y$ points, and Corollary \tb{2.22}, we have
\begin{cor}Let $E$ be an 2-dimensional Almgren minimal set in a domain $U\subset \R^n$. Then 

$1^\circ$ The set $E_P$ is open in $E$;

$2^\circ$ The set $E_Y$ is a countable union of $C^1$ curves. The endpoints of these curves are either in $E_T:=E_S\bs E_Y$, or lie in $\partial U$.  
\end{cor}

We also have a similar quantified version of the $C^{1,\a}$ regularity (cf. \cite{DJT} Corollary 12.25). In particular, we can use the distance between a minimal set and a $\P$ or a $\Y$ cone to controle the constants of the $C^{1,\a}$ parametrization. As a direct corollary, we have the following neighborhood deformation retract property for regular and $\Y$ points:

\begin{cor}There exists $\e_2=\e_2(n)>0$ such that the following holds : let $E$ be an 2-dimensional Almgren minimal set in a domain $U\subset \R^n$. Then 

$1^\circ$ For any $x\in E_P$, and any co-dimension 1 submanifold $M\subset U$ which contains $x$, such that $M$ is transversal to the tangent plane $T_xE+x$, if $r>0$ satisfies that $d_{x,r}(E, x+T_xE)<\e_2$, then $\H^1(B(x,r)\cap M\cap E)<\infty$, and $B(x,r)\cap M\cap E$ is a Lipschitz deformation retract of $B(x,r)\cap M$;

$2^\circ$ For any $x\in E_Y$, and any co-dimension 1 submanifold $M\subset U$ which contains $x$, such that $M$ is transversal to the tangent cone $C_xE+x$ and its spine, if $r>0$ satisfies that $d_{x,r}(E, x+C_xE)<\e_2$, then $\H^1(B(x,r)\cap M\cap E)<\infty$, and$B(x,r)\cap M\cap E$ is a Lipschitz deformation retract of $B(x,r)\cap M$.
\end{cor}

As for the regularity for minimal sets of higher dimensions, we know much less. But for points which admit a tangent plane (i.e. some blow up-limit on the point is a plane), we still have the $C^1$ regularity.

\begin{thm}[cf.\cite{2p} Proposition 6.4]\label{e1}For $2\le d<n<\infty$, there exists $\epsilon_1=\e_1(n,d)>0$ such that if $E$ is a $d$-dimensional reduced minimal set in an open set $U\subset\R^n$, with $B(0,2)\subset U$ and $0\in E$. Then if $E$ is $\epsilon_1$ near a $d-$plane $P$ in $B(0,1)$, then $E$ coincides with the graph of a $C^1$ map $f:P\to P^\perp$ in $B(0,\frac 34)$. Moreover $||\nabla f||_\infty<1$.
\end{thm}

\begin{rem}
$1^\circ$ This proposition is a direct corollary of Allard's famous regularity theorem for stationary varifold. See \cite{All72}.

$2^\circ$ After this proposition, a blow-up limit of a reduced minimal set $E$ at a point $x\in E$ is a plane if and only if the plane is the unique approximate tangent plane of $E$ at $x$.
\end{rem}

After Remark \tb{2.26}, for any reduced minimal set $E$ of dimension $d$, and for any $x\in E$ at which an approximate tangent $d$-plane exists (which is true for a.e. $x\in E$), $T_xE$ also denotes the tangent plane of $E$ at $x$, and the blow-up limit of $E$ at $x$. 

\subsection{Stability for 2-dimensional minimal cones}

In this subsection we discuss briefly the notions of stability for 2-dimensional minimal cones, in two senses: 
$1^\circ$ whether a minimal cone is still minimal if we allow to ''change'' slightly their boundary; $2^\circ$ whether the ''measure'' of a minimal cone stays the same if we translate it a little bit in a properly chosen domain.

In both cases, we need first to define the convex domain $\cU(K,\eta)$, associated to each minimal cone $K$ and small $\eta$. Briefly, if $K$ is a 2-dimensional minimal cone, then the domain is roughly the unit ball $B$, but we modify it slightly, so that at almost every point $x$ of $K\cap \partial \cU(K,\eta)$, $\vec{ox}$ is perpendicular to $\partial \cU(K,\eta)$ in a neighborhood of $x$.

So let us fix a 2-dimensional minimal cones $K$ in $\R^n$. Let $B$ denote the unit ball. Then by Theorem \tb{2.21}, $K\cap \partial B$ is  a union of arcs of great circles with only $\Y$ type junctions. Denote by $\{a_j, 1\le j\le m\}$ the set of $\Y$ points in $K\cap \pa B$. Again by Theorem \tb{2.21}, the length of each arc is no less than a certain $\eta_0$. So here let $\eta<10^{-2}\eta_0$, very small.

Define the $\eta$-convex domain for $K$ \be \cU(K,\eta)=\{x\in B: <x, y>< 1-\eta, \forall y\in K\mbox{ and }<x,a_j><1-2\eta, \forall 1\le j\le m\}\subset \R^n. \ee 

From the definition, we see directly that each $\cU(K,\eta)$ is obtained by "cutting off" some small part of the unit ball $B$. More precisely, we first take the unit ball $B$, then just like peelling an apple, we use a knife to peel a thin band (with width about $2\sqrt \eta$) near the net $K\cap \pa B$. Then after this operation, the ball $B$ stays almost the same, except that near the set $K$, the boundary surface will be a thin cylinderical surface. This is the condition ''$<x, y>< 1-\eta, \forall y\in K$''. Next we turn to the singular points $a_k$: they are isolated, so we make one cut at each point, perpendicular to the radial direction, to get a small planar surface near each $a_j$, of diameter about $4\sqrt \eta$. This follows from the condition ''$<x,a_j><1-2\eta, \forall 1\le j\le m $''. Note that since $\eta<10^{-2}\eta_0$, each cut at these $\Y$ points does not intersect each other. 

We will discuss the precise structure of $\cU(K,\eta)$ in more details in Section \tb{5}. Here we continue the definitions of stability.

\begin{defn}Let $U$ be an open subset of $\R^n$, let $E\subset \bar U$ be closed. For $\d>0$, a \tb{$\d$-sliding Lipschitz deformation} of $E$ in $\bar U$ is a set $F\subset \bar U$ that can be written as $F=\varphi_1(E)$, where $\varphi_t: \bar U\to \bar U$ is a family of continuous mappings such that
\be \varphi_0(x)=x\mbox{ for }x\in \bar U;\ee
\be \mbox{ the mapping }(t,x)\to \varphi_t(x)\mbox{ of }[0,1]\times \bar U\mbox{ to }\bar U\mbox{ is continuous};\ee
\be \varphi_1\mbox{ is Lipschitz },\ee
\be \varphi_t(\partial U)\subset\partial U,\ee
and
\be|\varphi_t(x)-x|<\d\mbox{ for all }x\in E\cap\partial U.\ee

Such a $\varphi_1$ is called a sliding deformation in $\bar U$, and $F$ is called a $\d$-sliding deformation of $E$ in $\bar U$.
\end{defn}

Let $\F_\d(E,\bar U)$ denote the set of all $\d$-sliding deformation of $E$ in $\bar U$, and let $\oF_\d(E,\bar U)$ be the family of sets that are Hausdorff limit of sequences in $\F_\d(E,\bar U)$ that do not converge to the boundary. That is: we set
\be \begin{split}\oF_d(E,U)&=\{F\subset\bar U: \exists \{E_n\}_n\subset\F_\d(E,U)\mbox{ such that }d_K(E_n,F)\to 0\\
&\mbox{ for all compact set }K\subset \R^n, 
\tb{(2.1)}\mbox{ holds for }F\mbox{, and }\H^d(F\cap\partial U\bs E)=0\}.\end{split}\ee

\begin{defn}[$\d$-sliding minimal sets]Let $\d>0$, $U\subset \R^n$ be open, and let $E\subset \bar U$ be closed. We say that $E$ is $\d$-sliding minimal in $\bar U$, if \tb{(2.1)} holds, and \tb{(2.2)} holds for all $F\in \oF_\d(E,\bar U)$.
\end{defn}


\begin{defn}[Stable minimal cones] Let $K$ be a 2-dimensional Almgren minimal cone in $\R^n$.

$1^\circ$ We say that $K$ is $(\eta, \d)$-sliding stable, if there exists $\eta\in (0,10^{-2}\eta_0)$, and $\d>0$, such that $K$ is $\d$-sliding minimal in $\bar \cU(K,\eta)$. We say that $K$ is sliding stable if it is $(\eta, \d)$-sliding stable for some $\eta>0$ and $\d>0$.

$2^\circ$ We say that $K$ is $(\eta, \d)$-measure stable, if there exists $\eta\in (0,10^{-2}\eta_0)$, and $\d>0$, such that for all $y\in \R^n$ with $||y||<\d$, we have
\be \H^2(K\cap \cU(K,\eta))=\H^2((K+y)\cap \cU(K,\eta)).\ee
We say that $K$ is measure stable if it is $(\eta, \d)$-measure stable for some $\eta>0$ and $\d>0$.

\end{defn}

\begin{rem}It is not hard to see from the definition that the sliding and measure stabilities does not depend on the ambient dimension. That is, if $K\subset \R^d$ is a 2-dimensional minimal cone, then it is $(\eta,\d)$-sliding (resp. measure) stable if and only if it is $(\eta,\d)$-sliding (resp. measure) stable in $\R^n$ for any $n\ge d$.
\end{rem}

We will see in \cite{stablePYT} and \cite{stableYXY} that up to now, all known 2-dimensional minimal cones are sliding stable (and hence measure stable). In addition, independent of the hypothesis of sliding stability, we can prove that every 2-dimensional minimal cone is measure stable, which is stated in the following theorem:

\begin{thm}[cf. \cite{stablePYT} Theorem \tb{3.1}]Let $K$ be a 2-dimensional minimal cone in $\R^n$. Then $K$ is measure stable.
\end{thm}

\begin{rem} Let $K^1$ and $K^2$ be two $(\eta,\d)$-measure stable minimal cones of dimension 2, with dist$(K^1\cap \pa B,K^2\cap \pa B)>\eta_0$, then even if the cone $C=K^1\cup K^2$ might not be minimal, we can define $\cU(C,\eta)$, and the measure stability in the sense of Definition \tb{2.29}. And directly from the definition, we know that $C$ is also $(\eta,\d)$-measure stable.
\end{rem}

\subsection{Unit simple $d$-vectors and related results}

Let $2\le d<n$ be two integers. Denote by $\wg_d(\R^n)$ the space of all $d$-vectors in $\R^n$. 
Set $I_{n,d}=\{I=(i_1,i_2,\cdots, i_d):1\le i_1<i_2<\cdots<i_d\le n\}$.
Let $\{e_i\}_{1\le i\le n}$ be an orthonormal basis of $\R^n$. For any $I=(i_1,i_2,\cdots, i_d)\in I_{n,d}$, denote by $e_I=e_{i_1}\wg e_{i_2}\wg\cdots\wg e_{i_d}$. Then the set $\{e_I,I\in I_{n,d}\}$ forms a basis of $\wg_d(\R^n)$. The standard scalar product on $\wg_d(\R^n)$ is: for $\xi=\sum_{I\in I_{n,d}}a_I e_I$ and $\zeta=\sum_{I\in I_{n,d}}b_I e_I$,
\be <\xi,\zeta>=\sum_{I\in I_{n,d}}a_Ib_I.\ee

Denote by $|\cdot|$ the norm induced by this scalar product.

Now given a unit simple $d$-vector $\xi$ (that is, $|\xi|=1$, and $\xi$ can be written as the exterior product of $d$ vectors), we can associate it to a $d$-dimensional subspace $P(\xi)\in G(n,d)$, where $G(n,d)$ denotes the set of all $d$-dimensional subspace of $\R^n$:
\be P(\xi)=\{v\in\R^n,v\wedge\xi=0.\}\ee
In other words, if $\xi=x_1\wg x_2\cdots\wg x_d$, $x_1,\cdots x_d$ being orthogonal, then $P(\xi)$ is the $d-$subspace generated by $\{x_i\}_{1\le i\le d}$.

From time to time, when there is no ambiguity, we also write $P=x_1\wg x_2\wg\cdots \wg x_d$, where $P\in G(n,d)$ and $\{x_i\}_{1\le i\le d}$  are $d$ unit vectors such that $P=P(x_1\wg\cdots \wg x_d)$.

Now if $f$ is a linear map from $\R^n$ to $\R^n$, then we denote by $\wg_d f$ (and sometimes by $f$ if there is no ambiguity) the linear map from $\wg_d(\R^n)$ to $\wg_d(\R^n)$ such that
\be \wg_df(x_1\wg x_2\cdots\wg x_d)=f(x_1)\wg f(x_2)\wg\cdots\wg f(x_d).\ee

And if $P\in G(n,d)$, then $P=P(\xi)$ for some unit simple $d$-vector $\xi$ (such a $d$ vector always exists), we define $|f(\cdot)|:G(n,d)\to\R^+\cup\{0\}$ by
\be |f(P)|=|\wg_df(x_1\wg\cdots \wg x_d)|.\ee 

One can easily verify that the value of $|f(P)|$ does not depend on the choice of the unit simple vector $\xi$ that generates $P$. Hence $|f(\cdot)|$ is well defined.

\begin{lem}\label{F1}(c.f. \cite{Mo84} Lemma 5.2)

Let $P,Q$ be two subspaces of $\R^n$ with 
\be dim(P\cap {Q}^\perp)\ge dim\ P-d+2\ee 
Let $\xi$ be a simple unit $d-$vector in 
$\wg_d(\R^n)$. Denote by $p,q$ the orthogonal projections from $\R^n$ onto $P$ and $Q$ respectively. Then the projections of $\xi$ verify
\be|p(\xi)|+|q(\xi)|\le 1.\ee
If moreover
\be dim(P\cap Q)<d-2,\ee 
then
\be |p(\xi)|+|q(\xi)|=1\mbox{ if and only if }\xi\mbox{ belongs to }P\mbox{ or }Q.\ee
\end{lem}

More generally, we have the following variation for projections of a simple $d$-vector to two almost orthogonal subspaces:

\begin{lem}When $P^1$ and $P^2$ are two subspaces in $\R^n$ of dimension $d_1$ and $d_2$, with minimal angle $\a=\a_{P^1,P^2}$, let $p^i$ denote the orthogonal projection to $P^i$, $i=1,2$, then there exists a map for any $d\le \min\{d_1,d_2\}$, and any unit simple $d$-vector $\xi$, we have
\be |p^1(\xi)|+|p^2(\xi)|\le (1+d\cos\a)|\xi|.\ee
\end{lem}

\nd The proof is the same as in \cite{2p} Proposition 2.19.\qed

\begin{lem}
 Let $n>d\ge 2$, and $P,Q$ be two subspaces in $\R^n$, $F\subset \R^n$ be a $d-$rectifiable set. Denote by $p,q$ the orthogonal projections on $P$ and $Q$ respectively. Let $\lambda\ge 0$ be such that for almost all $x\in F$, the approximate tangent plane $T_xF\in G(n,d)$ of $F$ verifies
\be |p(T_xF)|+|q(T_xF)|\le\lambda.\ee
Then
\be \H^d(p(F))+\H^d(q(F))\le\lambda H^d(F).\ee
\end{lem}

\nd

Denote by $f$ the restriction of $p$ to $F$; then $f$ is a Lipschitz function from a $d$-rectifiable set to a $d-$rectifiable subset of $P$. Since $F$ is $d$-rectifiable, for $\H^d-$almost all $x\in F$, $f$ has an approximate differential
\be ap Df(x): T_xF\to P \ee (c.f.\cite{Fe}, Thm 3.2.19). Moreover this differential is such that $||\bigwedge_d ap Df(x)||\le 1$ almost everywhere, because $f$ is $1-$Lipschitz.

Now we can apply the area formula to $f$, (c.f. \cite{Fe} Cor 3.2.20). For all $\H^d|_F$-integrable functions $g\ : \ F\to\bar\R$, we have
\be \int_F (g\circ f)\cdot ||\wedge_d apDf(x)||d\H^d=\int_{f(F)} g(z)N(f,z)d\H^dz,\ee where $ N(f,z)=\s\{f^{-1}(z)\},$
and for $z\in p(F)$ we have $N(f,z)\ge 1$.
Take $g\equiv1$, we get
\be\int_F ||\wedge_d apDf(x)||d\H^d=\int_{p(F)} N(f,z)d\H^dz\ge \int_{p(F)}d\H^d=\H^d(p(F)).\ee

Recall that $p$ is linear, hence its differential is itself. As a result $apDf(x)$ is the restriction of $p$ on the $d$-subspace $T_xF$, which implies that if $\xi$ is a unit simple $d$-vector with $T_xF=P(\xi)$, then
\be ||\wedge_dapDf(x)||=||\wedge_dapDf(x)(\xi)||=|p(T_xF)|\ee
by \tb{(2.21)}.
 Hence by \tb{(2.31)},
\be\int_F |p(T_xF)| d\H^d(x)\ge \H^d(p(F)).\ee
A similar argument gives also:
\be\int_F |q(T_xF)| d\H^d(x)\ge \H^d(q(F)).\ee
Summing \tb{(2.33) and (2.34)} we get
\be\begin{array}{ll} \H^d(p F)+\H^d(q F)&\le\int_F |p T_xF|+|q T_xF|d\H^d(x)\\
&\le\int_F\ \lambda\ d\H^d(x)=\lambda \H^d(F)\end{array}\ee
since $|p T_xF|+|q T_xF|\le \lambda$. \qed

As a corollary of Lemmas \tb{2.33 and 2.35}, we have

\begin{pro}Let $d\ge 2$, and $E_1,E_2$ two Almgren minimal sets of dimension $d$ in $\R^{n_1}$ and $\R^{n_2}$ respectively. Then the orthogonal union  $E_1\cup E_2$ is an Almgren minimal set in 
$\R^{n_1+n_2}$.
\end{pro}

\nd Let $F$ be a deformation of $E_1\cup E_2$ in $\R^{n_1+n_2}$, then there exists $R>0$ and $f$ a Lipschitz deformation in $\R^{n_1+n_2}$ such that  
\be f(B(0,R))\subset B(0,R);f|_{B(0,R)^C}=Id,\mbox{ and }f(E_1\cup E_2)=F.\ee

Denote by $p_i$ the projection on $\R^{n_i}$, $i=1,2$. Then $p_i\circ f(E_i)$ is a deformation of $E_i$ in $B(0,R)\cap \R^{n_i}$, $i=1,2$. By the Almgren minimality of $E_i$, $\H^d(p_i\circ f(E_i))\ge H^d(E_i)$, hence
\be \H^d(p_i(E))=\H^d(p_i\circ f(E_1\cup E_2))\ge \H^d(p_i\circ f(E_i))\ge \H^d (E_i),i=1,2.\ee

Then we apply Lemma \ref{F1}, and the following lemma, we obtain that
\be \H^d(E)\ge \H^d(p_1(E))+\H^d(p_2(E))\ge \H^d(E_1)+\H^d(E_2)=\H^d(E_1\cup E_2),\ee
and the conclusion follows.\qed

 Let $P^1$ and $P^2$ be two mutually orthogonal and complementary subspaces of $\R^n$ (hence dim$P^1+$dim$P^2=n$). Then for any $z\in \R^n\bs\{0\}$, it can be uniquely decomposed into $z=z^1+z^2$, with $z^i\in P^i, i=1,2$. 
 
 The angle $\a=\arccos\frac{|z^1|}{|z|}=\arcsin\frac{|z^2|}{|z|}$ is called the angle of $z$ between $P^1$ and $P^2$. In other words, a non-zero vector $z$ is of angle $\a$ between $P^2$ and $P^2$ if and only if there exists unit vectors $u\in P^1$ and $v\in P^2$ such that $z=|z|(\cos\a u+\sin\a v)$.

\begin{defn}[Angle of a vector between two mutually orthogonal complementary subspaces] The angle $\a=\arccos\frac{|z^1|}{|z|}=\arcsin\frac{|z^2|}{|z|}$ is called the angle of $z$ between $P^1$ and $P^2$. In other words, a non-zero vector $z$ is of angle $\a$ between $P^2$ and $P^2$ if and only if there exists unit vectors $u\in P^1$ and $v\in P^2$ such that $z=|z|(\cos\a u+\sin\a v)$.
\end{defn}

\begin{defn}[Analytic $d$-simple vectors and analytic $d$-subspaces] Let $P^1$ and $P^2$ be two mutually orthogonal and complementary subspaces of $\R^n$. 

$1^\circ$ Let $P$ be a $d$-dimensional subspace of $\R^n$. We say that $P$ is a $d$-analytic subspace between $P^1$ and $P^2$, if there exists $\a\in [0,\frac\pi 2]$, such that for all $z\in P$, the angle of $z$ between $P^1$ and $P^2$ is $\a$;

$2^\circ$ A simple $d$-vector $\xi$ in $\R^n$ is a analytic $d$-simple vector between $P^1$ and $P^2$ (with angle $\a$) if $P(\xi)$ is an analytic $d$-subspace between $P^1$ and $P^2$ (with angle $\a$).
\end{defn}

\begin{lem}Let $P^1$ and $P^2$ be as above. 

A $d$-simple vector $\xi$ is $d$-analytic between $P^1$ and $P^2$ with angle $\a$ if and only if there exists an orthonomal system $\{u_1,\cdots, u_d\}$ of $P^1$, and an orthonormal system $\{v_1,\cdots, v_d\}$ of $P^2$, such that 
\be \xi=\wg_{i=1}^d(\cos\a u_i+\sin\a v_i);\ee
\end{lem}

\nd 

The only if part follows directly by definition. Let us prove the converse: so let $\{u_1,\cdots, u_d\}$ and $\{v_1,\cdots, v_d\}$ be orthonormal systems of $P^1$ and $P^2$ respectively, such that $\xi=\wg_{i=1}^d(\cos\a u_i+\sin\a v_i)$. Obviously, the system $\{\cos\a u_i+\sin\a v_i,1\le i\le d\}$ forms an orthonormal basis of $P(\xi)$. Hence for any $z\in P(\xi)\bs\{0\}$, there exists $z_1,\cdots, z_d\in \R$ such that 
\be z=\sum_{i=1}^d z_i(\cos\a u_i+\sin\a v_i).\ee
This implies
\be z=\cos\a(\sum_{i=1}^d z_iu_i)+\sin\a(\sum_{i=1}^d z_iv_i),\ee
hence $z$ is of angle $\a$ between $P^1$ and $P^2$. This holds for arbitrary $z\in P(\xi)\bs \{0\}$. Hence $P(\xi)$ is $d$-analytic.
 \qed
 
 \begin{defn}[Minimal angles between spaces and cones]

$1^\circ$ For every pair of non trivial subspaces $P^1$ and $P^2$ of $\R^{n_1+n_2}$ of dimension $n_1$ and $n_2$ respectively, set 
\be \a_{P^1,P^2}=\inf\{\theta: |<u,v>|\le \cos\theta ||u||\cdot ||v||,\ \forall u\in P^1,v\in P^2\},\ee
$\a_{P^1,P^2}$ is called the minimal angle between $P^1$ and $P^2$.

$2^\circ$ For every pair of non trivial cones $A^1$ and $A^2$ in $\R^n$, the minimal angle $\a_{A^1,A^2}$ between $A^1$ and $A^2$ is defined similarly :
\be \a_{A^1,A^2}=\inf\{\theta: |<u,v>|\le \cos\theta ||u||\cdot ||v||,\ \forall u\in A^1,v\in A^2\}.\ee
And in fact, if $P^i$ is the subspace generated by $A_i$, $i=1,2$. Then
\be \a_{A_1,A_2}= \a_{P^1,P^2}.\ee

$3^\circ$ We say that two cones are orthogonal, if the minimal angles between them is $\frac\pi 2$.
\end{defn}

\begin{rem}When the two minimal cones $K^i,i=1, 2$ are planes, we use characteristic angles to define relative positions between two subspaces, as in \cite{2ptopo}.
\end{rem}

\section{Uniqueness of the orthogonal union}

In this section we discuss the uniqueness of the orthogonal union.

So let $d\ge 2$, and let $K^1\subset \R^{n_1}$ and $K^2\subset \R^{n_2}$ be two $d$-dimensional Almgren (resp. $G$-topological) minimal sets. We want to prove that the orthogonal union of $K^1$ and $K^2$ will be somehow the only minimal set under some natural assumption, such as spanning the same boundary, etc.

To this aim we will of course make some assumptions on the uniqueness of $K^1$ and $K^2$, otherwise the uniqueness of their orthogonal union fails obviously. So let $C$ be a $d$-dimensional Almgren minimal set in a convex domain $U$ that contains the origin, we say that 

$1^\circ$ $C$ is Almgren unique in $U$ if it is the only set in $\oF(C,U)$ that attains the minimal measure. That is:
\be \forall E\in \overline\F(C,U), \H^d(E)=\inf_{F\in \overline \F(C,U)}\H^d(F)\Rightarrow E=C.\ee

$2^\circ$ $C$ is $G$-topological unique in $U$ if $C$ is $G$-topologically minimal, and 
\be \begin{split}\mbox{For any }d\mbox{-dimensional $G$-topological competitor }E\mbox{ of }C\mbox{ in }U,\\
\H^d(E\cap U)=\H^d(C\cap U), \mbox{ implies }C=E;\end{split}\ee

$3^\circ$ We say that a $d$-dimensional minimal set $C$ in $\R^n$ is Almgren (resp. $G$-topological) unique, if it is Almgren (resp. $G$-topologial) unique in every domain $U\subset \R^n$.

\begin{rem}
$1^\circ$ note that the condition $\H^d(E)=\inf_{F\in \overline\F}\H^d(F)$ already implies that $E$ is itself a minimal set, since the class $\oF$ is stable under deformations. Also notice that $\H^d(E)=\inf_{F\in \overline\F}\H^d(F)$ is equivalent to the condition $\H^d(E)\le\inf_{F\in \overline\F}\H^d(F)$ since $E\in \overline \F$.

$2^\circ$ By Corollary \tb{4.7} of \cite{uniquePYT}, when $E$ is a 2-dimensional minimal cone, $\H^d(C\cap U)=\inf_{F\in \oF(C,U)}\H^2(F)$.

$3^\circ$ If $C$ is an Almgren unique minimal set in $U$, $V\subset U$ is a domain, then $C$ is also Almgren unique minimal in $V$.

$4^\circ$ If $C$ is a minimal cone, then its uniqueness does not depend on domains. Therefore we always say it is unique (without specifying in which domain).

$5^\circ$ So far for almost all of the 2-dimensional minimal cones we know, they are all Almgren unique. See \cite{uniquePYT}. We also remark that for minimal cones with calibrations or paired calibrations, or satisfy some topological condition, they automatically satisfy the sliding stability. Sets with calibrations and paired calibrations are very likely to satisfy the Almgren uniqueness as well. 
\end{rem}

 Due to the lack of information on the structure for minimal sets, especially those of dimension greater or equal to 3, we are kind of still far from a general conclusion for any minimal sets. However, from the very little information we get, here we can still treat some cases.

\begin{thm}\label{unicite}[Almgren uniqueness] Let $K^1\subset \R^{n_1}$ and $K^2\subset \R^{n_2}$ be two 2-dimensional Almgren minimal cones. Denote by $C_0=K^1\cup_\perp K^2$ their orthogonal union, which is Almgren minimal by Proposition \tb{2.36}. Then if $K^i,i=1,2$ are Almgren unique, so is $C_0$.

%
\end{thm}

The rest of this section will be devoted to the proof of Theorem \ref{unicite}. The proof will be composed of a sequence of lemmas and propositions. 

So let $K^1$ and $K^2$ be two 2-dimensional Almgren unique minimal cones in $\oB(0,1)$, and $C_0=K^1\cup_\perp K^2$. 

To save notations, let $P^i=\R^{n_i}, i=1,2$, and denote by $p^i$ the orthogonal projection to $P^i,i=1,2$.
Denote by $B=B(0,1)$, $\F=\F(C_0, B)$, $\oF=\oF(C_0,B)$. Also set $B^i=B\cap P^i$, and for any $x\in P^i$, $B^i(x,r)=B(x,r)\cap P^i, i=1,2$.
Let $E\in \oF$ be such that $\H^2(E)=\inf_{F\in \overline\F}\H^2(F)$. Then $E$ is Almgren minimal in $B$.

Still set $\Xi:=\{\xi\in\wedge_2\R^{n_1+n_2}\mbox{ unit simple },|p^1(\xi)|+|p^2(\xi)|=1\}$. Then we have

\begin{lem}
1) For $\H^d$-almost all $x\in E$, $T_xE\in P(\Xi)$.

2) $p^i(E)=K^i, i=1,2.$
 
3) For $\H^2$ almost all $x\in p^i(E)= K^i$, 
\be \s \{{p^i}^{-1}\{x\}\cap E\}=1.\ee 
\end{lem}

\nd Since $E\in \oF$, for $i=1,2$, the projection $p^i(E)\in \oF(K^i, B^i)$. As a result, we have 
\be \H^2(p^i(E))\ge \inf_{F\in \overline\F(K^i,B^i)}\H^2(F)=\H^2(K^i), i=1,2.\ee

By Lemma \tb{2.35}, with $\lambda=1$, we have
\be \H^2(E)\ge\H^2(p^1(E))+\H^2(p^2(E))\ge \H^2(K^1)+\H^2(K^2)=\H^2(C_0)=\H^2(E),\ee
hence all the inequalities above and in the proof of Lemma \tb{2.35} (with $\lambda=1$) are equalities. In particular, 1) follows, and for $i=1,2$,
\be \H^2(p^i(E))=\H^2(K^i)\mbox{, and for almost all }x\in p^i(E), \s \{{p^i}^{-1}\{x\}\cap E\}=1.\ee
This implies, by Almgren uniqueness of $K^i$, that $p^i(E)=K^i, i=1,2$, and hence 
\be \s \{{p^i}^{-1}\{x\}\cap E\}=1\mbox{ for }a.e. x\in p^i(E)=K^i.\ee\qed

\begin{cor}If $x\in E$ is of type $\P$, then $T_xE\in P(\Xi)$.
\end{cor}

\nd Suppose that $T_xE\not\in P(\Xi)$. By the $C^1$ regularity of $E$ around $x$, there exists a neighborhood $U$ of $x$, such that for all $y\in U\cap E$, $T_yE$ exists, and $T_yE\not\in P(\Xi)$, because $P(\Xi)$ is closed. But $U$ is open, thus by the Ahlfors regularity of minimal sets, $H^2(E\cap U)$ is of positive measure, which contradicts Lemma \tb{3.3} 1).\qed

Now let us decide the structure of $\Xi$.

\begin{lem}Let $x,y$ be two unit vectors in $\R^{n_1+n_2}=P^1\times P^2$, then $x\wedge y\in \Xi$ if and only if they verify the following property.

There exists unit vectors $u_1,u_2,v_1,v_2$, with $u_i\in P^1,v_i\in P^2$, $u_1\perp u_2,v_1\perp v_2$, and $\a\in[0,\frac\pi2]$, such that 
\be x=\cos\a u_1+\sin\a v_1,y=\cos\a u_2+\sin\a v_2.\ee 

In other words, $x\wedge y\in \Xi$ if and only if it is an analytic 2-simple vectors. Hence $P(\Xi)$ is exactly the set of all 2-analytic subspaces.
\end{lem}

\nd We will only prove the ''only if'' part. The converse is trivial. 

For all unit vectors $x,y\in \R^{n_1+n_2}$, we have 
\be \begin{split}|p^1(x\wedge y)|+|p^2(x\wedge y)|&=|p^1(x)\wedge p^1(y)|+|p^2(x)\wedge p^2(y)|\\
&\le|p^1(x)||p^1(y)|+|p^2(x)||p^2(y)|\\
&\le \frac12(|p^1(x)|^2+|p^1(y)|^2+|p^2(x)|^2+|p^2(y)|^2)=1.\end{split}\ee

The last equality is because $\R^{n_1+n_2}=P^1\oplus P^2$, and hence $|p^1(x)|^2+|p^2(x)|^2=|p^1(y)|^2+|p^2(y)|^2=1.$

Then $x\wedge y\in \Xi$ if and only if all inequalities in \tb{(3.9)} are equalities. The first means that $p^1(x)\perp p^1(y)$ and $p^2(x)\perp p^2(y)$; the second shows that $|p^1(x)|=|p^1(y)|$ and $|p^2(x)|=|p^2(y)|$. Denote by $u_1=p^1(x),u_2=p^1(y),v_1=p^2(x),v_2=p^2(y)$, and $\a=\arctan\frac{|p^2(x)|}{|p^1(x)|}$, we obtain the conclusion.\qed

\begin{cor}Let $x,y$ be two unit vectors in $\R^{n_1+n_2}=P^1\times P^2$, with $x\wedge y\in \Xi$. If $x\wg y\not\in P^1$, then $|p^2(x\wg y)|>0$. In other words, the linear map $p^2|_{P(x\wg y)}$ is of full rank (non degenerated).
\end{cor}

\nd By Lemma \tb{3.5}, there exists unit vectors $u_1,u_2,v_1,v_2$, with $u_i\in P^1,v_i\in P^2$, $u_1\perp u_2,v_1\perp v_2$, and $\a\in[0,\frac\pi2]$, such that $ x=\cos\a u_1+\sin\a v_1,y=\cos\a u_2+\sin\a v_2.$ The assumption $x\wg y\not\in P^1$ implies that $\a\ne 0$, hence $\sin\a\ne 0$. Therefore $|p^2(x\wg y)|=\sin^2\a>0.$\qed

\bigskip

Next let us look at tangent cones for $\Y$ points of $E$.

Recall that $E_Y$ denotes the set of all $\Y$ points of $E$. By Corollary \tb{2.23}, $E_Y$ is a locally finite union of $C^1$ curves: $E=\cup_{j\in \N}\gamma_j$, where the endpoints of the $\gamma_j$'s are either in $E_T$, or belong to the boundary $E\cap \pa B=C_0\cap \pa B$. Each $\gamma_j$ will be called a $\Y$ curve in the following text.

\begin{lem}[Position of tangent $\Y$ cones] Let $y\in E_Y$, then the tangent $\Y$ set $C_y=C_yE$ must belong to an analytic 3-subspace between $P^1$ and $P^2$.
\end{lem}

\nd

Denote by $P_i,1\le i\le 3$ the three closed half planes in $\R^n$ such that  $C_y=\cup_{i=1}^3 P_i$. Then the intersection of the three $P_i$ is a line $D$. Denote by $Q_i$ the plane containing $P_i$, $w$ the unit vector that generates $D$, and $w_i$ the unit vector such that  $w_i\perp w$ and $Q_i=P(w\wedge w_i)$. Thus the angle between $w_i$ and $w_j$, $i\ne j$, is $120^\circ$. And the three vectors $w_i,1\le i\le 3$ belongs to a 2-plane $W$.

By $C^1$ regularity for $\P$ points and $\Y$ points in $E$, we know that $Q_i\in P(\Xi), 1\le i\le 3$.

Suppose the angle of $w$ between $P^1$ and $P^2$ is $\a$. Then $Q_i=P(w\wg w_i)\in P(\xi),1\le i\le 3$. By Lemma \tb{3.5}, the angle of $w_i,1\le i\le 3$ between $P^1$ and $P^2$ are all $\a$. Hence there exists normal vectors $u_i,1\le i\le 3$ and $v_i,1\le i\le 3$ such that $u_i\in P^1, v_i\in P^2$, and $w_i=\cos\a u_i+\sin\a v_i$, $1\le i\le 3$.

We have
\be 0=\sum_{i=1}^3 w_i=\cos\a(\sum_{i=1}^3 u_i)+\sin\a(\sum_{i=1}^3 v_i),\ee
therefore
\be 0=\cos\a(\sum_{i=1}^3 u_i)=\sin\a(\sum_{i=1}^3 v_i),\ee
since $\cos\a(\sum_{i=1}^3 u_i)\in P^1$ and $\cos\a(\sum_{i=1}^3 v_i)\in P^2$.

If $\cos\a=0$ or $\sin\a=0$, then the four vectors $w$ and $w_i,1\le i\le 3$ belong to one of the $P^1$ or $P^2$. In particular the 3-subspace generated by them belong to one of the $P^1$ or $P^2$, hence is 3-analytic;

Otherwise, we have $0=\sum_{i=1}^3 u_i=\sum_{i=1}^3 v_i$. As result, $|u_1+u_2|=|u_3|=1$. Note that $|u_1|=|u_2|=1$. Hence the only possibility is that the angle between $u_i$ and $u_j$, $i\ne j$, is $120^\circ$. The same for $v_i, 1\le i\le 3$.

Now let $z=\frac{w_2-w_3}{\sqrt 3}\in W$ (recall that $W$ is the 2-plane generated by $w_i,1\le i\le 3$). Then $w_1$ and $z$ form an orthonormal basis of $W$. We have
\be \begin{split}z&=\frac{1}{\sqrt 3} w_2-\frac{1}{\sqrt 3} w_3=\frac{1}{\sqrt 3}(\cos\a u_2+\sin\a v_2)+\frac{1}{\sqrt 3}(\cos\a u_3+\sin\a v_3)\\
&=\cos\a(\frac{1}{\sqrt 3} u_2-\frac{1}{\sqrt 3} u_3)+\sin\a(\frac{1}{\sqrt 3} v_2-\frac{1}{\sqrt 3} v_3).
\end{split}\ee

Since the angles between $u_2, u_3$ and $v_2,v_3$ are both $120^\circ$, the two vectors $z_1=\frac{1}{\sqrt 3} u_2-\frac{1}{\sqrt 3} u_3\in P^1$ and $z_2=\frac{1}{\sqrt 3} v_2-\frac{1}{\sqrt 3}\in P^2$ are unit vectors, such that $z_1\perp u_1, z_2\perp v_1$. Thus $z=\cos \a z_1+\sin\a z_2\in \Theta$, and the angle of $z$ between $P^1$ and $P^2$ is $\a$.

Next suppose that $w=\cos\a u+\sin\a v$ with $u\in P^1$ and $v\in P^2$ unit vectors. We know that for $i=1,2,3$, $w\wg w_i\in \Xi$.
Since $w\perp w_i$, we have $u\perp u_i,v\perp v_i,1\le i\le 3$. As a result, $u\perp z_1$ and $v\perp z_2$ as well. 

Now we know that the $\Y$ set $C_y$ belong to the 3-subspace $V$ of which $\{w,w_1,z\}$ forms an orthonormal basis, and the angles of $w,w_1$ and $z$ between $P^1$ and $P^2$ are all $\a$. By Lemma \tb{2.39}, $V$ is 3-analytic between $P^1$ and $P^2$.\qed

\begin{rem}In particular, when one of the $P^i,i=1,2$, say, $P^1$, is of dimension 2 (or equivalently, $C^1$ is a plane), then $C_y$ must belong to $P^2$. 
\end{rem}

\begin{lem}[Projection of $E_Y$] Let $y\in E_Y$ such that the blow up limit $C_y\not\in P^2$. Then $p^1(y)\in K^1_S$. Similarly, if $C_y\not\in P^1$ then $p^2(y)\in K^2_S$.\end{lem}

\nd We keep the notation ($D,w, P_i,Q_i,w_i, 1\le i\le 3$) as in the proof of lemma \tb{3.7}. Let $\a$ be the angle of the vectors $w,w_i,1\le i\le 3$ between $P^1$ and $P^2$. Then since $C_y\not\in P^2$, $\a\ne\frac\pi 2$.

By Remark \tb{2.12} $2^\circ$, $p^1(C_y)$ is contained in a blow up limit of $p^1(E)=K^1$ at $p^1(y)$.


Let us look at the projections of $w$ and $w_i,1\le i\le 3$. Since $\a\ne \frac\pi 2$ the components $u\perp u_i,1\le i\le 3$, and $u_i,1\le i\le 3$ make $120^\circ$ between each pair. In other words, the vectors $\{u,u_1,u_2,u_3\}$ generates a 3-subspace in $P^1$. In particular, $p^1(C_y)$ spans a 3-subspace in $P^1$. Hence $p^1(C_y)$ is not a plane. As a result, the set $p^1(E)=K^1$ admit a blow up limit $p^1(C_y)$ at $p^1(y)$ which is not a plane. By Definition \tb{2.17}, $p^1(y)\in K^1_S$.\qed

\begin{lem}Let $x\in E_Y$ such that $p^1(x)\in K^1_P$ is a regular point of $K^1$. Let $j\in \N$ be such that $x\in \gamma_j$, then $\gamma_j$ is parallel to $P^2$. 
\end{lem}

\nd Let $x\in E_Y$ be such that $p^1(x)\in K^1_P$. Let $\gamma_x\subset E_Y$ be a $\Y$ curve that contains $x$. 

By Corollary \tb{2.23} $1^\circ$, the set $K^1_P$ is an open set in $K^1$. Hence there exists $r>0$ such that $B^1(p^1(x),r)\cap K^1\subset K^1_P$. Then for each $y\in \gamma_x$ such that $p^1(y)\in B^1(p^1(x),r)$, we also have $p^1(y)\in K^1_P$. Denote by $\gamma_0=\gamma_x\cap {p^1}^{-1}(B^1(p^1(x),r)$. Then $\gamma_0$ is a $C^1$ curve, and by Lemma \tb{3.9}, for each $y\in \gamma_0$, the blow up limit of $E$ at $y$ must satisfy 
$C_y\in P^2$. In particular, the tangent line to $\gamma_0$ at $y$, which is the spine of $C_y$, lies in $P^2$. As a result, the curve $\gamma_0$ is parallel to $P^2$.

The above argument yields that for each $x\in E_Y$ such that $p^1(x)\in K^1_P$, there exists a neighborhood of $y$ in $\gamma_x$ which is parallel to $P^2$. Therefore the set $\{y\in \gamma_x: T_y\gamma_x$ parallel to $P^2\}$ is open in $\gamma_x$. Obviously it is also closed, hence the whole $\Y$ curve $\gamma_x$ is parallel to $P^2$.\qed

%
%

\begin{pro}$p^1(E_Y)\subset K^1_S$.
\end{pro}

\nd We prove by contradiction. So let $x\in E_Y$ be such that $p^1(x)\in K^1_P$. 

Let $\gamma_x\subset E_Y$ be a $\Y$ curve that contains $x$. By Lemmas \tb{3.9}, for each $y\in \gamma_x$, $C_y\in P^2$. This implies that a blow up limit of $K^2=p^2(E)$ at $p^2(y)$ contains $p^2(C_y)=C_y$, which is a singular set. Hence $p^2(y)\in K^2_S$. As a result, $p^2(\gamma_x)\subset K^2_S$. On the other hand, by Lemma \tb{3.10}, $\gamma_x$ is parallel to $P^2$, in particular, $p^2(\gamma_x)$ is a $C^1$ curve in $K^2$. Hence $p^2(\gamma_x)$ is contained in one of the $\Y$-line $\xi$ in $K^2$. Note that $\xi$ is a segment with the origin as one end point, and the other endpoint belongs to $\pa B^2$.

Let $a,b$ be the endpoints of $\gamma_x$. Since $\gamma_x\subset \xi$, one of the $p^2(a)$ and $p^2(b)$ is further from the origin than the other one. Suppose $p^2(a)$ is further from the origin. In particular, $p^2(a)\in K^2_Y$.

Obviously $a\not\in \pa B$: since $p^1(a)=p^1(x)\in K^1_P$, we know that $p^1(a)\ne 0$. But $E\cap \pa B=(K^1\cup K^2)\cap \pa B\subset {p^1}^{-1}\{0\}\cup {p^1}^{-1}\{0\}$, therefore $a\not\in E\cap \pa B$. But $a\in E$, since $E$ is closed. Hence $a\not\in \pa B$.

As consequence, as an endpoint of a $\Y$ curve $\gamma_x$ of $E$, $a$ must be a $\T$ type singular point of $E$: $a\in E_T$. Thus a blow-up limit $C_a$ of $E$ at $a$ is a $\T$ type cone. By the structure Theorem \tb{2.21} and the bi-H\"older regularity Theorem \tb{2.16} for the minimal cone $C_a$, there exists finitely many $\Y$ curves $\xi_1,\cdots, \xi_l$ of $E$, that meet at $a$, and $E_S\cap B(a,t)=[\{a\}\cup(\cup_{1\le j\le l}\xi_j)]\cap B(a,t)$. Since $p^1(x)\in K^1_P$, and $K^1_P$ is open in $K^1$, there exists $r>0$ such that $B^1(p^1(x), r)\cap K^1\subset K^1_P$. Then for each $1\le j\le l$, $\xi_j\cap {p^1}^{-1}(B^1(p^1(x), r))\subset K^1_P$. By Lemma \tb{3.10}, we have $\xi_j$ is parallel to $P^2$. Since one of the endpoints of $\xi_j$ is $a$, we know that $p^1(\xi_j)=p^1(a)=p^1(x),\forall 1\le j\le l$. As a result, $E_S\cap B(a,r)$ is parallel to $P^2$. As a result, $C_a=p^2(C_a)$, and thus by Remark \tb{2.12} $2^\circ$, $C_a$ is contained in a blow up limit of $K^2=p^2(E)$ at $p^2(a)$. But we know that $p^2(a)\in K^2_Y$, hence $C_a$ is in fact a $\Y$-set (since a $\Y$ contains no minimal cones other than itself). This contradicts the fact that $a\in E_T$. 

The contradiction yields that $p^1(E_Y)\subset K^1_S$.\qed

\begin{cor} We have $p^1(E_S)\subset K^1_S$. That is, if $x\in K^1_P$ is a regular point of $K^1$, then ${p^1}^{-1}(x)\cap E\subset E_P$.
\end{cor}

\nd This follows from the fact that $E_S=\overline{E_Y}$.\qed

\bigskip

Next we will study regular points of $E$.

Now let $R$ be a connected component of $K^1_P$. Then by Corollary \tb{2.22}, $R$ is an open planar sector (can be a disk) of dimension 2. Set $E(R)=E\cap {p^1}^{-1}(R)$. By Corollary \tb{3.12}, $E(R)\subset E_P$. Hence $E(R)$ is a 2-dimensional minimal surface. In particular, it is an analytic submanifold. Moreover, by Corollary \tb{3.4}, we know that for every $x\in E(R)$, the tangent plane $T_xE(R)=T_xE\in P(\Xi)$. 

Denote by $Q^1\subset P^1$  the 2-subspace containing $R$.

Set $A=\{x\in E(R):T_xE\subset P^1\}=\{x\in E(R):T_xE=Q^1\}$. Then $A$ is closed in $E(R)$. 

%

\begin{lem}Suppose that $E(R)\bs A\ne\emptyset$. Then for each connected component $\Omega$ of $E(R)\bs A$, there exists a 2-subspace $Q^2\subset P^2$, such that $p^2(\Omega)\subset Q^2$. 

As a result, $\Omega\subset Q^1\times Q^2$.
\end{lem}

\nd  
Set $E_0=\{x\in E(R)\bs A: p^2(x)=0\}$. Let us first prove that $E_0$ is finite. In fact, suppose that $K^2$ contains $k_2$ planes. Then if there are more than $k_2$ points in $E_0$, then there exists at least one plane in $K^2$ of which a neighborhood at $0$ will be covered more than once by the disks $D_x:x\in E_0$. This contradicts Lemma \tb{3.3} 3).

Hence $E_0$ is finite. As a result, $\Omega\bs E_0$ is connected.

Now let $x$ be any point in $\Omega\bs E_0$. Since $x\not\in A$, $T_xE\not\subset P^1$. Then since $T_xE\in P(\Xi)$, we know by Corollary \tb{3.6} that $p^2|_{T_xE}$ is of full rank. Hence there exists $r=r_x>0$ such that $(E_0\cup A)\cap B(x,r)=\emptyset$, and $p^2$ maps a neighborhood $B(x,r)\cap E(R)$ of $x$ to an open 2-submanifold in $P^2$. On the other hand by Lemma \tb{3.3}, $p^2(B(x,r)\cap E(R))\subset K^2$. As a 2-manifold in $K^2$, obviously $p^2(B(x,r)\cap E(R))$ is an open subset of the connected component of $K^2_P$ that contains $p^2(x)$.
This holds for any $x\in \Omega\bs E_0$.

Let $O$ be a connected component of $K^2_P$ that conains $p^2(B(x_0,r_{x_0})\cap E(R))$ for some $x_0\in \Omega\bs E_0$. By the above property for points in $\Omega\bs E_0$, the subset $\{x\in \Omega\bs E_0: p^2(x)\in O\}$ is both open and closed in $\Omega\bs E_0$. Since $\Omega\bs E_0$ is connected, this implies that $\{x\in \Omega\bs E_0: p^2(x)\in O\}=\Omega\bs E_0$. That is, $p^2(\Omega\bs E_0)\subset O$.

Let $Q^2$ be the 2-plane containing $O$. Then $p^2(x)\in Q_2$, $\forall x\in \Omega\bs E_0$. Since $E_0$ is a finite subset of $\Omega$, $\Omega\subset \overline{\{\Omega\bs E_0\}}$. Hence 
\be p^2(\Omega)\subset p^2(\overline{\{\Omega\bs E_0\}})\subset \overline {p^2(\Omega\bs E_0)}\subset \overline O\subset Q^2.\ee\qed

\begin{pro}Suppose that $E(R)\bs A\ne\emptyset$. Let $\Omega$ be a connected component of of $E(R)\bs A$, and let $Q^2$ be the 2-plane containing $p^2(\Omega)$. Let $\{e_1,e_2=ie_1\}$ be an orthonormal basis of $Q^1$, and $\{e_3,e_4=ie_3\}$ be an orthonormal basis of $Q^2$. Then 

$1^\circ$ For any $x\in \Omega$ there exists $r=r_x>0$ such that in $B(x,r), \Omega$ is (under the given basis)  the graph of an injective complex analytic or anti-analytic function 
$\varphi=\varphi_x : Q^2\to Q^1$. More precisely,
\be \Omega\cap B(x,r)= \mbox{graph}(\varphi)\cap B(x,r).\ee

$2^\circ$ Similarly, for each $x\in \Omega$, near $x$, $E$ is the graph of a complex analytic or anti-analytic function from $Q^1$ to $Q^2$. 

$3^\circ$ $\Omega$ is the graph of a complex conformal map from $p^1(\Omega)$ to $p^2(\Omega)$.
\end{pro}

\nd The proof is similar to that of Lemma 3.4 in \cite{2p}, with minor difference. So here we will only sketch the proof here. The readers may refer to Lemma 3.4 of \cite{2p} for more details.

$1^\circ$ 

Let $x\in \Omega$ be as stated in the proposition. Assume all the hypotheses in the lemma. Since $x\in \Omega\subset E_P$, and $T_xE\in P(\Xi)\bs \{Q^1\}$, we know that $Dp^2$ is non degenerate on $T_xE$. Since $\Omega$ is a smooth manifold, there exists $r_x>0$ such that in $B(x,r_x)$, $p^2\lfloor_\Omega$ is injective. So we can define the smooth map
\be\begin{array}{rll}\varphi=\varphi_x: \omega^2_x&\to & P^1,\\
 x&\mapsto & p^1\circ {p^2}^{-1}(x),
 \end{array}\ee
 where $\omega^2_x=p^2(\Omega\cap B(x,r_x))$ is an open set containing $p^2(x)$ in $Q^2$.
 
Then $\Omega$ coincides with the graph of $\varphi$ in $B(x,r_x)$. Moreover, by Corollary \tb{3.4}, $T_{x'}\Omega\in P(\Xi)\bs \{Q^1\}$, for all $x'\in \Omega\cap B(x,r_x)$. Equivalently, for each $y\in \omega^2_x$, $(Id, D\varphi(y))(Q^2)\in P(\Xi)\bs \{Q^1\}$. By the characterisation of $P(\Xi)$ in Lemma \tb{3.5}, this means that the differential $D\varphi(y)$ is complex analytic or anti-analytic (under the given basis).
 
So we have proved that $\varphi$ is a complex smooth function such that at each point, its differential is either complex analytic, or anti-analytic.

Set $B_1=\{y\in \omega^2_x, \frac{d\varphi}{d z}(y)\ne 0\}=\{y\in\omega^2_x: D\varphi(y)$ is analytic. Then $B_1$ is open, since $\varphi$ is $C^1$. If $B_1=\emptyset$, then $\varphi$ is anti-analytic. Otherwise, $B_1$ is not empty, and set $g=\frac{\partial \varphi}{\partial z}$. Then $g$ is continuous on $\omega^2_x$, and $B_1=\{y\in \omega^2_x: g(y)\ne 0\}$. Moreover, since $\frac{d\varphi}{d z}(y)\ne 0$ on $B_1$, $\frac{d\varphi}{d \bar z}(y)=0$ on $B_1$ (since at each point $\varphi$ is either analytic or anti-analytic), and hence $\varphi$ is analytic on the open set $B_1$, so that its derivative $g$ is analytic, too. Then the conclusion $1^\circ$ of Proposition \tb{3.14} will follow from the following theorem (c.f.\cite{Ru} Thm 12.14) :

\begin{thm}[Rad\'o's theorem] Let $U\subset \C$ be an open domain, and $f$ be a continuous function on $\overline U$. Set $\Omega=\{z\in U:f(z)\ne 0\}$, and suppose that $f$ is holomorphic on $\Omega$. Then $f$ is holomorphic on $U$.
 \end{thm} 
 
In fact, we apply the theorem to $g$, and obtain that $g$ is complex analytic on $\omega^2_x$. But since $B_1\ne\emptyset$, $g\not\equiv 0$. As a result $B_1^C=\{y\in B: g(y)=0\}$ does not have any limit point. Notice that $B_1^C\supset B_2:=\{y\in \omega^2_x,\frac{d\varphi}{d\bar z}(y)\ne 0 \}$, and $B_2$ is open, so it is an open set without limit point. Therefore $B_2=\emptyset$, which means that $\varphi$ is complex analytic on $B$.

Hence $\varphi$ is complex analytic or anti-analytic on $\omega^2_x$, which yields $1^\circ$. 

$2^\circ$. 

We can see easily that for any point $x\in \Omega$ with $T_x\Omega\ne Q^2$ the proof will be the same as that of $1^\circ$. So we only have to prove that for each $x\in \Omega$, $T_x\Omega=Q^2$.

We prove by contradiction. Let $x\in \Omega$ with $T_x\Omega=Q^2$. By the conclusion of $1^\circ$, there exists $r=r_x>0$ such that in $B(x,r), \Omega$ is (under the given basis)  the graph of an injective complex analytic or anti-analytic function 
$\varphi=\varphi_x : Q^2\to Q^1$. Modulo change of basis, we suppose that $\varphi$ is complex analytic. $T_x\Omega=Q^2$ means that $\varphi'(y)=0$, where $y=p^2(x)$. Hence $y$ is a zero of order at least 2 for the map $\psi=\varphi-\varphi(y)$. Thus in a punctured neighborhood $\omega^1$ of $\psi(y)\in Q^1$, each point has at least 2 pre-images. So $\omega^1\subset \{z\in P^1:\sharp\{{p^1}^{-1}\{z\}\cap E\}\ge 2\}$.  But $\omega^2$ is of positive measure because it is open. This contradicts Lemma \tb{3.3} 3).

Hence $T_x\Omega\ne Q^2$ for each $x\in \Omega$. As mentioned before, this yields $2^\circ$. 

$3^\circ$. 

Let us prove that $p^2$ is injective on $\Omega$. Suppose there exists $x_1,x_2\in \Omega$ such that $p^2(x_1)=p^2(x_2)$. Take $r<\min\{r_{x_1},r_{x_2}\}$ such that $B(x_1,r)\cap B(x_2,r)=\emptyset$. Then in $B(x_i,r)$, $\Omega$ coincides with the graph of $\varphi_{x_i}:\omega^2_{x_i}\to Q^1$, $i=1,2$. Therefore $p^2((B(x_i,r)\cap\Omega)$ is an open set containing $p^2(x_i)$, $i=1,2$. Let $\omega=p^2((B(x_1,r)\cap\Omega)\cap p^2((B(x_2,r)\cap\Omega)$. Then for each $y\in \omega$, ${p^2}^{-1}(y)$ contains the two points $\varphi_{x_1}^{-1}(y)\in B(x_1,r)$ and $\varphi_{x_2}^{-1}(y)\in B(x_2,r)$. These two points are distinct, since $B(x_1,r)\cap B(x_2,r)=\emptyset$. As consequence, $\omega\subset \{z\in Q^2:\sharp\{{p^2}^{-1}\{z\}\cap E\}\ge 2\}$. But $\omega$ is open and hence is of positive measure. This contradicts Lemma \tb{3.3} 3).

We have thus proved that $p^2$ is injective on $\Omega$. Therefore, $\Omega$ is the graph of a map $f$ from $p^2(\Omega)$ to $p^1(\Omega)$, with $p^i(\Omega)\subset Q^i$, $k=1,2$. And $f$ coincides with each $\varphi_x$ locally near $x$.

By $1^\circ$, the two sets $\{x\in \Omega: \varphi_x$ is analytic$\}$ and $\{x\in \Omega: \varphi_x\}$ is anti-analytic$\}$ are both open, and their union is $\Omega$. But obviously they are also closed in $\Omega$, hence one of them is empty, and the other is the whole $\Omega$. As a result, $f$ an analytic or anti-analytic map itself. Again by the proof of $2^\circ$, $f'$ is never 0. Hence $f$ is a conformal map from $p^1(\Omega)$ to $p^2(\Omega)$.\qed

\begin{cor}$E(R)=R$.
\end{cor}

\nd We would like to prove that $A=E(R)$. We prove by contradiction. So suppose that $A\subsetneq E(R)$. 

Set $\Gamma=\overline E(R)\cap \pa B^1=\overline R\cap\pa B^1$. Then $\Gamma$ is a non trivial part of the boundary of $R$ in $Q^1$.

Since $E(R)$ is an analytic manifold, $A$ is an analytic subvariety of $E(R)$. Since $A\ne E(R)$, the dimension of $A$ is at most 1.
As a result, there exists $x\in \Gamma$ and $r>0$ such that the simply connected open subset $\Delta=B_{Q^1}(x,r)\cap E(R)$ does not meet $A$. As a result, $\Delta$ is connected in $E(R)\bs A$. Let $\Omega$ be the connected component of $E(R)\bs A$ that contains $\Delta$. They by Proposition \tb{3.14} $3^\circ$, $\Omega$ is the graph of a conformal map $\varphi: p^1(\Omega)\to p^2(\Omega)\subset Q^2$. Then the restriction $\psi=\varphi\lfloor_{p^1(\Delta)}$ is also a conformal map, from $p^1(\Delta)$ to $p^2(\Delta)$. Since $\Delta$ is simply connected, so does $p^1(\Delta)$. Note that the boundary of $p^1(\Delta)$ is piecewise smooth, hence the map $\psi$ can be extended continuously to a homeomorphism $\bar\psi$ from $\overline{p^1(\Delta)}$ to $\overline{p^2(\Delta)}$. On the other hand, we know that $\overline E(R)\cap \pa B^1=\overline R\cap\pa B^1=\Gamma$, hence in particular, $\overline\Delta\cap \pa B^1\subset \Gamma\subset Q^1$. This means that the restriction of $\bar\psi$ on $\overline\Delta\cap \pa B^1$ is constantly 0. But $\overline\Delta\cap \pa B^1$ is a non trivial part of the boundary of $p^1(\Delta)$, this contradicts the fact that $\bar\psi$ is a homeomorphism.

Hence $A=E(R)$. By definition of $A$, it is easy to see that in this case, $E(R)=R$.\qed

\noindent\textbf{Proof of Theorem \ref{unicite}.} By Corollary \tb{3.16}, we know that for every connected component $R$ of $K^1_P$, $E\cap {p^1}^{-1}(R)=E(R)=R$. In particular, $R\subset E$. As a result, $K^1_P\subset E$. Since $\overline{K^1_P}=K^1$, and $E$ is closed, we know that $K^1\subset E$. Similarly, $K^2\subset E$. Therefore $C_0=K^1\cup_\perp K^2\subset E$. But we know that 
\be\H^2(E)=\inf_{F\in \overline\F}\H^2(F)\le \H^2(C_0),\ee
hence $E=C_0$, which yields the conclusion of Theorem \ref{unicite}. \qed

Note that almost the same argument gives that:

\begin{thm}\label{unicite1}[Almgren uniqueness] Let $K^1\subset \R^{n_1}$ and $K^2\subset \R^{n_2}$ be two 2-dimensional Almgren minimal cones. Let $p_1,p_2\in \bar B(0,\frac 12)$. Denote by $C_{p_1,p_2}=(K^1+p_1)\cup_\perp (K^2+p_2)$ a translated orthogonal union. (Then $C_{p_1,p_2}$ is Almgren minimal, by Proposition \tb{2.36})
Then if $K^i,i=1,2$ are Almgren unique, so is $C_{p_1,p_2}$.
\end{thm}

\nd We follow the proof of Theorem \ref{unicite}. The only small difference is that Proposition \tb{3.11} and Corollary \tb{3.12} may no be true, if $p^1(q)$ is a regular point of $K^1$. But since this only results in one point of exception, one can easily verify that this does not make any difference to the arguments that follows (which are essentially properties of analytic maps).\qed

\section{A converging sequence of minimal competitors}

Now we begin to prove Theorem \ref{main}. So let $K^i\subset \R^{n_i}, i=1,2$ be 2-dimensional Almgren unique minimal cones in, as in Theorem \ref{main}. Since they are cones, we will only look at what happens in the unit ball $B=B(0,1)$. So in the following text, we suppose that $K^i,i=1,2$ are minimal cones in $\oB(0,1)\cap \R^{n_i}$.

Let $n=n_1+n_2$. In $\R^n=\R^{n_1}\times \R^{n_2}$, set $B^1(x,r)=B(x,r)\cap \R^{n_1}\times \{0\}$ for $x\in \R^{n_1}, r>0$, and set $B^2(x,r)=B(x,r)\cap \{0\}\times\R^{n_1}$ for $x\in \R^{n_2}, r>0$. Set $B^i=B^i(0,1)$, $i=1,2$. Let $K_0^1\subset \oB^1\times\{0\}$ and $K_0^2\subset \{0\}\times \oB^2$ be copies of $K^1$ and $K^2$. Set $P_0^1=\R^{n_1}\times\{0\}, P_0^2=\{0\}\times \R^{n_2}$. Set $C_0=K_0^1\cup K_0^2$.

We prove by contradication, so suppose the conclusion of Theorem \ref{main} is not true. Therefore, there exists $\theta_k\in (\frac\pi 2-\frac 1k,\frac\pi 2)$ (hence $\lim_{k\to\infty}\theta_k=\frac\pi 2$), and copies $K^1_k$ of $K^1_0$, $K^2_k$ of $K^2_0$ in $\R^n$, with $\a_{K^1_k,K^2_k}\ge\theta_k$, such that $C_k=K^1_k\cup K^2_k$ is not Almgren minimal for each $k$. Let $P^i_k$ be the subspace generated by $K^i_k$, $i=1,2, k\in \N$. Without loss of generality, we can suppose that $K_k^1=K_0^1$ and $P_k^1=P_0^1, \forall k\in \N$, and that $d_H(C_k, C_0)\to 0$ as $k\to\infty$.

Define $A^i=\H^2(K_k^i), i=1,2$, and $A=A_1+A_2=\H^2(C_k),\forall k\in \N$.

\begin{pro}\label{existence}For each $k$, there exists $F_k\in \oF(C_k, B)$, such that

$1^\circ$ $\H^2(F_k)=\inf_{F\in \oF(C_k,B)}\H^2(F)<\H^2(C_k)$;

$2^\circ$  $F_k$ is Almgren minimal in $B$;

$3^\circ$ $F_k$ is contained in the convex hull $\Delta_k$ of $C_k\cap \oB$ (which is also the convex hull of $C_k\cap\pa B$). In particular, $F_k\cap\pa B=C_k\cap\pa B$.

\end{pro}

\nd Fix $k$. 

We apply the upper semi continuity Theorem \tb{4.1} of \cite{uniquePYT}, with $U=\R^n\bs (C_k\cap \pa B)$ and $\mathfrak F=\F(C_k,U)$ and $d=2$, and get the existence of a set $F_k\in \oF(C_k, U)$, which is Almgren minimal in $U$, with $\H^2(F_k)\le \inf_{F\in \F(C_k,U)}\H^2(F)$. By the convex hull property for minimal sets, we know that $F_k$ is contained in the convex hull $\Delta_k$ of $F_k\cap \pa U=C_k\cap \pa B$. Hence $F_k\subset \Delta_k\subset \oB$, $F_k\cap\pa B=C_k\cap\pa B$, and 
\be H^2(F_k)\le \inf_{F\in \F(C_k,U)}\H^2(F)\le \inf_{F\in \F(C_k,B)}\H^2(F)<\H^2(C_k),\ee
the last inequality is because $C_k$ is not minimal in $B$.

 Since $F_k$ is minimal in $U$, it is also minimal in $B$. 

 Note that $B$ is uniformly convex, by Corollary \tb{4.7} of \cite{uniquePYT}, we have 
 \be\H^2(F_k)\le \inf_{F\in \oF(C_k,B)}\H^2(F).\ee
 
 By the exact same argument as that after \tb{(4.8)} of \cite{uniquePYT} , we get $F_k\in \oF(C_k,B)$. \qed
 
 Note that the $F_k$ in the above proposition satisfies that $\H^2(F_k)<A=\H^2(C_k)$, since $C_k$ is not minimal in $B$.

Now since $\overline B$ is compact, we can extract a converging subsequence of $\{F_k\}$, denoted still by $\{F_k\}$ for short. Denote by $F_\infty$ its limit. Then $F_\infty$ is contained in $\cap_n\cup_{k>n} \Delta_k$, so that $F_\infty\cap \partial B\subset (\cap_n\cup_{k>n} \Delta_k)\cap\partial B=C_0\cap\partial B$. On the other hand $F_\infty\cap \partial B\supset\lim_{k\to\infty} (\Delta_k\cap\partial B)=C_0\cap\partial B$. Hence 
\be F_\infty\cap\partial B=C_0\cap\partial B.\ee

\begin{pro} $F_\infty=C_0$.
\end{pro} 

\nd We would like to use the uniqueness theorem \ref{unicite}, to prove that $F_\infty$ is in fact $C_0$. So we have to check all the conditions:

1) First, we know that $K^1_0$ and $K^2_0$ are Almgren unique minimal cones in $B$;

2) Since for each $k$, $F_k$ is Almgren minimal, hence the lower semi continuity of Hausdorff measure for minimal sets (cf.\cite{GD03} Thm 3.4) holds, that is, 
\be \H^2(F_\infty)\le\liminf_{k\to\infty}\H^2(F_k)\le A=\inf_{F\in \oF(C_0,B)}\H^2(F).\ee

3) The rest is to check that $F_\infty\in \oF(C_0, B)$. For each k, there exists a linear homeomorphism $\varphi_k: \R^n\to \R^n$ that maps $K_k^1$ to $K_0^1$ and $K_k^2$ to $K_0^2$. We can also suppose that $\varphi_k(B)\subset B$, and $\varphi_k\to Id$ as $k\to\infty$. Then it is easy to see that $\varphi_k(F_k)\in \oF(C_0,B)$. Since $\oF(C_0,B)$ is closed under Hausdorff limit, $F_\infty\in \oF(C_0,B)$ as well.

Hence $F_\infty$ verifies all the assumptions in Theorem \ref{unicite}. By Theorem \ref{unicite}, $F_\infty=C_0$.\qed

\section{Stopping times argument and $\e$-decomposition}

Recall that in \cite{2p}, when each $K^i$ is a plane, we cut each $F_k$ into two pieces (this is called an $\e$-decomposition). One piece is inside a small ball near the origin, where something complicated happens, and we can only estimate its measure by a projection argument; the other piece is outside the small ball, where $F_k$ is very near the planes. By the regularity of minimal sets near planes, this part of $F_k$ are graphs of the two planes, so that we can get more precise estimates for the measure of this part.

In this section we will do the $\e$-decomposition for the sequences $F_k$. Here the stopping time argument works well, but due to the existence of singularities, we cannot use balls as in \cite{2p}, but to use $\cU(K^i,\eta)$ instead. We will describe the convex domain $\cU(K^i,\eta)$ in more details in the first subsection. In the second subsection, we sketch the idea of the using stopping time arguments, and we will prove the existence of the $\e$-decomposition in the last subsection. 

\subsection{The convex domain $\mU$}

Recall that $K^i,i=1,2$ are 2-dimensional minimal cones in $B$, hence by the structure theorem \tb{2.21},  $K^i\cap \pa B^i$ is a union of circles $\{s^i_j,1\le j\le \mu_i\}$, and arcs of great circles with only $\Y$ type junctions. By Theorem \tb{2.21}, the length of each of these arcs is no less than a certain $\eta_0$. So here we fix some $\eta<10^{-2}\eta_0$, very small, such that $K^i,i=1,2$, $\Y$, and the plane are all $(\eta,\d_0)$-sliding stable and $(\eta,\d_0)$-measure stable for some $\d_0>0$.

Denote by $\{a_j^i, 1\le j\le m_i\}$ the set of $\Y$ points in $K^i\cap \pa B$, $i=1,2$.

For $i=1,2$, set
\be \cU^i=\{x\in B^i: <x, y>< 1-\eta, \forall y\in K^i\mbox{ and }<x,a^i_j><1-2\eta, \forall 1\le j\le m_i\}\subset P^i. \ee 
Then $\cU^i$ is the $\eta$-convex domain $\cU(K^i, \eta)$ for $K^i$ in $P^i$, $i=1,2$.
%

Now for $i=1,2$, $1\le j,l\le m_i$, let $\gamma^i_{jl}$ denote the arc of great circle that connects $a^i_j$ and $a^i_l$, if it exists; otherwise set $\gamma^i_{jl}=\emptyset$. Set $J^i=\{(j,l): 1\le j,l\le m_i\mbox{ and }\gamma^i_{jl}\ne\emptyset\}$, which is exactly the set of pairs $(j,l)$ such that the $\Y$ points $a^i_j$ and $a^i_l$ are connected directly by an arc of great circle on $K^i$. 
 
Denote by $A^i_j$ the $n_i-1$-dimensional planar part centered at $(1-2\eta)a^i_j$ of $\partial \cU^i$. That is,  
\be A^i_j=\{x\in B^1: <x,a_j^i>=1-2\eta\mbox{ and }<x,y>\le 1-\eta, \forall y\in K^i\}.\ee

The shape of this planar region (of dimension $n_i-1$) will be obained by cutting off 3 small planar part of a $n_i-1$-dimensional ball . Take $a^1_i$ for example, suppose, without loss of generality, that the three $\Y$ points in $K^1\cap \pa B^1$ that are adjacent to $a_1^1$ are $a_l^1,l=2,3,4$. Then the planar region centered at $(1-2\eta)a_1^1$ is obtained by:

Firstly, take the $n_i-1$-dimensional ball $\Omega_1^1$ perpendicular to $\overrightarrow{oa_1^1}$ and centered at $(1-2\eta)a_1^1$ (hence the radius of $\Omega_1^1$ is $R=\sqrt{1-(1-2\eta)^2}$). For $l=2,3,4$, denote by $x_l$ the intersection of $(1-\eta)\gamma^1_{1l}$ with $\Omega^1_1$. Then the $x_l,l=2,3,4$ will situated on a 2-plane passing through the center of the ball $\Omega_1^1$, hence they belong to a same great circle. Let $L^1_{1l}$ be the $n_i-1$ subspace containing $x_l$ and orthogonal to $\overrightarrow{(1-2\eta)a^1_1, x_l}$.  We then cut off the small part of  $\Omega_1^1$ that is on the other side of $L^1_{1l}, l=2,3,4$, and get the planar region $A^1_1$ . This forms part of the boundary of $\cU^1$.

Let us look at the boundary of $A_1^1$: By definition, it is composed of three disjoint $n_i-2$-dimensional small balls : 
$$I^1_{1l}=\{x\in B^1: <x,a_1^1>=1-2\eta\mbox{ and }<x,y>=1-\eta\mbox{ for some }y\in \gamma_{1l}^1\},2\le l\le 4$$
and the rest of the boundary of $\Omega^1_1$. Note that the diameter of $I^1_{1l}$ is $R_1=\sqrt{1-(1-\eta)^2}$.

See the figure below.  

\centerline{\includegraphics[width=0.7\textwidth]{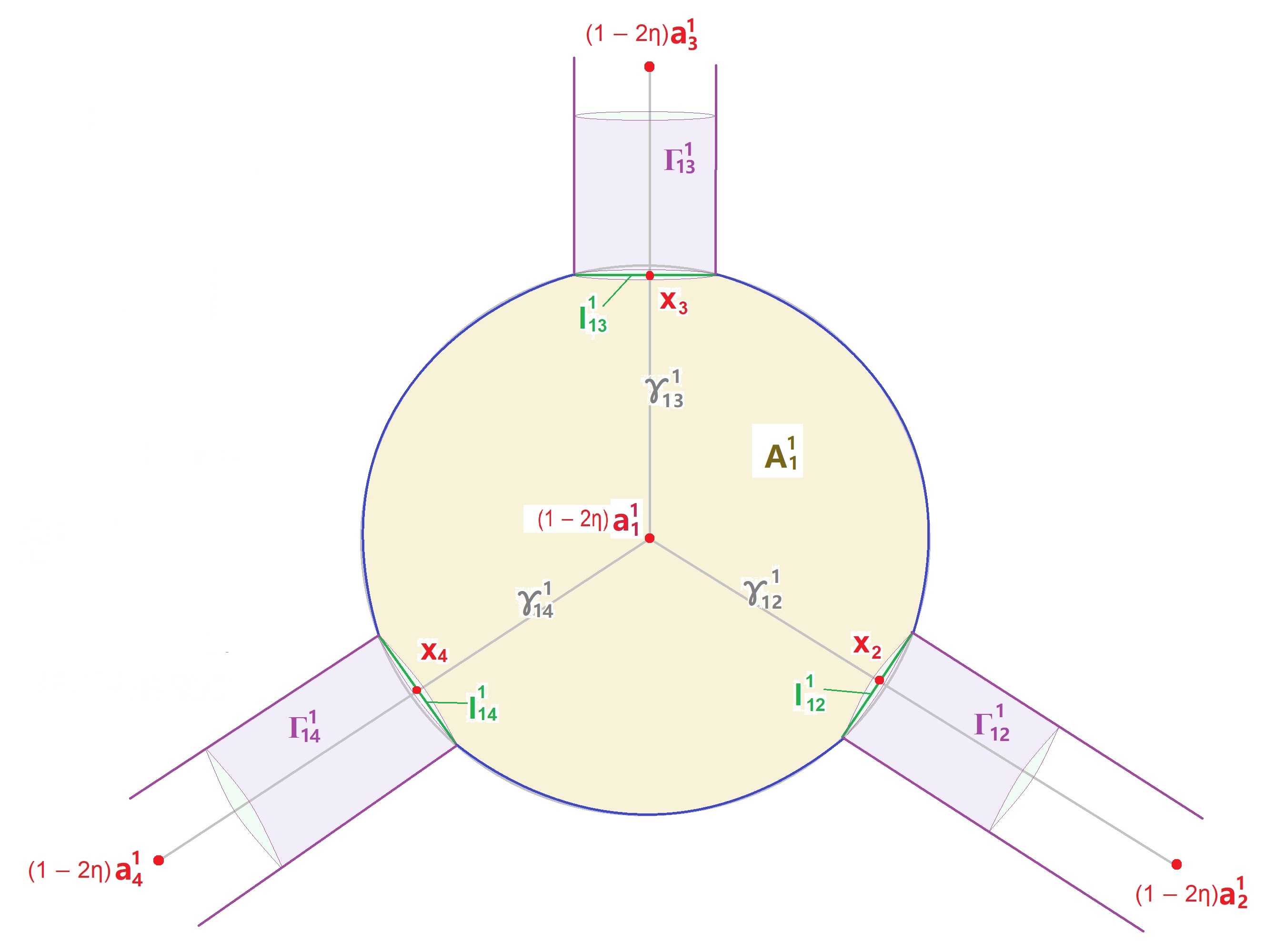} }

We define similarly, for $i=1,2$ and any $(j,l)\in J^i$: $I^i_{jl}$. Then the boundary of $A^i_j, j\in m_i$ is the union of the balls $I^i_{jl},(j,l)\in J^i$, and the rest of the sphere $\pa\Omega_1^1$.

Set $A^i=\cup_{1\le j\le m_i}A_j^i$. This is the whole planar part of $\pa \cU^i$. 
 and will be used to get estimates of the measures of $F_k$ near the singular set of $K^i$.

For the rest of $\pa \cU^i$ that is not spherical, they are obtained by the equation
\be x\in B^i, <x,y>\le 1-\eta, \forall y\in K^i.\ee

So set 
\be\Gamma^i_{jl}=\{x\in B^i, <x,y>=1-\eta\mbox{ for some }y\in \gamma^i_{jl}\}\bs \mA^i,\ee
with $\mA^i$ being the cone over $A^i$ centered at 0,
and
\be S^i_j=\{x\in B^i,<x,y>=1-\eta\mbox{ for some }y\in s^i_j\}.\ee
Then $\Gamma^i_{jl}$ is the band like part of $\pa U^i$ near each $(1-\eta)\gamma^i_{jl}$, and similar for $S^i_j$. The union
$\Gamma^i=\cup_{1\le j,l\le m_i}\Gamma^i_{jl}$ together with $S^i=\cup_{1\le j\le \mu_i}S^i_j$ is the whole cylinderical part of $\pa \cU^i$, and will be used to get estimates of the measures of $F_k$ near the regular set of $K^i$.

Set $\mathfrak c^i_{jl}$, $\mC^i_{jl}$, $\mC^i$, $\ms^i_j$, $\mS^i_j$, $\mS^i$, $\mA^i_j$ and $\mA^i$ the cone (centered at 0) over $\gamma^i_{jl}$, $\Gamma^i_{jl}$, $\Gamma^i$, $s^i_j$, $S^i_j$, $S^i$, $A^i_j$ and $A^i$ respectively, where for any set $S\subset \R^n$, the cone over $S$ is defined to be $\{ts: x\in S, t\ge 0\}$.

Next, for any subset $S\subset \R^n$, and any $x\in \R^n$, $r>0$, set $S_{x,r}=x+rS$.

Last, we define the decomposition region for each $k\in \N$. Recall that $P^i_k$ is the $n_i$-dimensional subspace of $\R^n$ that contains $K^i_k$, $i=1,2$. Denote by $p^i_k$ the orthogonal projection from $\R^n$ to $P^i_k$ for $i=1,2$.

For $i=1,2$, let $\psi_k^i: \R^{n_i}\to P_k^i$ be the isometry that maps $K^i\to K^i_k$. Set $\cU_k^i(0,r)=\psi_k^i(\cU^i(0,r)), i=1,2$. Denote by
$a^i_{j,k}$, $\gamma^i_{jl,k}$, $\Gamma^i_{jl,k}$, $\Gamma^i_k$, $s^i_{j,k}$, $S^i_{j,k}$, $S^i_k$, $A^i_{j,k}$, $I^i_{jl,k}$, $\mI^i_{jl,k}$, $A^i_k$, $\mathfrak c^i_{jl,k}$, $\mC^i_{jl,k}$, $\mC^i_k$, $\ms^i_{j,k}$, $\mS^i_{j,k}$, $\mS^i_kl$, $\mA^i_{j,k}$ and $\mA^i_k$ the image of $a^i_j$, $\gamma^i_{jl}$, $\Gamma^i_{jl}$, $\Gamma^i$, $A^i_j$, $I^i_{jl}$,$\mI^i_{jl}$ ,$A^i$, $\mathfrak c^i_{jl}$, $\mC^i_{jl}$, $\mC^i$, $\mA^i_j$ and $\mA^i$ under $\psi_k^i$ respectively. These are subsets of $P_k^i$.

Define the cylinders
\be \mathfrak U_k^i(x,r)={p_k^i}^{-1}(\cU_k^i(0,r))+x, \forall x\in \R^n, r>0, i=1,2,\ee
and set
Set
\be\mU_k(x,r)=\mU_k^1(x,r)\cap \mU_k^2(x,r)=\{y\in \R^n: p_k^i(y-x)\in \cU_k^i(0,r), i=1,2\}.\ee

Note that $\mU_k(x,r)$ is convex, and that $B(x,\frac12 r)\subset \mU_k(x,r)\subset B(x,2r)$ when $k$ is large.

Note that all these sets depend on $\eta$.

In what follows, when $k$ is fixed, we always write $\mU$ instead of $\mU_k$, to save notations. 

\subsection{The stopping time argument}

The construction of the stopping time procedure is similar to the one that we used in \cite{2p}. So here we will sketch the precedure. The readers may refer to \cite{2p} for more details. On the other hand, the proof of the fact that the procedure has to stop at a finite step (Proposition \tb{5.4}) is different.

\bigskip

For any $\e>0$, we say that two sets $E,F$ are $\e r$ near each other in an open set $U$ if 
\be d_{r,U}(E,F)<\e,\ee
where 
\be d_{r,U}(E,F)=\frac 1r\max\{\sup\{d(y,F):y\in E\cap U\},\sup\{d(y,E):y\in F\cap U\}\}.\ee


We set also
\be \begin{split}&d^k_{x,r}(E,F)=d_{r,\mU_k(x,r)}(E,F)\\
     &=\frac 1r\max\{\sup\{d(y,F):y\in E\cap \mU_k(x,r)\},\sup\{d(y,E):y\in F\cap \mU_k(x,r)\}\}.\end{split}\ee

\begin{rem}We should be clear about the fact that 
\be d_{r,U}(E,F)\ne\frac 1rd_H(E\cap U,F\cap  U).\ee
To see this, we can take $U=\mU_k(x,r)$, and set $E_n=\partial \mU_k(x,r+\frac 1n)$ and $F_n=\partial \mU_k(x,r-\frac 1n)$. Then we have 
\be d^k_{x,r}(E_n,F_n)\to 0\ee and \be \frac 1rd_H(E_n\cap \mU_k(x,r),F_n\cap \mU_k(x,r))=\frac 1rd_H(E_n\cap \mU_k(x,r),\emptyset)=\infty.\ee
So $d_{r,U}$ measures rather how the part of one set in the open set $U$ could be approximated by the other set, and vice versa. However we always have
\be d^k_{x,r}(E,F)\le\frac 1rd_H(E\cap \mU_k(x,r),F\cap \mU_k(x,r)).\ee
\end{rem}

\bigskip

Now we start our stopping time argument. Recall that $\{F_k\}$ is a sequence of sets as in Proposition \ref{existence},  with $\theta_k>\frac\pi2-\frac1k$, which converges to $C_0\cap\oB$. 

We fix a small $\e<\min\{10^{-5}, 10^{-3}\eta\}$ and a large $k$, and we set $s_i=2^{-i}$ for $i\ge 0$. Denote by $\mU(x,r)=\mU_k(x,r), d_{x,r}=d^k_{x,r}$ for short. Then we proceed as follows.

Step 1: Denote by $q_0=q_1=0$, then in $\mU(q_0,s_0)$, the set $F_k$ is $\e s_0$ near $C_k+q_1$ when $k$ is large, because $F_k\to C_0$ and $C_k\to C_0$ implies that $d_{q_0,s_0}(F_k,C_k+q_1)=d_{0,1}(F_k,C_k)\to 0$.

Step 2: If in $\mU(q_1,s_1)$, there is no point $q\in\R^n$ such that $F_k$ is $\e s_1$ near $C_k+q$, we stop here; otherwise, there exists a point $q_2$ such that $F_k$ is $\e s_1$ near $C_k+q_2$ in $\mU(q_1,s_1)$. Here since $\e$ is small ($\e<10^{-5}$), such a $q_2$ is automatically in $\mU(q_1,\frac12 s_1)$, by the conclusion of the step 1. Then in $\mU(q_1,s_1)$ we have simultaneously  
\be d_{q_1,s_1}(F_k,C_k+q_1)\le s_1^{-1}d_{q_0,s_0}(F_k,C_k+q_1)\le 2\e;\ d_{q_1,s_1}(F_k, C_k+q_2)\le\e.\ee
This implies (by definition of $d_{x,r}$) that $d_{q_1,\frac12 s_1}(C_k+q_1,C_k+q_2)\le 12\e$ when $\e$ is small. And hence $d(q_1,q_2)\le 6\e$.

Now we are going to define our iteration process. Notice that this process depends on $\e$, hence we also call it a $\e-$process.

Suppose that $\{q_i\}$ are defined for all $i\le m$, with 
\be d(q_1,q_{i+1})\le 12 s_i\e=12\times 2^{-i}\e\ee
for $0\le i\le m-1$, and hence
\be d(q_i,q_j)\le 24\e_{\min(i,j)}=2^{-\min(i,j)}\times 24\e\ee
for $0\le i,j\le m$, and that for all $i\le m-1$, $F_k$ is $\e s_i$ near $C_k+q_{i+1}$ in $\mU(q_i,s_i)$. We say in this case that the process does not stop at step $m$. Then

Step $m+1$: We look inside $\mU(q_m,s_m)$.

If $F_k$ is not $\e s_m$ near any $C_k+q$ in this ``ball'' of radius $s_m$, we stop. In this case, since $d(q_{m-1},q_m)\le 12\e s_{m-1}$, we have $\mU(q_m,2s_m(1-12\e))=\mU(q_n,s_{m-1}(1-12\e))\subset \mU(q_{m-1},s_{m-1})$, and hence
\be\begin{split}
    d_{q_m,2s_m(1-12\e)}(C_k+q_m,F_k)&\le (1-12\e)^{-1}d_{q_{m-1},s_{m-1}}(C_k+q,F_k)\\
&\le\frac{\e}{1-12\e}.
   \end{split}
\ee
Moreover
\be d(q_m,0)=d(q_m,q_1)\le 2^{-\min(1,m)}\times 24\e=12\e.\ee

Otherwise, we can find a $q_{m+1}\in\R^n$ such that $F_k$ is still $\e s_m$ near $C_k+q_{m+1}$ in $\mU(q_m,s_m)$, then since $\e$ is small, as before we have $d(q_{m+1},q_m)\le 12\e s_m$, and for $i\le m-1$,
\be d(q_i,q_m)\le\sum_{j=i}^n 12\times 2^{-j}\e\le 2^{-\min(i,m)}\times 24\e.\ee
Thus we get our $q_{m+1}$, and say that the process does not stop at step $m+1$.

From the $\e$-process, when $\e$ is small enough, and $k$ is large enough so that $P^1_k$ and $P^2_k$ are almost orthogonal, we have directly the following:

\begin{lem}There exists $\e_0\in (0,\min\{10^{-5}, 10^{-3}\eta\})$, such that for all $\e<\e_0$, $k$ large enough, and for every $m$ such that the $\e-$process does not stop before $m$ (which means in particular that there exists $q_m\in B(q_{m-1},\frac12s_{m-1})$ such that $F_k$ is $\e s_{m-1}$ near $C_k+q_m$ in $\mU(q_{m-1},s_{m-1})$), we have

$1^\circ$  $F_k\cap (\mU(0,1)\backslash \mU(q_m,s_{m+2}))=F_{k,m}^1\cup F_{k,m}^2$, 
where $F_{k,m}^1, F_{k,m}^1$ are $\frac12 s_{m+2}$ far from each other, that is
\be \mbox{dist}(F_{k,m}^1,F_{k,m}^2)\ge \frac12 s_{m+2};\ee

$2^\circ$ Set $G_{k,m}^i=F_{k,m}^i\cap \mU(q_{m-1},s_{m-1})\bs \mU(q_m,s_{m+2})$, then
 $d_{q_m, s_m}(G_{k,m}^i, K_k^i+q_m))<\e$, $i=1,2$.
\end{lem}

\subsection{Existence of the $\e$-decomposition}

Now let us prove that the $\e$-process has to stop at a finite step.

\begin{pro}For each $\e<\e_0$, and $k$ large, if the $\e$ process does not stop at step $m+1$, then there exists a deformations $\varphi^1_m$ in $B^1\times \R^{n_2}=B^1\times P_0^2$ (which contains $B$), satisfying
\be \varphi^1_m(F_k)\subset [F_{k,m}^1]\cup [(\mU(q_m,s_m)\bs \mU(q_{m+1},s_{m+1}))\cap B(P_k^2+q_{m+1}, 100\e s_{m+1})]\cup [\mU(q_{m+1},s_{m+1})]\ee
and
\be \varphi^1_m=id\mbox{ on }[B^1\times P_0^2]^C\cup [F_{k,m}^1]\cup \mU(q_m,s_m).\ee

Similarly, 
there exists a deformations $\varphi^2_m$ in $\R^{n_1}\times B^2=P_0^1\times B^2$ (which contains $B$), satisfying
\be \varphi^2_m(F_k)\subset [F_{k,m}^2]\cup [(\mU(q_m,s_m)\bs \mU(q_{m+1},s_{m+1}))\cap B(P_k^1+q_{m+1}, 100\e s_{m+1})]\cup [\mU(q_{m+1},s_{m+1})]\ee
and
\be \varphi^2_m=id\mbox{ on }[P_0^1\times B^2]^C\cup [F_{k,m}^2]\cup \mU(q_m,s_m).\ee

\end{pro}

\nd Fix any $k$ large. Suppose the $\e$ process does not stop at step $m+1$. We will only prove the existence of $\varphi^1_m$, since the argument for $\varphi^2_m$ is exactly the same.

We will construct, by recurrence, a sequence of deformations in $B^1\times \R^{n_2}=B^1\times P_0^2$ (which contains $B$), that maps $F_k$ to $[F_{k,m}^1]\cup [(\mU(q_m,s_m)\bs \mU(q_m,s_{m+1}))\cap B(P_k^2+q_m, 100\e s_m)]\cup [\mU(q_m,s_{m+1})]$.

For each $j\le m$, let $\pi_j: \R^n\to \overline{\mU(q_j,s_j)}$ denote the shortest distance projection to the convex set $\overline{\mU(q_j,s_j)}$, and set $h_j: B(P_k^2+q_j, 100s_j\e)\to B(P_k^2+q_{j+1}, 100s_{j+1}\e)$ the orthogonal projection from the $100s_j\e$-neightborhood of $P_k^2+q_j$ to the $100s_{j+1}\e$-neightborhood of $P_k^2+q_{j+1}$. Note that since dist$(q_j,q_{j+1})<24s_{j+1}\e$, $B(P_k^2+q_{j+1}, 100s_{j+1}\e)\subset B(P_k^2+q_j, 100s_j\e)$.

Let us first map $F_k$ to $[F_{k,1}^1]\cup [(\mU(q_1,s_1)\bs \mU(q_1,s_2))\cap B(P_k^2+q_2, 100\e s_2)]\cup [\mU(q_1,s_2)]$. In fact, we decompose $F_k$ two three disjoint parts $F_{k,1}^1\bs \mU(q_2,s_2)$, $F_{k,1}^2\bs \mU(q_2,s_2)$, and $\mU(q_2,s_2)$, and we will define the deformation on these three parts : set
\be g_1(x)=\left\{\begin{array}{rcl}h_1\circ\pi_1(x)&\ if\ & x\in F_{k,1}^2\bs \mU(q_2,s_2);\\
 x&\ if\ & x\in F_{k,1}^1\bs \mU(q_2,s_2);\\ x&\ if\ & x\in \mU(q_2,s_2).\end{array}\right.\ee
 
So in fact we have only moved points in $F_{k,1}^2\bs \mU(q_2,s_2)$, and clearly $g_1(F_k)\subset [F_{k,1}^1]\cup [(\mU(q_1,s_1)\bs \mU(q_2,s_2))\cap B(P_k^2+q_2, 100\e s_2)]\cup [\mU(q_2,s_2)]$. 

Let us check that $g_1$ is Lipschitz. Note that by definition, $g_1$ is Lipschitz on $ [F_{k,1}^1]\cup [\mU(q_2,s_2)]$, and on $F_{k,1}^2\bs \mU(q_2,s_2)$ respectively, so we only to deal with points that are close to both parts. Apply Lemma \tb{5.2} to step $m$, we know that the parts $F_{k,1}^1$ and $F_{k,1}^2$ are far from each other; on the other hand, since the $\e$-process does not stop at step $m+1$, we know that $F_{k,1}^2\cap \mU(q_1,s_1)$ is already contained in $\mU(q_1,s_1)\cap B(P_k^2+q_2,100\e s_2)$, hence $h_1\circ\pi_1=id$ on $F_{k,1}^2\cap \mU(q_1,s_1)$, which is a neighborhood of the compact set $\overline {F_{k,1}^2\bs \mU(q_2,s_2)}\cap \overline {[F_{k,1}^1]\cup [\mU(q_2,s_2)]}$. 

As a result, $g_1$ is a Lipschitz map on $F_k$. Moreover, it only moves points in $F_{k,1}^2\bs \mU(q_1,s_1)$, which is a compact set in $B^1\times P_0^2$. So we can extend $g_1'$ to a Lipschitz map to $\R^n$, with 
\be g_1'=id\mbox{ on }[B^1\times P_0^2]^C\cup [F_{k,1}^1]\cup \mU(q_1,s_1),\ee
\be g_1'(F^2_{k,1}\bs \mU(q_2,s_2))\subset (\mU(q_1,s_1)\bs\mU(q_2,s_2))\cap B(P_k^2+q_2, 100\e s_2),\ee
and 
\be g_1'(B^1\times P_0^2)\subset B^1\times P_0^2.\ee

Thus $g_1'$ is a Lipschitz deformation in $B^1\times P_0^2$, and 
\be g_1'(F_k)\subset [F_{k,1}^1]\cup [(\mU(q_1,s_1)\bs \mU(q_2,s_2))\cap B(P_k^2+q_2, 100\e s_2)]\cup [F_k\cap \mU(q_2,s_2)].\ee


Set $\tilde g_1=g_1'$.

Let us next define, for each $j\le m$: $\tilde g_j: F_k\to [F_{k,j}^1]\cup [(\mU(q_j,s_j)\bs \mU(q_{j+1},s_{j+1}))\cap B(P_k^2+q_{j+1}, 100\e s_{j+1})]\cup [\mU(q_{j+1},s_{j+1})]$ by recurrence. So suppose that the Lipschitz deformation map 
\be\begin{split}\tilde g_{j-1}:F_k&\to [F_{k,j-1}^1\bs \mU(q_j,s_j)]\cup [(\mU(q_{j-1},s_{j-1})\bs \mU(q_j,s_j))\cap B(P_k^2+q_j,100 \e s_j)]\cup [F_k\cap\mU(q_j,s_j)]\\
&=[F_{k,j}^1\bs \mU(q_j,s_j)]\cup [(\mU(q_{j-1},s_{j-1})\bs \mU(q_j,s_j))\cap B(P_k^2+q_j,100 \e s_j)]\cup [F_k\cap\mU(q_j,s_j)]\end{split}\ee
 is already defined, that satisfies 
 
 \be \tilde g_{j-1}=id\mbox{ on }[B^1\times P_0^2]^C\cup [F_{k,j-1}^1]\cup \mU(q_{j-1},s_{j-1}),\ee
\be \tilde g_{j-1}(F^2_{k,j-1}\bs \mU(q_j,s_j))\subset (\mU(q_{j-1},s_{j-1})\bs \mU(q_j,s_j)\cap B(P_k^2+q_j, 100\e s_j),\ee
and 
\be g_{j-1}(B^1\times P_0^2)\subset B^1\times P_0^2.\ee
%
 
Define $g_j': \tilde g_{j-1}(F_k)\to [F_{k,j}^1\bs \mU(q_{j+1},s_{j+1})]\cup [(\mU(q_j,s_j)\bs \mU(q_{j+1},s_{j+1}))\cap B(P_k^2+q_{j+1}, 100\e s_{j+1})]\cup [\mU(q_{j+1},s_{j+1})]$: 

\be g_j'(x)=\left\{\begin{array}{rcl}h_j\circ\pi_j(x)&\ if\ & x\in (\mU(q_{j-1},s_{j-1})\bs \mU(q_j,s_j))\cap B(P_k^2+q_j, 100\e s_j);\\
 x&\ if\ & x\in F_{k,j}^1;\\ x&\ if\ & x\in \mU(q_{j+1},s_{j+1}).\end{array}\right.\ee
 
 Then a similar discuss gives that $g_j$ is the identity map on $F_{k,j}^1\cup \mU(q_j,s_j)$, and 
 \be \begin{split}&g_j'\circ \tilde g_{j-1}(F^2_{k,j}\bs \mU(q_{j+1},s_{j+1}))\\
 &\subset g_j'\circ\tilde g_{j-1}[F^2_{k,j}\bs \mU(q_j,s_j)]\cup g_j'\circ\tilde g_{j-1}[F^2_{k,j}\cap \mU(q_j,s_j)\bs \mU(q_{j+1},s_{j+1})]\\
 &\subset g_j'[\mU(q_{j-1},s_{j-1})\bs \mU(q_j,s_j))\cap B(P_k^2+q_j, 100\e s_j)]\cup g_j'(F^2_{k,j}\cap \mU(q_j,s_j)\bs \mU(q_{j+1},s_{j+1}))\\
 &\subset \{[\mU(q_j,s_j)\bs \mU(q_{j+1},s_{j+1})]\cap B(P_k^2+q_{j+1}, 100\e s_{j+1})\}\cup (F^2_{k,j}\cap \mU(q_j,s_j)\bs \mU(q_{j+1},s_{j+1}))\\
 &\subset [\mU(q_j,s_j)\bs \mU(q_{j+1},s_{j+1})]\cap B(P_k^2+q_{j+1}, 100\e s_{j+1})
 \end{split}\ee
 
 The second inclusion is due to \tb{(5.32) and (5.33)}; the third inclusion is due to the definition of $g_j'$, and the last one is because $F_{k,j}^2$ is $s_{j+1}\e$ close to $P_k^2$ in $\mU(q_j,s_j)\bs \mU(q_{j+1},s_{j+1})$.
 
 Similar to $g_1'$, we know that $g_j'$ is Lipschitz, and we can extend it to a global Lipschitz deformation $g_j$  in $B^1\times P_0^2$. Set $\tilde g_j=g_j\circ \tilde g_{j-1}$, then it verifies 

\be \tilde g_j=id\mbox{ on }[B^1\times P_0^2]^C\cup [F_{k,j}^1]\cup \mU(q_j,s_j),\ee
\be \tilde g_j(F^2_{k,j}\bs \mU(q_{j+1},s_{j+1}))\subset (\mU(q_j,s_j)\bs \mU(q_{j+1},s_{j+1})\cap B(P_k^2+q_{j+1}, 100\e s_{j+1}),\ee
and 
\be g_j(B^1\times P_0^2)\subset B^1\times P_0^2.\ee

Thus by recurrence, we get the deformation $\tilde g_m$, and set $\varphi_m^1=\tilde g_m$.\qed

As a direct corollary, we have

\begin{pro}For each $\e<\e_0$, $k$ large enough, there exists $m=m_k$ such that the $\e$-process does not stop at steps before $m$, but stops at step $m$.
\end{pro}

\nd We prove by contradiction. So suppose that for some $\e<\e_0$, and a fixed $k$, the $\e$-process does not stop at any finite step. Then by Proposition \tb{5.3}, for each $m\in \N$, there exists a deformation $\varphi^1_m$ such that \tb{(5.22) and (5.23)} hold. Thus we know that
\be F^1_{k,m}\subset \varphi^1_m(F_k)\subset F^1_{k,m}\cup \mU(q_m,s_m).\ee

We take the Hausdorff limit, and get
\be \lim_{m\to\infty} \varphi^1_m(F_k)\stackrel{d_H}{\to}\cup_{m\in \N}F^1_{k,m}.\ee

Similarly, we have
\be  \lim_{m\to\infty} \varphi^2_m(F_k)\stackrel{d_H}{\to}\cup_{m\in \N}F^2_{k,m}.\ee

Set $F^i=\cup_{m\in \N}F^i_{k,m}$, $i=1,2$, then we can see from Lemma \tb{5.2} that 
$ F^1\cap F^2=\emptyset.$ And from \tb{(5.41) and (5.42)}, we know that $F^1\in\oF(F_k, B^1\times P_k^2)$ and $F^2\in\oF(F_k, P_k^1\times B^2)$.

Since $F_k$ is a Hausdorff limit of deformations of $C_k$ in $B$, it contains a Hausdorff limit of deformation of $K_k^1$ in $B$, which is also a Hausdorff limit of deformations of $K_k^1$ in $B^1\times P_k^2$. As a result, since $F^1\in\oF(F_k, B^1\times P_k^2)$, it contains a Hausdorff limit of deformations of $K_k^1$ in $B^1\times P_k^2$. Therefore

\be \H^2(F^1)\ge\inf_{F\in \oF(K_k^1,B^1\times P_k^2)}\H^2(F).\ee

On the other hand, note that $K_k^1$ is minimal in $B^1$, for any $F_0\in \oF(K_k^1,B^1\times P_k^2)$ we have $p_k^1(F_0)\in \oF(K_k^1,B^1)$, hence
\be \H^2(F_0)\ge \H^2(p_k^1(F_0))\ge\inf_{F\in \oF(K_k^1,B^1)}\H^2(F)=\H^2(K_k^1),\ee
where the last equation is due to upper semi continuous Theorem \tb{4.1} of \cite{uniquePYT}.

As a result, we have
\be \H^2(F^1)\ge\inf_{F_0\in \oF(K_k^1,B^1\times P_k^2)}\H^2(F_0)\ge \H^2(K_k^1).\ee

The same argument gives
\be \H^2(F^2)\ge \H^2(K_k^2).\ee

But then, since $F^1$ and $F^2$ are disjoint subsets of $F_k$, we have
\be \H^2(F_k)\ge \H^2(F^1)+\H^2(F_2)\ge \H^2(K_k^1)+\H^2(K_k^2)=\H^2(C_k),\ee
which contradicts our hypothesis that $\H^2(F_k)<\H^2(C_k)$. This contradiction yields that the $\e$-process has to stop at a finite step.\qed

\begin{rem}In fact, the arguments in this and the previous subsections hold for arbitrary dimension $d\ge 2$.
\end{rem}

For each $k$, if the process stops at step $m=m_k$, we define $o_k=q_m,r_k=s_m$. Then $\mU_k:=\mU_k(o_k,r_k)$ is the critical ball that we look for, because inside the small ball, by definition we know that $F_k$ is $\e s_m$ far from any translation of $C_k$, but outside it, things are near. We also have, by \tb{(5.17)}, $d(o_k,0)\le 12\e$, hence the center $o_k$ of the critical ball is near the origin.


\section{Structure of $F_k$}

In this section we will give two structure theorems for $F_k$ inside and outside $\mU_k(q_m,s_m)$, when the $\e$-process does not stop at step $m$. The results of this section will be in the next sections to estimate the measure of $F_k$.

\subsection{Structure of $F_k\bs \mU_k$ in regular and singular regions}

In this subsection we will analyse the geometry and regularity for $F_k$ outside $\mU_k(q_m,s_m)$ when the $\e$-process does not stop at step $m$.

The main result of this section will be a theorem on local geometric structure of $F_k\bs\mU_k(q_m,s_m)$: in regions far from singularities of $C_k+q_m$, things are flat and regular, while near the $\Y$ points of $C_k$, $F_k$ satisfy some topological ''separation'' condition in higher dimension. With the help of this structure theorem, we give estimates of measures of $F_k$ near and far from singularities in the next section.

A big part of the results of this section will highly rely on topology. So let us first give some notations. 

First, for any subset $A$ of $\R^n$ that admit a locally finite triangulation, and any integer $d$, $H_d(A)$ stands for the $d$-th simplicial homology group on $A$ with coefficient in $\bf{\Z_2}$. For any simplicial $d$-cycle $\sigma$ in $A$, $[\sigma]$ denotes the corresponding element in $H_d(A)$. We write $\sigma\sim \sigma'$ if the two $d$-cycles $\sigma$ and $\sigma'$ are homologic. When $A$ is an open subset of $\R^n$, for any smooth closed $d$-surface $\sigma\subset A$, $\sigma$ also stands for the simplicial cycle in the usual sense. (Note that we are using $\Z_2$, hence there is no orientation.)

Let $\tilde C_k$, $\t K_k^1$, and $\t K_k^2$ denote the cone over $C_k$, $K_k^1$ and $K_k^2$, which are in fact the minimal cones in $\R^n$ that coincide with $C_k$, $K_k^1$ and $K_k^2$ respectively in $B(0,1)$. For $i=1,2$, $(j,l)\in J^i$, let $\xi^i_{jl,k}$ denote the midpoint of $\gamma^i_{jl,k}$. Let $Q^i_{jl,k}$ be the $n-2$-linear plane orthogonal to $\mc^i_{jl,k}$ (or equivalently, orthogoonal to the $2$-plane containing $\gamma^i_{jl,k}$). Set $\sigma^i_{jl,k}=\{x\in Q^i_{jl,k}+\xi^i_{jl,k}: ||x-\xi^i_{jl,k}||=\eta\}$, where $\eta$ is the one we used to define the convex domain $\mU$ in the beginning of Section 5.1. Then $\sigma^i_{jl,k}$ is a $n-3$ sphere that links the
 planar part $\mc^i_{jl,k}$. By structure of 2-dimensional minimal cones Theorem \tb{2.21}, we know that 
 
 \be [\sigma^i_{jl,k}]\ne 0\mbox{ and for different triples }(i,j,l), [\sigma^i_{jl,k}]\mbox{ are different }\ee 
in $H_{n-3}(\R^n\bs \t C_k)$.
 
 Set $\t F_k=F_k\cup (\t C_k\bs B(0,1))$. Then $\t F_k\in \oF(\t C_k, B)$. 
 
%
%
%
  
We can easily see that $\sigma^i_{jl,k}\cap (C_k\cup B)=\emptyset$. By Proposition \tb{2.7}, one can see that \tb{(6.1)} holds in $H_{n-3}(\R^n\bs F)$ for any $F\in \oF(C_k,B)$. In particular, \tb{(6.1)} holds in $H_{n-3}(\R^n\bs \t F_k)$.

\begin{pro}There exists $\e_0\in (0,\min\{10^{-5}, 10^{-3}\eta\})$, such that for all $\e<\e_0$, $k$ large enough, and for every $m$ such that the $\e-$process does not stop before $m$, we have

$1^\circ$ For any $t\in [s_{m+2}, s_{m-1}]$, the $n-3$-sphere $t\sigma^i_{jl,k}+q_m$ is homologic to $\sigma^i_{jl,k}$ in $H_{n-3}(\R^n\bs \t F_k)$; 

$2^\circ$ For each $i=1,2$, and $(j,l)\in J^i,i=1,2$, let $\pi^i_{jl,k}$ denote the orthogonal projection from 
$(\mc^i_{jl,k}\bs \mU(q_m,s_m))\times Q^i_{jl,k}\times {P^i_k}^\perp$ to $\mc^i_{jl,k}+q_m$. Then for each  $(j,l)\in J^i$, we have
\be (q_m+(\mc^i_{jl,k}\bs\mA^i_k))\cap \mU(q_{m-1},s_{m-1}\bs\mU(q_m,s_{m+2})=\pi^i_{jl,k}(G^i_{k,m}\cap (\mC^i_{jl,k}\times {P_k^i}^\perp+q_m));\ee

$3^\circ$ For each $i=1,2$, and $1\le j,l \le m_i$, let $P$ be the 2-plane containing $q_m+\mc^1_{jl,k}$. Then $G_{k,m}^i\cap (\mC^1_{jl,k}\times {P_k^i}^\perp+q_m)$ is a graph of a map $h=h^1_{jl,k}: (q_m+(\mc^i_{jl,k}\bs\mA^i_k))\cap \mU(q_{m-1},s_{m-1}\bs\mU(q_m,s_{m+2})\to P^\perp$, with $|\nabla h|<1$;

$4^\circ$ For each $i=1,2$, and $1\le j\le \mu_i$, let $P$ be the 2-plane containing $q_m+\ms^1_{jl,k}$. Then $G_{k,m}^i\cap ((\mS^1_{jl,k}\times {P_k^i}^\perp)+q_m)$ is a graph of a map $h=h^i_{j,k}: (q_m+\ms^i_{jl,k})\cap \mU(q_{m-1},s_{m-1}\bs\mU(q_m,s_{m+2})\to P^\perp$, with $|\nabla h|<1$;

$5^\circ$ For each $i=1,2$, $1\le j\le m_i$, and for any $t\in [s_{m+2},s_{m-1}]$, the intersection of $G^i_{k,m}$ with the $n-1$-dimensional planar part $t(A_{j,k}^i\times {P_k^i}^\perp)+q_m$ is a set that connects the three flat parts $(tI^i_{jl,k}\times {P_k^i}^\perp)+q_m$, $(l,j)\in J^i$, of the boundary of $t(A_{j,k}^i\times {P_k^i}^\perp)+q_m$. That is, $\{G^i_{k,m}\cap[t(A_{j,k}^i\times {P_k^i}^\perp)+q_m]\}\cup [\cup_{(j,l)\in J^i}((tI^i_{jl,k}\times {P_k^i}^\perp)+q_m)]$ contains a connected part that contains the three parts $tI^i_{jl,k}\times {P_k^i}^\perp)+q_m, (j,l)\in J^i$.
\end{pro}

\nd $1^\circ$ First, by definition of $\sigma^i_{jl,k}$, we know that the distance dist$(\sigma^i_{jl,k}, \t C_k)\ge \eta>10^3\e$, hence for all $t>0$ we have
\be dist(t\sigma^i_{jl,k}+q_m, \t C_k+q_m)\ge t\eta>10^3t\e,\ee 
since $C_k$ is a cone.

We prove by induction. Suppose that the $\e$-process does not stop at step $m$. We will prove that for each $m'\le m$,  $1^\circ$ holds. 

$\bf{m'=1:} $ We know that the $\e-$process does not stop at step 1. Hence $F_k+q_1\subset B(C_k+q_1, \e )$, where for each set $E\subset \R^n$ and $r>0$, $B(E, r)$ denotes the $r$ neighborhood of $E$: $B(E,r):=\{x\in \R^n: dist (x,E)<r\}$. Therefore $\t F_k+q_1=(\t C_k\bs B)\cup F_k\subset (\t C_k\bs B)\cup B(C_k+q_1, \e)\subset B(\t C_k+q_1, \e)$. \tb{(6.3)} tells us that for each $t\in (s_0, s_3)$, 
\be dist(t\sigma^i_{jl,k}+q_1, \t F_k+q_1)\ge dist (t\sigma^i_{jl,k}+q_1, \t C_k+q_1)-\e\ge 10^3t\e-\e>0.\ee
Therefore, for any $t_0\in [s_3,s_0]$, the tube $T_{t_0}:=\{t\sigma^i_{jl,k}+q_1: t_0\le t\le s_0\}$ does not intersect $\t F_k+q_1$. Since, for any $t_0\in
(s_0, s_3)$,  $\partial T_{t_0}=(\sigma^i_{jl,k}+q_1)\cup (t_0\sigma^i_{jl,k}+q_1)$, we have $[t_0\sigma^i_{jl,k}+q_1]=[\sigma^i_{jl,k}+q_1]$ in $H_{n-3}(\R^n\bs \t F_k)$, that is, for each $t\in (s_3, s_0)$, 
\be [t\sigma^i_{jl,k}+q_1]=[s_0 \sigma^i_{jl,k}+q_1]=[s_0 \sigma^i_{jl,k}]=[\sigma^i_{jl,k}]\mbox{ in }H_{n-3}(\R^n\bs \t F_k),\ee
because $q_1=0$ and $s_0=1$.

\textbf{Recurrence:} Now suppose that $1^\circ$ holds for some all $1,2,\cdots, m'-1$, where $m'\le m$, and we shall prove it for $m'$. Since the $\e-$process does not stop at step $m'$, a similar argument as above yields that for each $t\in (s_{m'+2}, s_{m'-1})$, 
\be [t\sigma^i_{jl,k}+q_{m'}]=[s_{m'-1}\sigma^i_{jl,k}+q_{m'}]\mbox{ in }H_{n-3}(\R^n\bs \t F_k);\ee 

On the other hand, by hypothesis of recurrence, we know that $[\sigma^i_{jl,k}]=[s_{m'-1}\sigma^i_{jl,k}+q_{m'-1}]$ in $H_{n-3}(\R^n\bs \t F_k)$. Hence we only have to prove that 
\be [s_{m'-1}\sigma^i_{jl,k}+q_{m'-1}]=[s_{m'-1}\sigma^i_{jl,k}+q_{m'}]\mbox{ in }H_{n-3}(\R^n\bs \t F_k).\ee

 By \tb{(5.17)}, we know that $d(q_{m'-2}, q_{m'})\le 24s_{m'-2}\e$, and $d(q_{m'-2}, q_{m'-1})\le 24s_{m'-2}\e$,
hence both $s_{m'-1}\sigma^i_{jl,k}+q_{m'-1}$ and $s_{m'-1}\sigma^i_{jl,k}+q_{m'}$ belong to $\mU(q_{m'-2}, s_{m'-1})$. Since the $\e$-process does not stop at step $m'-1$, we know that $\t F_k\cap \mU(q_{m'-2}, s_{m'-1})\subset B(\t C_k+q_{m'-2},\e s_{m'-2})$. In the mean time, by \tb{(6.3)}, 
\be dist(s_{m'-1}\sigma^i_{jl,k}+q_{m'-2}, \t C_k+q_{m'-2})\ge t\e>10^3s_{m'-1}\e,\ee 
hence
\be \begin{split}&dist(s_{m'-1}\sigma^i_{jl,k}+q_{m'-1}, \t C_k+q_{m'-2})\\
\ge &dist(s_{m'-1}\sigma^i_{jl,k}+q_{m'-2}, \t C_k+q_{m'-2})-d(q_{m'-2}, q_{m'-1})\\
\ge &10^3s_{m'-1}\e-24 s_{m'-2}\e>10^2s_{m'-2}\e,\end{split}\ee
and similarly
\be dist(s_{m'-1}\sigma^i_{jl,k}+q_{m'}, \t C_k+q_{m'-2})>10^2s_{m'-2}\e.\ee
 
Now for $t\in [0,1]$, set $R_t=\{t(x+q_{m'})+(1-t)(x+q_{m'-1}): x\in s_{m'-1}\sigma^i_{jl,k}, t\in [0,1]\}$. Then the straight line homotopy $R=\cup_{t\in [0,1]} R_t$ is a 2-surface included in $\mU(q_{m'-2}, s_{m'-1})$, and  $\partial R=(s_{m'-1}\sigma^i_{jl,k}+q_{m'})\cup (s_{m'-1}\sigma^i_{jl,k}+q_{m'-1})$. On the other hand, for each $x\in s_{m'-1}\sigma^i_{jl,k}$, the segment $I_t:=\{t(x+q_{m'})+(1-t)(x+q_{m'-1}): t\in [0,1]\}$ does not intersect $B(\t C_k+q_{m'-2},\e s_{m'-2})$. In fact, the length of $I_t$ is $d(q_{m'},q_{m'-1})\le 24s_{m'-2}\e$, hence it is included in $B(x+q_{m'}, 24s_{m'-2}\e)$. But by \tb{(6.10)}, $d(x+q_{m'},C_k+q_{m'-2})>10^2s_{m'-2}\e$, hence 
\be dist(I_t, C_k+q_{m'-2})\ge d(x+q_{m'},C_k+q_{m'-2})-24s_{m'-2}\e>50s_{m'-2}\e,\ee
which implies that $I_t$ does not intersect $B(\t C_k+q_{m'-2},\e s_{m'-2})$. As a result, $R=\cup_{t\in [0,1]}I_t$ does not intersect $B(\t C_k+q_{m'-2},\e s_{m'-2})$. Since $R\subset \mU(q_{m'-2}, s_{m'-1})$, and $\t F_k\cap \mU(q_{m'-2}, s_{m'-1})\subset B(\t C_k+q_{m'-2},\e s_{m'-2})$, we know that $R\cap \t F_k=\emptyset$. Since $\partial R=(s_{m'-1}\sigma^i_{jl,k}+q_{m'})\cup (s_{m'-1}\sigma^i_{jl,k}+q_{m'-1})$, we get \tb{(6.7)}. Thus $1^\circ$ holds for $m'$.

By recurrence, $1^\circ$ holds for all steps $m$ where the $\e$-process does not stop.

\bigskip

\noindent $2^\circ$ Without loss of generality we only prove it for $i=1$. Fix $j,l$, and $k$.

Take any $p\in (q_m+(\mc^1_{jl,k}\bs \mA^1_k))\cap \mU(q_{m-1},s_{m-1})\bs\mU(q_m,s_{m+2}).$

Let $Q$ denote the $Q^1_{jl,k}$, for short. and let $t_p=||p-q_m||\in (\frac 12 s_{m+2}, 2s_{m-1})$. Set $D=\{x\in p+Q: ||x-p||\le  \eta t_p$ an $n-2$-dimensional disk centered at $p$. Then $\partial D=t_p\sigma^1_{jl,k}+(p-t_p\xi^1_{jl,k})$. Note that the center  $t_p\xi^1_{jl,k}+q_m$ of  $t_p\sigma^1_{jl,k}+q_m$, and the center $p$ of $D$, both lie in $\mc^1_{jl,k}+q_m$, and are of same distance to $q_m$, hence they are connected by an arc of circle $\gamma$ of center $q_m$ and radius $t_p$. Since $q_m+(\mc^i_{jl,k}\bs \mA^1_k)$ is a sector centered at $q_m$, the whole $\gamma$ lies in it. Let 
$T$ be the tube $T=\cup_{y\in \gamma}\{x \in y+Q: ||x-y||\le  \eta t_p\}$. Then for any $x\in T$, dist$(x, \t C_k)=$dist$(x,q_m+(\mc^1_{jl,k}\bs \mA^1_k))=\eta t_p$, and the boundary of $T$ is $(t_p\sigma^1_{jl,k}+q_m)\cup\partial D$. Since in $\mU(q_{m-1},s_{m-1})\bs\mU(q_m,s_{m+2})$, $\t F_k$ is contained in the $\e s_{m-1}$ neighborhood of $\t C_k+q_m$, we know that $T\cap \t F_k=\emptyset$. As a result, $[t_p\sigma^1_{jl,k}+q_m]=\partial[D]$ in $H_{n-3}(\R^n\bs \t F_k)$. The conclusion of $1^\circ$ tells that $[t_p\sigma^1_{jl,k}+q_m]\ne 0$ in $H_{n-3}(\R^n\bs \t F_k)$. Hence $\partial[D]\ne 0$ in $H_{n-3}(\R^n\bs \t F_k)$. In particular, $D\cap \t F_k\ne\emptyset$. 

On the other hand, notice that $D\subset \mU(q_{m-1},s_{m-1}\bs\mU(q_m,s_{m+2})$, hence $D\cap \t F_k\cap \mU(q_{m-1},s_{m-1})\bs\mU(q_m,s_{m+2})\ne\emptyset$. By definition of $G^i_{k,m}$, the intersection $\t F_k\cap \mU(q_{m-1},s_{m-1})\bs\mU(q_m,s_{m+2})=G^1_{k,m}\cup G^2_{k,m}$. We know that $G^2_{k,m}$ is contained in the $\e s_m$-neighborhood of $K_k^2+q_m$. In particular, the distance between $G^2_{k,m}$ and $q_m+(\mc^1_{jl,k}\bs \mA^i_k))$ is larger than $\e s_{m-1}$. On the other hand, 
\be D\subset B(p,\eta t_p)\subset B(p, 10^{-1}s_m)\subset B((q_m+(\mc^1_{jl,k}\bs \mA^1_k)), 10^{-1}s_{m-1}),\ee
hence $D\cap G^2_{k,m}=\emptyset$. As a result, $D\cap G^1_{k,m}\ne\emptyset$. That is, 
\be p\in \pi^1_{jl,k}(G^1_{k,m}).\ee

On the other hand, since $p\in (q_m+(\mc^1_{jl,k}\bs \mA^1_k))$, $p\in {\pi^1_{jl,k}}^{-1}(q_m+(\mc^1_{jl,k}\bs \mA^1_k))=(q_m+(\mc^1_{jl,k}\bs \mA^1_k))\times Q\times {P_k^1}^\perp$. Thus 
\be p\in \pi^1_{jl,k}[G^1_{k,m}\cap ((q_m+(\mc^1_{jl,k}\bs \mA^1_k))\times Q\times {P_k^1}^\perp)]
.\ee

Since \tb{(6.14)} holds for all $p\in (q_m+(\mc^1_{jl,k}\bs \mA^1_k))\cap \mU(q_{m-1},s_{m-1})\bs\mU(q_m,s_{m+2})$, we know that
\be (q_m+(\mc^1_{jl,k}\bs \mA^1_k))\cap \mU(q_{m-1},s_{m-1})\bs\mU(q_m,s_{m+2})\subset \pi^1_{jl,k}[G^1_{k,m}\cap ((q_m+(\mc^1_{jl,k}\bs \mA^1_k))\times Q\times {P_k^1}^\perp)].\ee

Now we know that $G^1_{k,m}$ is contained in the $\e s_m$-neighborhood of $K_k^1+q_m$, therefore $G^1_{k,m}\cap ((q_m+(\mc^1_{jl,k}\bs \mA^1_k))\times Q\times {P_k^1}^\perp)\subset G^1_{k,m}\cap ((q_m+\mC^1_{jl,k})\times {P_k^1}^\perp)$. Altogether we have
\be (q_m+(\mc^1_{jl,k}\bs\mA^1_k))\cap \mU(q_{m-1},s_{m-1}\bs\mU(q_m,s_{m+2})\subset\pi^1_{jl,k}(G^1_{k,m}\cap ((q_m+\mC^1_{jl,k})\times {P_k^1}^\perp));\ee 

 The inverse inclusion in (6.2) follows directly from the definition of $G^1_{k,m}, \mC^i_{jl,k}$, etc. Hence equality in (6.2) holds for $i=1$. 

\bigskip

\noindent $3^\circ$ Again we only prove for $i=1$. The proof for $i=2$ is exactly the same.

Notice that in the proof of $1^\circ$ and $2^\circ$, we only used the fact that $\e<10^{-3}\eta$. Here we will make a further constraint on $\e_0$: we ask that $\e_0$ be small than the $10^{-3}\e_1(n,2)\eta$, $\e_1$ being the one in Theorem \tb{2.25}.

Set $U=(\mC^1_{jl,k}\times {P_k^1}^\perp)\cap (\mU(q_{m-1},s_{m-1})\bs \mU(q_m,s_{m+2}))$, and $G=G^1_{k,m}\cap \mC^1_{jl,k}\times {P_k^1}^\perp$. Then $G=F_k\cap U$. Let $U_1$ be the $\frac 12\eta s_{m+2}$-neighborhood of $U$.

Since $F_k$ is $\e s_{m-2}$ near $C_k+q_{m-1}$ in $\mU(q_{m-2},s_{m-2})$, by \tb{(5.17)} $F_k$ is also $50 \e s_{m-2}$ near $C_k+q_m$ in $\mU(q_{m-2},s_{m-2})$, and hence in $U_1$. But in $U_1$, $C_k+q_m$ coincides with the 2-plane $P$ which contains the planar part $q_m+\mc^1_{jl,k}$, therefore $F_k$ is $50 \e s_{m-2}$ near $P$ in $U_1$. 

Now take any $x\in G$. Set $r=\frac 12\eta s_{m+2}$. Then since $G\in U$, the ball $B(x,r)\subset U_1$. Hence in $B(x,r)$, $G$ is $50 \e s_{m-2}$ near $P$. That is, 
$$d_{x,r}(G, P)<\frac{50 \e s_{m-2}}{r}=400\frac{\e}{\eta} <\e_1(n,2).$$ 
Therefore by Theorem \tb{2.25}, $G$ coincides with the graph of a $C^1$ map $f_x:P\to P^\perp$ in $B(x, \frac 34 r)$. Moreover $||\nabla f_x||<1$.

Now $G$ is locally $C^1$ graphs on $P$ in $B(x,r)$ for every $x\in G$. We still have to show that $G$ is an entire graph on $(q_m+(\mc^i_{jl,k}\bs\mA^i_k))\cap \mU(q_{m-1},s_{m-1})\bs\mU(q_m,s_{m+2})$. Let $\pi$ denote the orthogonal projection onto $P$. Note that by $2^\circ$, $\pi$ coincides with $\pi^1_{jl,k}$ on $G$, hence we know that $\pi|_G: G\to (q_m+(\mc^i_{jl,k}\bs\mA^i_k))\cap \mU(q_{m-1},s_{m-1})\bs\mU(q_m,s_{m+2})$ is surjective. To prove $G$ is a graph, we only have to show that $\pi$ is injective on $G$.

Suppose it is not, then there exists two points $x,y\in G$ with $\pi(x)=\pi(y)$. Since in $B(x,r)$, $G$ coincides with a $C^1$ graph on $P$, hence $\pi$ is injective in $B(x,r)$. Hence $y\not\in B(x,r)$, that is, $|x-y|>r$. On the other hand, since $G$ is $50 \e s_{m-2}$ near $P$ in $U$, we know that dist$(x,P)<50 \e s_{m-2}$, and $dist(y,P)<50 \e s_{m-2}$. That is, $|\pi(x)-x|<50 \e s_{m-2}$, and $|\pi(y)-y|<50 \e s_{m-2}$. Since $\pi(x)=\pi(y)$, this means that $|x-y|<50 \e s_{m-2}<\e_1{n,2}r<r$, contradiction. Hence $\pi$ is injective on $G$. And we get $3^\circ$ for $i=1$.

\bigskip

\noindent $4^\circ$ The proof of $4^\circ$ is exactly the same as $3^\circ$.

\bigskip

\noindent $5^\circ$ Fix $i,j$ and $k$. Fix $t\in [s_{m+2}, s_{m-1}]$.

Let $\Omega$ denote the $n_i-1$ dimensional region $tA^i_{j,k}+q_m$, whose center is $o_t:=q_m+t(1-2\eta)a^i_{j,k}$. Let $p_1,p_2,p_3$ denote the three intersection points of $(\t C_k+q_m)\cap (\Omega\times {P^i_k}^\perp)$.
Then $p_1,p_2,p_3$ are exactly the intersection of $\partial tA^i_{j,k}+q_m$ with the three faces $t\mc^i_{jl,k}+q_m: (j,l)\in J^i$. Denote by $Q_1$ the $n-1$ dimensional subspace in $\R^n$ orthogonal to $\overrightarrow {o_t p_1}$ and passing through $p_1$, and let $D_1$ be the $n-2$ dimensional ball in $Q_1$, centered at $p_1$ and of radius $tR$, where $R=\sqrt{1-(1-2\eta)^2}$. Define $Q_2,Q_3,D_2,D_3$ similarly. Then the three $n_i-2$-dimensional disks $D_\mu\cap (P_k^i+q_m),\mu=1,2,3$ forms the three planar part of $\partial (tA^i_{j,k}+q_m)$, (they are in fact the three scaled and translated version $tI^i_{jl,k}+q_m$ of $I^i_{jl,k}:(j,l)\in J^i$) and the rest of $\partial tA^i_{j,k}+q_m$ coincides with the $n_i-1$-dimensional ball $\partial t\Omega^i_{j,k}+q_m$. Note that $o_t$ is also the center of $t\Omega^i_{j,k}+q_m$, and $p_1,p_2,p_3$ belong to a 2-dimensional plane passthing through $o_t$, hence they lie on a same great circle of the ball $t\Omega^i_{j,k}+q_m$.

Let $S_i,1\le i\le 3$ denote the segment with endpoints $o$ and $p_i$, and let $Y$ be the $\Y$ shape set $\cup_{i=1}^3 S_i$. Denote by $Q$ the subspace ${P_k^i}^\perp$. 

Then $Y=(\Omega\times Q)\cap (\t C^k+q_m)$. Since $\t C^k+q_m$ is a cone centered at $q_m$, and $G^i_{k,m}$ is $\e s_{m-2}$ near $\t C^k+q_m$ in $\mU(q_{m-2},s_{m-2})$, the intersection $G^i_{k,m}\cap (\Omega\times Q)$ is contained in the $\e s_{m-2}$-neighborhood of $Y$.

To prove that $G^i_{k,m}\cap (\Omega\times Q)$ connects the three connected parts $t(I^i_{jl,k}\times {P_k^i}^\perp)+q_m$, $(j,l)\in J^i$ it is enough to prove that $G^i_{k,m}\cap (\Omega\times Q)$ connectes the three subsets $D_\mu$, $1\le \mu\le 3$, since they are subsets of the three connected sets $t(I^i_{jl,k}\times {P_k^i}^\perp)+q_m$, $(j,l)\in J^i$.

By definition, we would like to show that $(G^i_{k,m}\cap (\Omega\times Q))\cup (\cup_{\mu=1}^3 D_\mu)$ contains a connected set that contains $\cup_{\mu=1}^3 D_\mu$. For each $s\in [0,1]$, let $sD_\mu$ denote the $n-2$ dimensional ball contained in $D_\mu$ of the same center and $s$ times the radius of $D_\mu$.
Since $G^i_{k,m}\cap (\Omega\times Q)$ is contained in the $\e s_{m-2}$-neighborhood $B(Y,\e s_{m-2})$ of $Y$, and $D_\mu, \mu=1,2,3$ are all connected, it is enough to prove that $(G^i_{k,m}\cap \Omega\times Q)\cup (\cup_{\mu=1}^3 \overline{\frac12 D_\mu})$  contains a connected subset that contains $\cup_{\mu=1}^3 \overline{\frac12 D_\mu}$. Set $D'_\mu=\overline{\frac12 D_\mu}$, $\mu=1,2,3$. Note that the $D'_\mu$ are closed.

%

\textbf{Step 1} Let us first prove that the $n-3$-sphere $\sigma_1=\{x\in Q_1\times Q:|x-p_1|=tR\}=\partial D_1$ represents a non zero element in $H_{n-3}(\R^n\bs \t F_k)$. Here note that $p_1$ is the intersection $(\partial tA^i_{j,k}+q_m)\cap (\mc^i_{jl,k}+q_m)$. 

Let $p$ denote the intersection of $\partial A^i_{j,k}\cap \mc^i_{jl,k}$. Let $Q'$ be the $n-2$ subspace orthogonal to $\mc^i_{jl,k}$ and passing through the origin. Let $\sigma$ be the sphere $\{x\in Q'+p: |x-p|=R\}$. Then $\sigma$ is an $n-3$ sphere that links the planar part $\mc^i_{jl,k}$. Obviously, $\sigma$ is homologic to $\sigma^i_{jl,k}$ in $\R^n\bs \t C_k$. Moreover, let $p(s)=sp+(1-s)\xi^i_{jl,k}$, and $r(s)=sR+(1-s)\eta$, $s\in [0,1]$, and let $\sigma(s)=\{x\in Q'+p(s): |x-p(s)|=r(s)$. Then $\sigma(s)$ is a $n-3$-sphere orthogonal to $\mc^i_{jl,k}$. Let $T$ denote the tube $\cup_{s\in [0,1]}\sigma(s) $. Then the boundary of $T$ is $\sigma^i_{jl,k}\cup \sigma$, and $T$ is in fact the image of the line homotopy between $\sigma^i_{jl,k}$ and $\sigma$. Note that $p$ and $\xi^i_{jl,k}$ both belong to the convex set $\mc^i_{jl,k}$, hence for each $s\in [0,1]$, $p(s)\in \mc^i_{jl,k}$. Then for each $x\in T$, let $s$ be such that $x\in \sigma(s)$. Then the shortest distance projection of $x$ to $\mc^i_{jl,k}$ is $p(s)$. Hence dist$(T, \mc^i_{jl,k})=r(s)\ge \eta$. On the other hand, from the definition of $p$ and $\xi^i_{jl,k}$, we know that the distance of $T$ to all other parts of $\t C^k$ is larger than $\eta$. Hence dist$(T,\t C_k)\ge \eta$.

Now as in the condition in $4^\circ$, since $t\in [s_{m+2}, s_{m-1}]$, and $\t C_k$ is a cone, the tube $tT+q_m$ is of distance larger than $t\eta$ to $\t C_k+q_m$, and $tT+q_m\subset \mU(q_{m-2},s_{m-2})\bs\mU(q_m,s_{m+2})$. On the other hand, we know that $\t F_k$ is $\e s_{m-2}$ near $\t C_k+q_{m-1}$, hence by \tb{(5.17)}, $\t F_k$ is $25 \e s_{m-2}$ near $\t C_k+q_m$ in $\mU(q_{m-2},s_{m-2})\bs\mU(q_m,s_{m+2})$. This means that $\t F_k\cap \mU(q_{m-1},s_{m-1})\bs\mU(q_m,s_{m+2})$ is included in the $25 \e s_{m-2}$ neighborhood of $\t C_k$. Since $\eta>10^3 \e$, and $s_{m+2}\le t\le s_{m-1}$, we know that the tube $tT+q_m$ does not intersect $\t F_k$. This means that, as the boundary of $tT+q_m$, the two spheres $t\sigma+q_m$ and $t\sigma^i_{jl,k}+q_m$ are homologic in $\R^n\bs\t F_k$. By $1^\circ$, we know that $[t\sigma+q_m]=[t\sigma^i_{jl,k}+q_m]=[\sigma^i_{jl,k}]\ne 0$ in $H_{n-3}(\R^n\bs\t F_k)$.

Note that $t\sigma+q_m=\sigma_1$, hence we get 
\be[\sigma_1]\ne 0\mbox{ in }H_{n-3}(\R^n\bs\t F_k).\ee

Similarly, we have
\be[\sigma_\mu]\ne 0\mbox{ in }H_{n-3}(\R^n\bs\t F_k),i=1,2,3.\ee

\textbf{Step 2} 
Now let us prove that $(G^i_{k,m}\cap (\Omega\times Q))\cup (\cup_{\mu=1}^3 D'_\mu)$ containes a connected subset that contains $\cup_{\mu=1}^3 D'_\mu$.  Suppose not.

We know that $G$ is contained in the $\e s_{m-2}$-neighborhood of $Y$, and $G$ is relatively closed in $\Omega\times Q$, hence the limit points of $G$ are contained in $ \overline B(Y,\e s_{m-2})\cap \partial\overline\Omega\times Q\subset \cup_{\mu=1}^3 D'_\mu$. On the other hand, the three $D'_\mu$ are also closed, as a result, $G\cup (\cup_{\mu=1}^3 D'_\mu)$ is closed and bounded, and hence compact.

Write $H=G\cup (\cup_{\mu=1}^3 D'_\mu)$ for short. Then $H$ is a compact subset, containing $\cup_{\mu=1}^3 D'_\mu$. But the three $D'_\mu$ do not belong to the same connected component of $H$. Suppose, without loss of generality, that the connected component $V$ of $H$ that contains $D'_1$ does not intersect $D'_2\cup D'_3$. Then there exists a separation of $H$ that separates $V$ and $D'_2\cup D'_3$, i.e. a pair of disjoint relatively closed non-empty subsets $A$ and $B$ of $H$, whose union is $H$, with $V\subset A$, and $D'_2\cup D'_3\subset B$. Since $H$ is compact, $A$ and $B$ are compact as well. Then the distance $d_{A,B}$ between $A$ and $B$ is strictly positive.

We want to construct a smooth manifold in $\Omega\times Q$ which separates $A$ and $B$, and whose boundary is homologic to $\sigma_1$.

Fix any $0<r_1<r_2<\min\{d_{A,B},10^{-1}R\}$. Let $r_3=\frac {r_1+r_2}{2}$, and let $\e_3=\min\{\frac{r_3-r_1}{8},10^{-2}tR\}$. 


Let $Q'$ be the $n-1$ subspace containing $\Omega\times Q$. Let $\sigma':=\partial (\frac12 +r_3) D_1$. We will construct a surface in $Q'\bs H$, whose boundary is $\sigma'$. Hence in the following text, without precision, in stead of $\R^n$, we will be situated in the $n-1$ dimensional space $Q'$. The words ''open'', ''closed'', ''neighborhood'', ''smooth'' etc. are with respect to  $Q'$.

Set $f:Q'\to \R_+\cup \{0\}$: $f(x)=$dist$(x,A)$. Then $f$ is Lipschitz.

Let $M_0=\{x\in Q'\bs (\Omega\times Q): f(x)=r_3\}$. Then it is a topological $n-2$ ball, whose boundary is $\sigma'$, with $\bar M_0=M_0\cup \sigma'$.

Note that for any point $x$ in a small neighborhood $W_1:=B(\sigma', 5\e_3)\cap Q'$ of $\sigma'$, its distance to $Y$ is larger than $r_3+\frac 12 tR$, and since $G$ is contained in the $s_{m-2}\e$ neighborhood of $Y$, we know that dist$(x,G)\ge r_3+\frac 12 tR-s_{m-2}\e\ge r_3+\frac 14tR$. But dist$(x,D'_\mu)\le r_3+4\e_3$, hence dist$(x,A)=$dist$(x, D'_\mu)$, which is a smooth function on ${D'_\mu}^C$. Therefore $f$ is smooth on $W_1$. For the same reason, $f$ is also smooth in $W_2=B(M_0,5\e_3)$. Note that $W_1\subset W_2$. 
Let $W_3=B(M_0,4\e_3)$, and $W_4=B(M_0,3\e_3)$. Then there exists a smooth map $g:Q'\to \R$ such that 
\be g|_{\bar W_3}=f|_{\bar W_3} \mbox{ and } ||g-f||_\infty<\e_3.\ee


Moreover, since $D'_\mu$ is convex and for any $x\in \bar W_2$, the closest point in $A$ to $x$ belongs to $D'_\mu$,we know that for each $x\in W_2$, there exists a unique nearest point $\xi(x)$ in $A$, i.e. $f(x)=||x-\xi(x)||$.

By \cite{Fe59} Theorem 4.8, $D f(x)=2(x-\xi(x))\ne 0$. Hence the restriction of $f:W_2\to \R$ is smooth, and transverse to the 0-dimensional submanifold $\{r_3\}$ of $\R$, for any $x$ in $W_2$. Since $g$ coincides with $f$ in $W_3$, we know that $g$ is transverse to the 0-dimensional submanifold $\{r_3\}$ of $\R$ at every point of $W_3$, and hence for every point of the closed subset $\bar W_4$. By Thom transversality Theorem (cf. \cite{Bre} Chapt II, Thm 15.2), there exists a map $h$ defined on $Q'$, such that $h$ is transverse to $\{r_3\}$, and 
\be h|_{\bar W_4}=g|_{\bar W_4}=f|_{\bar W_4}\mbox{, and } ||h-g||_\infty<\e_3. \ee

Since $h$ is transverse to $\{r_3\}$, we know that $h^{-1}\{r_3\}$ is a $n-2$ dimensional smooth submanifold $M$ of $Q'$. And since $h$ coincides with $f$ in $\bar W_4$, we know that $M\cap \bar W_4\bs (\Omega^\circ\times Q)=\bar M_0$.

On the other hand, for any $x\in Q'\bs(\bar W_4\cup (\Omega^\circ\times Q))$, we know that $h(x)\ge g(x)-\e_3\ge f(x)-\e_3-\e_3\ge 2\e_3>\e_3$, hence $M\bs(\bar W_4\cup (\Omega\times Q))=\emptyset$. As a result, 
\be M\bs (\Omega^\circ\times Q)=\bar M_0.\ee 

Next we look at $M\cap (\Omega^\circ\times Q)$. Since $h(x)=r_3$, we know that $f(x)\in (r_3-2\e_3,r_3+2\e_3)$, which means that dist$(x,A)\in (r_3-2\e_3,r_3+2\e_3)\subset (r_1,r_2)$. Therefore $x\not\in A\cup B$. Hence $M\cap H=M\cap (A\cup B)=\emptyset$. Therefore, $M$ is a smooth $n-2$-dimensional submanifold in $Q'$, $M$ does not touch $H$, and $M\bs (\Omega^\circ\times Q)=\bar M_0$, which is a smooth $n-2$ ball, whose boundary is $\sigma'$. As consequence, the boundary of the part $M_1=M\cap (\Omega^\circ\times Q)$ is also $\sigma'$. 

This yields that $[\sigma']=0$ in $H_{n-3}(\Omega\times Q\bs (A\cup B))$. But $\Omega\times Q\bs (A\cup B)\subset \R^n\bs \t F_k$, hence $[\sigma']=0$ in $H_{n-3}(\R^n\bs \t F_k)$.

On the other hand, note that $D_1\bs (\frac 12+r_3)D_1$ does not meet $\t F_k$, and its boundary is $\sigma'\cup \sigma_1$, hence $[\sigma_1]=[\sigma']=0$ in $H_{n-3}(\R^n\bs \t F_k)$. But this contradicts \tb{(6.17)}.

This finishes the proof of $5^\circ$. 

\qed

\subsection{Structure of $F_k$ inside $\mU(q_m,s_m)$ if the $\e$-process does not stop at step $m$}

In this subsection we discuss the structure of $F_k$ in $\mU(q_m,s_m)$ if the $\e$-process does not stop at step $m$.

\begin{pro}There exists $\e_0\in (0,\min\{10^{-5}, 10^{-3}\e_1(n,2)\eta,\e_2(n)\})$, ($\e_2$ being the one in Corollary \tb{2.24}, such that for all $\e<\e_0$, $k$ large enough,  if the $\e-$process does not stop before $m$, then there exists a Lipschitz map
 $g_m: B\to \bar \mU(q_m,s_m)$ such that 
 \be g_m=id\mbox{ on }\bar\mU(q_m,s_m),\ee
 and
 \be g_m(F_k)=F_k\cap \mU(q_m,s_m).\ee
\end{pro}

\nd Let us fix $k$, $\e<\e_0$ small, to be decided. Suppose the $\e$-process does not stop at step $m$. 

We will prove by induction. So for $m=1$, we know that $F_k$ is $\e$ close to $C_k+q_1$ in $\bar B(0,1)$. We claim that,
\be \mbox{for each }x\in F_k\cap \partial \mU(q_1,\frac 12), x\mbox{ is either a regular point, or a }\Y\mbox{ point for }F_k.\ee

In fact, we know that in $B(x,\frac 18\eta)$, $F_k$ is $\e$-close to $C_k+q_m$. But $C_k+q_m$ is either a plane or a $\Y$ set in $B(x,\frac 18\eta+\e)$, therefore $F_k$ is $\e$-close to a plane or a $\Y$ set in $B(x,\frac 18\eta+\e)$. 

If $F_k$ is $\e$-close to a plane in $B(x,\frac 18\eta+\e)$, and since $\e<10^{-3}\e_1(n,2)\eta$, hence $F_k$ is $\e_1(n,2)[\frac 18\eta+\e]$ near a plane in $B(x, \frac 18\eta+\e)$. Then Theorem \tb{2.25} of \cite{2p} yields that $x$ must be a regular point of $F_k$;

If $F_k$ is $\e$-close to a 2-dimensional $\Y$ set $Y$ centered at $y$ in $B(x,\frac 18\eta+\e)$.By Lemma 16.43 of \cite{DJT}, when $\e$ is small enough, the measure of $F_k$ in $B(x,\frac 18\eta)$ should be smaller than the measure of $(Y+y)\cap B(x,\frac 18\eta))<d_T(\frac 18\eta)^2$, where $d_T$ is the smallest density of 2-dimensional minimal cone in $\R^n$ other than a plane or a $\Y$ set. Since $F_k$ is minimal in $B(x, \frac 14\eta)$, by monotonicity of density for minimal sets (cf. \cite{DJT} Proposition 5.16, or apply directly the monotonicity of density for stationary varifold \cite{Al76}), we know that the density of $F_k$ at $x$ is smaller than $d_T$ ($d_T$ is the one in Remark \tb{2.14} $2^\circ$), and hence $x$ can only be a regular point, or a $\Y$ point. 

Thus we get Claim \tb{(6.24)}. 

Moreover, if $C_x$ is the tangent cone of $F_k$ at $x$, then $C_x$ is transversal to $\partial \mU(q_1,\frac 12)$, since $C_x$ should also be very close to $C_k+q_1$ in $B(x,\frac 18\eta)$. As a result, since $\e<\e_2$, by Corollary \tb{2.24}, for each $x\in F_k\cap \partial \mU(q_1,\frac 12)$, there exists a 2-Lipschitz deformation retract $\varphi_x$ from $\partial \mU(q_1,\frac 12)\cap B(x,r_1)$ to $F_k\cap \partial \mU(q_1,\frac 12)\cap B(x,r_1)$, where $r_1=\frac 18\eta$. Also, $(C_k+q_1)\cap \partial \mU(q_1,\frac 12)$ is included in the $5\e$ neighborhood of $F_k\cap \partial \mU(q_1,\frac 12)$.

Since $E=B(F_k\cap \partial \mU(q_1,\frac 12),\frac 12 r_1)$ is compact, we know that it can be covered by finitely many $B(x_l,r_1), 1\le l\le l_1$. Then by partition of unity, we know that there exists a 2-Lipschitz neighborhood deformation retract $\varphi_1$ from $B(F_k,\frac12 r_1)\cap \partial \mU(q_1,\frac 12)$ to $F_k\cap \partial \mU(q_1,\frac 12)$.

Let $\pi_1$ denote the radial projection from $B(0,1)$ to $\mU(q_1,s_1)$. Then $\pi_1(F_k\bs \mU(q_1,s_1))\subset \partial \mU(q_1,s_1)$. Since $F_k$ is $\e$ near $C_k+q_1$ in $B(0,1)$, and $C_k+q_1$ is a cone centered at $q_1$, we know that $\pi_1(F_k\bs \mU(q_1,s_1))$ is included in the $2\e$ neighborhood of $(C_k+q_1)\cap \partial \mU(q_1,s_1)$, and hence is included in the $10\e$ neighborhood of $F_k\cap \partial \mU(q_1,s_1)$, which is inclued in $B(F_k,\frac12 r_1)\cap \partial \mU(q_1,\frac 12)$. As a result, $\varphi_1(\pi_1(F_k\bs \mU(q_1,s_1))\subset F_k\cap \partial \mU(q_1,s_1)$. We extend $\varphi_1$ to a global Lipschitz deformation, such that $\varphi_1=id$ on $B(0,1)^C$ and $F_k\cap \bar\mU(q_1,s_1)$. Then since $\pi_1=id$ on $\mU(q_1,s_1)$, we get
\be \varphi_1\circ\pi_1(F_k)=F_k\cap \bar\mU(q_1,s_1).\ee

Set $g_1=h_1=\varphi_1\circ\pi_1$.

We repeat this process, and get, a deformation $h_2$ in $\bar \mU(q_1,s_1)$ that maps $F_k\cap \bar \mU(q_1,s_1)$ to $\bar\mU(q_2,s_2)$. Set $g_2=h_2\circ h_1$.

By recurrence, for each $m$ such that the $\e$-process does not stop, we get a Lipschitz deformation $g_m$ from $B(0,1)$ to $B(0,1)$ that maps $F_k$ to $F_k\cap \mU(q_m,s_m)$. \qed

A similar argument gives that 

\begin{pro}There exists $\e_0\in (0,\min\{10^{-5}, 10^{-3}\eta, \e_2\})$, such that for all $\e<\e_0$, $k$ large enough, and for every $m$ such that the $\e-$process does not stop before $m$, we have: 
for each $i=1,2$, and $t\in [\frac 14, 1]$, the set $F_k\cap \mU(q_m, ts_m)\in \oF_{4\e s_m}((C_k+q_m)\cap \bar\mU(q_m, ts_m), \bar\mU(q_m, ts_m))$
\end{pro}

\nd In fact, a similar argument as in the proof of the previous proposition gives that: for every $m$ such that the $\e-$process does not stop before $m$, we have: 
for each $i=1,2$, and $t\in [\frac 14, 1]$, there exists a Lipschitz map $g=g_{m,t}:\mU(0,1)\to \mU(q_m, ts_m)$, such that $g(F_k)=F_k\cap \mU(q_m, ts_m)$.

Let $h:\bar\mU(q_m, ts_m)\to \bar\mU(0,1): h(x)=\frac{x-q_m}{ts_m}$, then $g\circ h:\bar\mU(q_m,ts_m)\to \bar\mU(q_m,ts_m)$ is a Lipschitz deformation, and $g\circ h(h^{-1}(F_k))=F_k\cap \bar\mU(q_m, ts_m))$. But $F_k\in \oF(C_k,\mU(0,1))$, hence $g\circ h(h^{-1}(F_k))\in \oF((C_k+q_m)\cap \bar\mU(q_m, ts_m), \bar\mU(q_m, ts_m))$. Moreover, since $g(F_k\bs\mU(q_m,ts_m))\subset F_k\cap \mU(q_m,ts_m)\subset B(C_k+q_m, \e s_m)$, we know that $F_k\cap \bar\mU(q_m,ts_m)\in \oF_{4\e s_m}((C_k+q_m)\cap \bar\mU(q_m, ts_m), \bar\mU(q_m, ts_m))$.\qed

\section{Measure estimates outside the central ball}

With the help of Proposition \tb{6.1}, we give estimates of measures of $F_k$ outside the critical ball in this section. Now set $\e_0=\min\{10^{-5}, 10^{-3}\eta, \e_2, 10^{-2}\d_0\}$, where $\d_0$ is the one selected before \tb{(5.1)}.

Recall that for each fixed $\e<\e_0$ and $k$ large, the $\e-$process stops at step $m$, and we set $o_k=q_m, r_k=s_m$. 

We will prove the following theorem:

\begin{thm}For each $k\in \N$, 
\be\H^2(F_k\bs \mU(o_k,\frac14 r_k))\ge \H^2(C_k\bs \mU(o_k,\frac14 r_k))+C(\e, K^1,K^2)r_k^2.\ee
\end{thm}

The rest of this section will be devoted to prove Theorem \tb{7.1}. The proof will be decomposed into a series of lemmas and propositions. 

\subsection{Measure estimates for zones far from the center}

This subsection will be devoted to the proof of the following proposition:

\begin{pro}For every $m$ such that the $\e$-process does not stop before $m$, 
\be \H^2(F_k\cap \mU(q_{m-1}, s_{m-1})\bs \mU(q_m,s_m))\ge \H^2((C_k+q_m)\cap \mU(q_{m-1}, s_{m-1})\bs \mU(q_m,s_m)).\ee
\end{pro}

\nd
Let $\Omega$ denote $\mU(q_{m-1}, s_{m-1})\bs \mU(q_m,s_m)$. We would like to decompose it into intersections with regular regions $\Omega\cap (\mC^i_{jl,k}\times {P_k^i}^\perp+q_m), (j,l)\in J^i$, $\Omega\cap (\mS^i_{j,k}\times {P_k^i}^\perp+q_m), 1\le j\le \mu_i$, $i=1,2$, and singular regions $\Omega\cap (\mA^i_{j,k}\times {P_k^i}^\perp+q_m),i=1,2, 1\le j\le m_i$, where we can get the measure estimates for each part:

$1^\circ$ In regular regions: fix either $i=1,2$ and $(j,l)\in J^i$. Then Proposition \tb{6.1} $3^\circ$ yields that $F_k\cap \Omega\cap(\mC^i_{jl,k}\times {P_k^i}^\perp+q_m)$ is a graph on $(\mc^i_{jl,k}+q_m)\cap \Omega\cap(\mC^i_{jl,k}\times {P_k^i}^\perp+q_m)$, which is just $C_k+q_m)\cap \Omega\cap(\mC^i_{jl,k}\times {P_k^i}^\perp+q_m)$. Hence
\be \H^2[F_k\cap \Omega\cap(\mC^i_{jl,k}\times {P_k^i}^\perp+q_m)]\ge\H^2[C_k+q_m)\cap \Omega\cap(\mC^i_{jl,k}\times {P_k^i}^\perp+q_m)].\ee

Similarly by Proposition \tb{6.1} $4^\circ$, we have, for $i=1,2$ and $1\le j\le \mu_i$,
\be \H^2[F_k\cap \Omega\cap(\mS^i_{j,k}\times {P_k^i}^\perp+q_m)]\ge\H^2[C_k+q_m)\cap \Omega\cap(\mS^i_{j,k}\times {P_k^i}^\perp+q_m)].\ee

$2^\circ$ In singular regions: fix either $i=1,2$, and $j\le m_i$. 

\begin{lem} For each $1\le j\le \mu_i$, $i=1,2$, let $E\subset A^i_j\times {P_0^i}^\perp$ be a closed 1-dimensional set, such that $E$ connects the three $I^i_{jl}\times {P_0^i}^\perp$.  Let $Y=C_0\cap A^i_j\times {P_0^i}^\perp$. Then
\be \H^1(Y)\le \H^1(E).\ee
\end{lem}
\nd Fix $i$ and $j$. Let $I_1,I_2,I_3$ denotes the three $n_i-1$-disks $I^i_{jl}, (j,l)\in J^i$ for short, and let $a_l$ denotes the center of $I_\a,\a=1,2,3$. Then $a_\a,\a=1,2,3$ is also the three points of intersection of $Y$ with $\partial A^i_j$. Let $P$ denote the 2-plane containing $Y$, and hence $a_\a\in P, \a=1,2,3$. Let $\pi$ denote the orthogonal projection from $A^i_j\times {P_0^i}^\perp$ to $P$. 

Now we project everything to $P$, and see what happens. Let $b$ denote $\pi(a^i_j)$. The projection of $A^i_j\times P_0^i$ is the set 
\be \pi(A^i_j\times P_0^i)=\pi(A^i_j)=B(b,R)\bs (\cup_{\a=1}^3\{x\in B(b,R): <x,a_\a>>1-\eta\}) ,\ee
and hence the boundary of $\pi(A^i_j\times {P_0^i}^\perp)$ is a union of three segments $S_\a=\pi(I_\a)$ with the same length $R_1=2\sqrt{1-(1-\eta)^2}$, the midpoint of $S_\a$ being $a_\a$, and the union of the three arcs of circles $\xi_\a,\a=1,2,3$, which are parts of $\partial B(b,R)$. For $\a=1,2,3$, the arc of circle $\xi_\a$ is on the opposite side of $S_\a$.
See the picture below. 

\centerline{\includegraphics[width=0.8\textwidth]{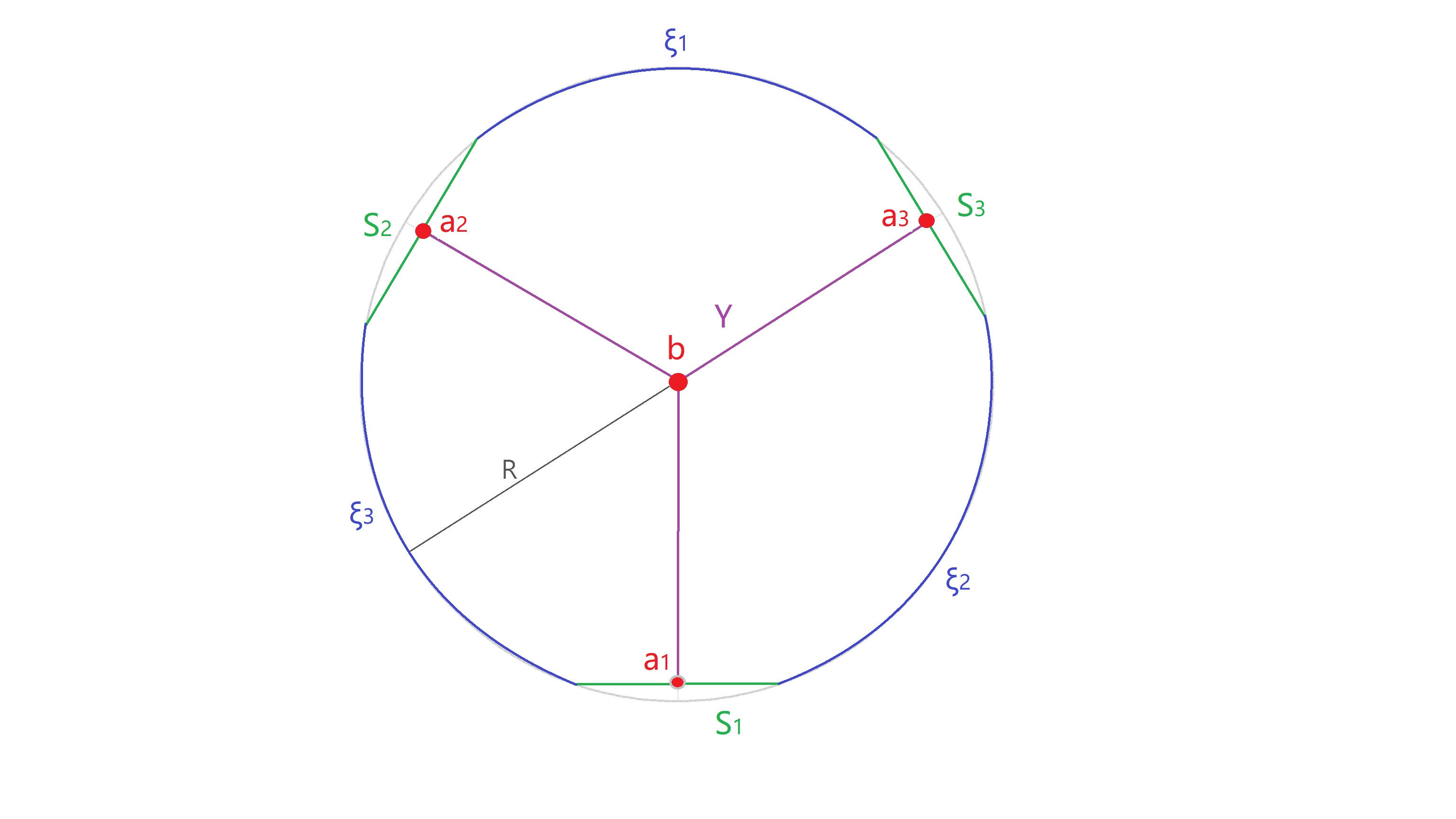} }

We also know that $\pi(Y)=Y$, and $\pi(a_\a)=a_\a,\a=1,2,3$, since they belong to $P$.

Since $E$ connects the three $I_\a$, the projection $\pi(E)$ also connects the three $S_\a=\pi(I_\a), \a=1,2,3$. Let $E_0\subset \pi(E)$ be a subset of $\pi(E)$, such that $E_0\cup(\cup_{\a=1}^3S_\a)$ is connected. Then $E_0$ intersects every $S_\a$. For each $\a=1,2,3$, fix $b_\a\in S_\a\cap E_0$. Then $E_0$ is a connected set that contains $b_\a$. Let $c$ be the Fermat point of the three points $b_\a,1\le \a\le 3$. Then we have 
\be\H^1(E_0)\ge \sum_{\a=1}^3\H^1([cb_\a]).\ee 

\centerline{\includegraphics[width=0.5\textwidth]{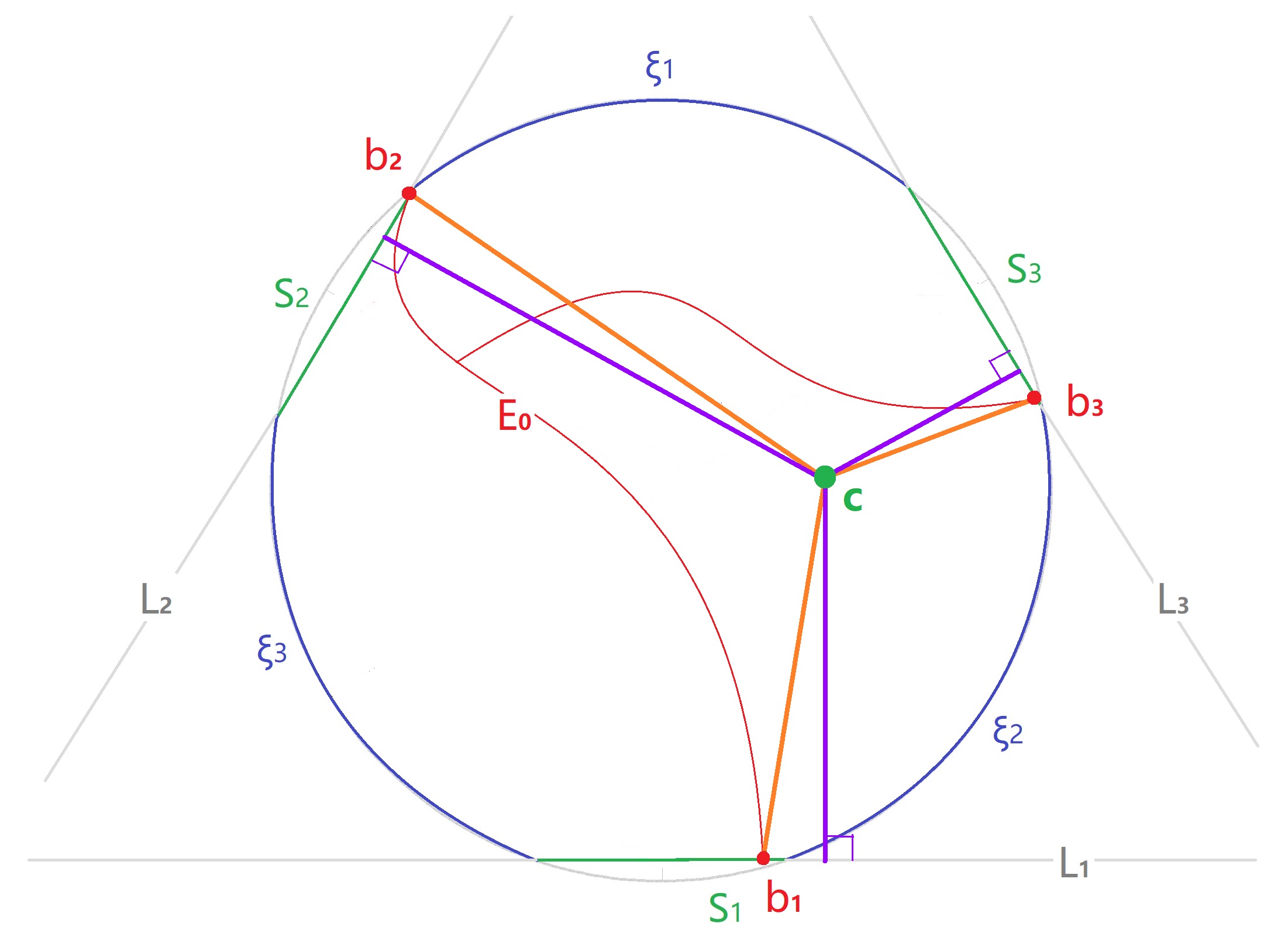} }

And since the triangle $\Delta_{b_1b_2b_3}$ has no angle larger or equal to $120^\circ$, we know that $c\in \Delta_{b_1b_2b_3}^\circ\subset \pi(A^i_j\times P_0^i)$, and the angles between each pair of $[cb_\a],1\le \a\le 3$ is $120^\circ$. 

Denote by $L_\a$ the line containing $S_\a$, $1\le \a\le 3$. Then we have
\be \sum_{\a=1}^3\H^1([cb_\a])\ge\sum_{\a=1}^3\mbox{ dist}(c,S_\a)\ge\sum_{\a=1}^3\mbox{ dist}(c,L_\a).\ee

Note that $c$ is contained in $\pi(A^i_j\times P_0^i)$, which is contained in the equilateral triangle $\Delta$ enclosed by the three $L_\a,1\le \a\le 3$. And it is easy to notice that since $\Delta$ is equilateral, for all $x\in \Delta$,  the quantity 
\be \sum_{\a=1}^3\mbox{ dist}(x,L_\a)\ee
are the same. Hence
\be \sum_{\a=1}^3\mbox{ dist}(c,L_\a)=\sum_{\a=1}^3\mbox{ dist}(\pi(b),L_\a)=\H^1(Y).\ee

Combining \tb{(7.7), (7.8) and (7.10)}, we have
\be \H^1(E_0)\ge \H^1(Y).\ee

But $E_0$ is a subset of $\pi(E)$, hence 
\be \H^1(E)\ge\H^1(\pi(E))\ge \H^1(E_0)\ge \H^1(Y),\ee
which yields the conclusion of the Lemma.\qed

As a corollary, we get that
\begin{cor}For each $i=1,2$, $k$ large, and $j\le m_i$, we have
\be \H^1(F_k\cap t(A^i_{j,k}\times {P^i_k}^\perp)+q_m)\ge \H^1((C_k+q_m)\cap t(A^i_{j,k}\times {P^i_k}^\perp)+q_m),\ee
for each $t\in [s_{m+1}, s_{m-1}]$.
\end{cor}

\nd Fix a $i,k$, and $j$. Let $\psi$ be an isometry from $\R^n=\R^{n_i}\times \R^{n-n_i}$ to $\R^n=P_k^i\times {P_k^i}^\perp$, such that $\psi|_{\R^{n_i}}=\psi^i_k(x)$, and $\psi|_{\R^{n-n_i}}$ is any linear isometry.  Here $\psi^i_k:\R^{n_i}\to P_k^i$ is the isometry that maps $K^i$ to $K^i_k$, defined just before \tb{(5.6)}. Set $f(x)=\psi(\frac 1t (x-q_m))$.
Then it is easy to see that for any set $F\subset \R^n$, 
\be\H^1(\psi(F))=t\H^1(F).\ee

Let $E=f[F_k\cap t(A^i_{j,k}\times {P^i_k}^\perp)+q_m]$, then by Proposition \tb{6.1} $5^\circ$, $E$ connects the three $I^i_{jl}\times {P_0^i}^\perp$. Hence $\H^1(E)\ge \H^1(Y)$. On the other hand, the set $Y$ is just $f((C_k+q_m)\cap t(A^i_{j,k}\times {P^i_k}^\perp)+q_m)$. Thus the conclusion of the corollary follows by \tb{(7.14)}.\qed

We would like to apply the Corollary \tb{7.4} to prove 
\be\begin{split}\H^2(F_k&\cap (\mU(q_{m-1}, s_{m-1})\bs \mU(q_m,s_m))\cap ( (\mA^i_{j,k}\times {P^i_k}^\perp)+q_m))\\
&\ge \H^2((C_k+q_m)\cap (\mU(q_{m-1}, s_{m-1})\bs \mU(q_m,s_m))\cap ( (\mA^i_{j,k}\times {P^i_k}^\perp)+q_m)),
\end{split}\ee
using slicing method. In fact, this is ok if we replace $q_{m-1}$ by $q_m$ in the above equation. In the above equation, due to the very tiny shift of the center, the shape of the intersection of $\mA^i_j\times {P_k^i}^\perp+q_m$ with $\Omega=\mU(q_{m-1}, s_{m-1})\bs \mU(q_m,s_m)$ is a little bit different near the boundary of $\mU(q_{m-1}, s_{m-1})$. In particular,  the slices $(\mU(q_{m-1}, s_{m-1})\bs \mU(q_m,s_m))\cap ( t(A^i_{j,k}\times {P^i_k}^\perp)+q_m)))$ may be only part of $ t(A^i_{j,k}\times {P^i_k}^\perp)+q_m)$, which makes it inconvenient to apply Corollary \tb{7.4}. See the picture below.

\centerline{\includegraphics[width=0.8\textwidth]{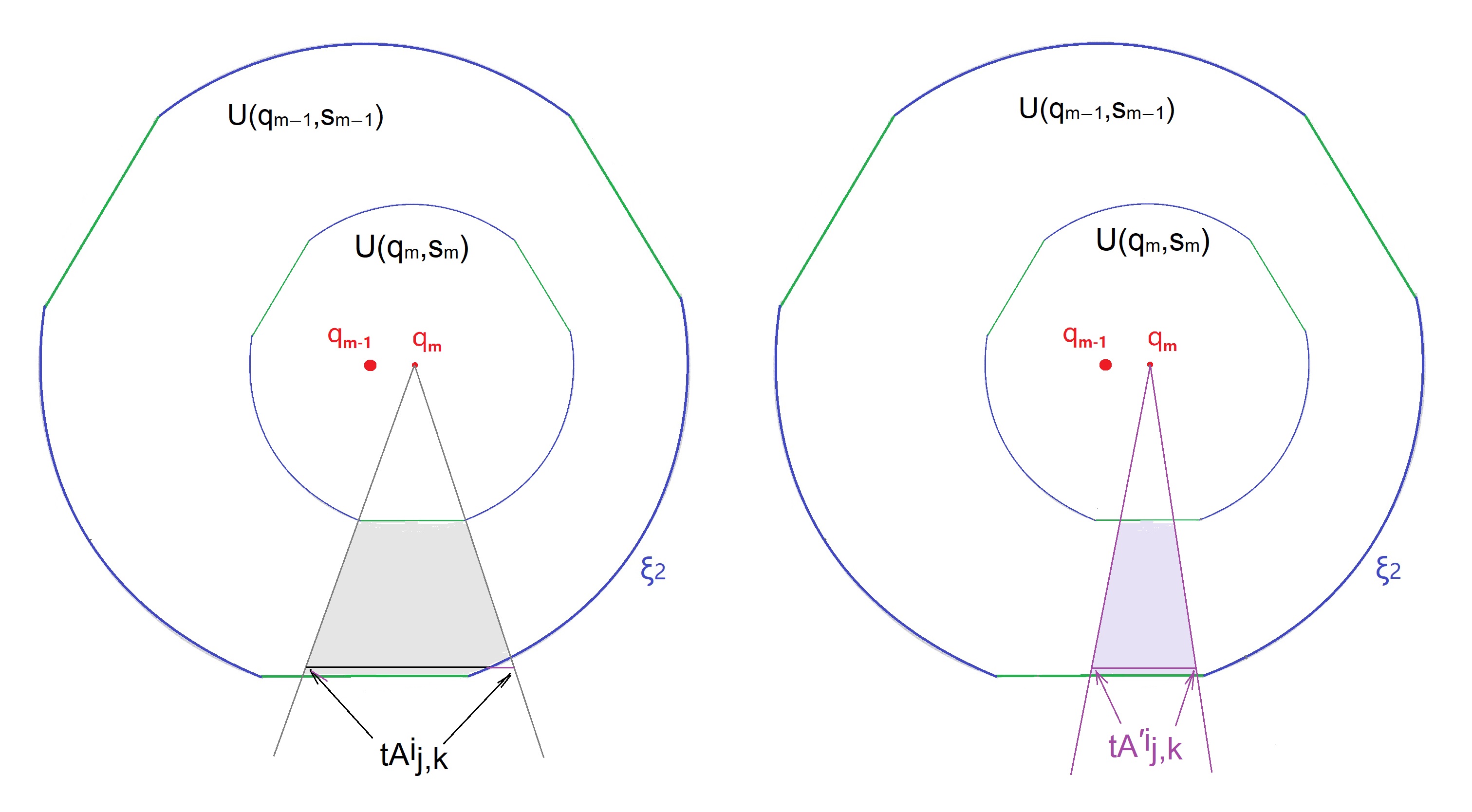} }

 So we need a slightly different decomposition of $\Omega\cap F_k$. More precisely, instead of taking $A^i_j$ around singular part, we take a similar but smaller region:

Let $ {A'}^i_j=\{a^i_j+\frac 12 x:a^i_j+x\in A^i_j\}$. Then ${A'}^i_j\subset A^i_j$ is of the same shape with the same center as $A^i_j$, but with half the diameter; set ${\Gamma'}^i_{jl}=\{x\in B^i,<x,y>=1-\eta\mbox{ for some }y\in \gamma^i_{jl}\}\bs {A'}^i_{j}$. Let ${\mA'}^i_j$ and ${\mC'}^i_{jl}$ denote the cone over ${A'}^i_j$ and ${\Gamma'}^i_{jl}$, and ${A'}^i_{j,k}$, ${\Gamma'}^i_{jl,k}$, ${\mA'}^i_{j,k}$ and ${\mC'}^i_{jl_k}$ denote the isometry copy in $P^i_k$ as before. Then this is a similar decomposition as before, in particular, the exactly same argument gives that
\be \H^2[F_k\cap \Omega\cap({\mC'}^i_{jl,k}\times {P_k^i}^\perp+q_m)]\ge\H^2[(C_k+q_m)\cap \Omega\cap({\mC'}^i_{jl,k}\times {P_k^i}^\perp+q_m)].\ee
and
\be \H^2(F_k\cap t({A'}^i_{j,k}\times {P^i_k}^\perp)+q_m)\ge \H^2((C_k+q_m)\cap t({A'}^i_{j,k}\times {P^i_k}^\perp)+q_m),\ee
 moreover, it satisfies that, for any $t$, 
\be(\mU(q_{m-1}, s_{m-1})\bs \mU(q_m,s_m))\cap ( t({A'}^i_{j,k}\times {P^i_k}^\perp)+q_m)\mbox{ is either } t({A'}^i_{j,k}\times {P^i_k}^\perp)+q_m)\mbox{ or }\emptyset.\ee

As we said before, we use slicing method to prove \tb{(7.15)}, replacing $\mA^i_{j,k}$ by ${\mA'}^i_{j,k}$: 

Set $f:\mU(q_{m-1}, s_{m-1})\bs \mU(q_m,s_m)\cap ({\mA'}^i_{j,k}\times {P^i_k}^\perp+q_m)\to \R$, $f(x)=t$ if $<f(x)-q_m,a^i_{j,k}>=t$. Then $f$ is essentially an orthogonal projection to the line parallel to $\vec{a^i_{j,k}}$ and passing through $q_m$. Hence $f$ is 1-Lipschitz.  By the coarea formula \cite{Fe} 3.2.22, we have
\be\begin{split}&\H^2(F_k\cap (\mU(q_{m-1}, s_{m-1})\bs \mU(q_m,s_m))\cap ( ({\mA'}^i_{j,k}\times {P^i_k}^\perp)+q_m))\\
&\ge\int_\R \H^1(F_k\cap (\mU(q_{m-1}, s_{m-1})\bs \mU(q_m,s_m))\cap ( ({\mA'}^i_{j,k}\times {P^i_k}^\perp)+q_m)\cap f^{-1}(t))dt\\
&=\int_\R \H^1(F_k\cap (\mU(q_{m-1}, s_{m-1})\bs \mU(q_m,s_m))\cap (t({A'}^i_{j,k}\times {P^i_k}^\perp)+q_m))dt\\
&=\int_T \H^1(F_k\cap (\mU(q_{m-1}, s_{m-1})\bs \mU(q_m,s_m))\cap  (t({A'}^i_{j,k}\times {P^i_k}^\perp)+q_m))dt,
\end{split}\ee
Where $T=\{t\in \R: (\mU(q_m, s_{m-1})\bs \mU(q_m,s_m))\cap ( t({A'}^i_{j,k}\times {P^i_k}^\perp)+q_m)\ne\emptyset$. By (7.18), and then (7.17), we have
\be\begin{split}&\int_T \H^1(F_k\cap (\mU(q_{m-1}, s_{m-1})\bs \mU(q_m,s_m))\cap  (t({A'}^i_{j,k}\times {P^i_k}^\perp)+q_m))dt\\
&=\int_T \H^1(F_k\cap (t({A'}^i_{j,k}\times {P^i_k}^\perp)+q_m))dt\\
&\ge \int_T \H^1((C_k+q_m)\cap  (t({A'}^i_{j,k}\times {P^i_k}^\perp)+q_m))dt\\
&=\int_T \H^1((C_k+q_m)\cap (\mU(q_{m-1}, s_{m-1})\bs \mU(q_m,s_m))\cap  (t({A'}^i_{j,k}\times {P^i_k}^\perp)+q_m))dt\\
&=\int_\R \H^1((C_k+q_m)\cap (\mU(q_{m-1}, s_{m-1})\bs \mU(q_m,s_m))\cap ( ({\mA'}^i_{j,k}\times {P^i_k}^\perp)+q_m)\cap f^{-1}(t))dt\\
&=\H^2((C_k+q_m)\cap (\mU(q_{m-1}, s_{m-1})\bs \mU(q_m,s_m))\cap ( ({\mA'}^i_{j,k}\times {P^i_k}^\perp)+q_m)).
\end{split}\ee
Thus \tb{(7.19) and (7.20)} give
\be\begin{split}\H^2(F_k&\cap (\mU(q_{m-1}, s_{m-1})\bs \mU(q_m,s_m))\cap ( ({\mA'}^i_{j,k}\times {P^i_k}^\perp)+q_m))\\
&\ge \H^2((C_k+q_m)\cap (\mU(q_{m-1}, s_{m-1})\bs \mU(q_m,s_m))\cap ( ({\mA'}^i_{j,k}\times {P^i_k}^\perp)+q_m)),
\end{split}\ee

We apply \tb{(7.16), (7.21) and (7.4)}, and sum up over $1\le j\le m_i$, $(j,l)\in J^i$ and $1\le j\le \mu_i$ respectively, this yields 
\tb{(7.2)}. This finishes the proof of Proposition \tb{7.2}. \qed

\begin{cor}Let $\e<\e_0$. For every $m$ such that the $\e$-process does not stop before $m$, 
\be \H^2(F_k\bs \mU(q_{m},s_m))\ge \H^2(C_k\bs \mU(0,s_m)).\ee
\end{cor}

\nd 

Note that for all $m'\le m$, we have $d(q_{m-1},q_m)\le 24\e s_{m-1}<\d_0s_{m-1}$. By Remark \tb{2.32}, we know that $C_k$ is $(\eta, \d_0)$-measure stable. Hence by Proposition \tb{7.2},
\be \begin{split}&\H^2((C_k+q_{m'})\cap (\mU(q_{m'-1}, s_{m'-1})\bs \mU(q_{m'},s_{m'}))\\
&=\H^2(C_k\cap \mU(q_{m'-1}-q_{m'}, s_{m'-1})\bs \mU(0,s_{m'}))\\
&=\H^2(C_k\cap \mU(0, s_{m'-1})\bs \mU(0,s_{m'})).
\end{split}\ee

We sum up over all $m'\le m$, and get \tb{(7.22)}.\qed

\subsection{Measure estimates near the critical zone}

Now let us look at the critical zone $\mU(o_k, r_k)\bs\mU(o_k,\frac14 r_k)$. Recall that we defined $o_k=q_{m_k}$, and $r_k=s_{m_k}$, where $m_k$ is the step where our $\e$-process stops. We know that in $\mU(o_k, r_k)$, $F_k$ is not $\e r_k$ near any translation of $C_k$, but is $2\e r_k$ near $C_k+o_k$, because the $\e$-process does not stop at step $m-1$. Then the next  proposition will tell us that $F_k$ is also far from any translation of $C_k$ in $\mU(o_k, r_k)\bs\mU(o_k,\frac14 r_k)$. For future use, we will need it in a more general setting. So let recall that: 

For any $p_1,p_2\in \bar B(0,\frac 14)$, the set $C_{p_1,p_2}$ denote the translated orthogonal union $(K^1+p_1)\cup_\perp (K^2+ p_2)$. By Proposition \tb{2.36}, $C_{p_1,p_2}$ is minimal in $\R^n$; by Theorem \ref{unicite1}, it is Almgren unique. 

Set $C_{k,p_1,p_2}=(K^1_k+p_1)+(K^2_k+p_2)$.

\begin{pro}For each $0<\e<\frac13 \e_0$, there exists a $\d(\e)\in (0,\e)$, and $k_0\in \N$, such that  the following holds: 

For any $\d<\d(\e)$, any $k>k_0$, and any $p_1,p_2\in B(0,\frac 14)$, suppose that $E\in\oF_\d(C_{k,p_1,p_2},\mU_k(0,1))$, $E$ is minimal in $\mU(0,1)$, and is $\e_0$-near $C_{k,p_1,p_2}$ in $\mU_k(0,1)$, $\d$-near $C_{k,p_1,p_2}$ in $\mU_k(0,1)\bs \mU_k(0,\frac 12)$, then $E$ must be $\e$-near $C_{k,p_1,p_2}$ in $\mU_k(0,1)$.
\end{pro}

\nd The proof is a simple compactness argument. 

Suppose the conclusion of the proposition is not true. Then there exists $\e\in (0,\frac13 \e_0)$, and four sequences $\d_l\to 0$, $k_l\to\infty$, $p_{1,l},p_{2,l}\in \bar B(0,\frac 14)$, and a sequence of minimal sets $E_l\in \oF_\d(C_{k_l,p_{1,l}, p_{2,l}},\mU_{k_l}(0,1))$, such that $E_l$ is $\d_l$ near $C_{k_l,p_{1,l}, p_{2,l}}$ in $\mU(0,1)\bs \mU(0,\frac 12)$, but not $\e$ near $C_{k_l,p_{1,l}, p_{2,l}}$ in $\mU_{k_l}(0,1)$.

Modolu extracting a subsequence, we can suppose that $p_{i,l}$ converges to $p_i$, $i=1,2$.  Then $C_{k_l,p_{1,l}, p_{2,l}}$ converges to the orthogonal union $C_{0,p_1,p_2}$. Since $E_l$ is $\d_l$ near $C_{k_l,p_{1,l}, p_{2,l}}$ with $\d_l\to 0$, there exists a sequence $a_l\to 0$ such that $E_l$ is $a_l$ near $C_{0,p_1, p_2}$ in $\mU_{k_l}(0,1)\bs \mU_{k_l}(0,\frac 12)$, but not $\e$ near $C_{0,p_1,p_2}$ in $\mU_{k_l}(0,1)$.

Modulo extracting a subsequence, we can suppose that $E_l$ converges to a limit $E_\infty$. Then $E_\infty$ is minimal in $\mU(0,1)$, not $\e$ near $C_{0,p_1,p_2}$ in $\mU_0(0,1)$, and
\be E_\infty\cap \mU_0(0,1)\bs\mU_0(0,\frac 12)=C_{0,p_1,p_2}\cap \mU_0(0,1)\bs\mU_0(0,\frac 12).\ee

Since for each $l$, $E_l\in \oF_\d(C_{k_l,p_{1,l}, p_{2,l}},\mU_{k_l}(0,1))$, we know that $E_\infty\in \oF_{2\d}(C_{0,p_1,p_2}, \mU_0(0,1))$. But by \tb{(7.24)}, we know that 
\be E_\infty\in \oF(C_{0,p_1,p_2}, B(0,1)).\ee

On the other hand, for each $l$, $E_l$ is $\e_0$-close to $C_{k_l,p_{1,l}, p_{2,l}}$ in $\mU_{k_l}(0,1)$, hence $E_\infty$ is $\e_0$-near $C_{0,p_1,p_2}$ in $B(0,1)$. Note that by regularity of 2-dimensional minimal cones, there exists a Lipschitz neighborhood deformation retract $\varphi$ from the $\e_0$ neighborhood of $C_{0,p_1,p_2}$ to $C_{0,p_1,p_2}$ in $B(0,1)$. Then $\varphi(E_\infty)\in \oF(C_{0,p_1,p_2}, B(0,1))$, and $\varphi(E_\infty)\subset C_{0,p_1,p_2}$. This means that 
\be \H^2(C_{0,p_1,p_2})\ge \H^2(\varphi(E_\infty)).\ee
But $E_\infty$ is minmal in $B(0,1)$, hence 
\be \H^2(E_\infty)\le \H^2(\varphi(E_\infty))\le \H^2(C_{0,p_1,p_2}).\ee

Now by Theorem \ref{unicite1}, \tb{(7.25) and (7.27)}, we know that $E_\infty=C_{0,p_1,p_2}$. But this contradicts the fact that $E_\infty$ is not $\e$-near $C_{0,p_1,p_2}$ in $\mU_0(0,1)$.
\qed

\begin{cor}For each $0<\e<\frac13 \e_0$, there exists a $\d_1=\d_1(\e)\in (0,\e)$, such that  the following holds: 

For any $k>k_0$ ($k_0$ is the one in the above proposition), suppose that $E\in\oF_\d(C_k,\mU_k(0,1))$, $E$ is minimal in $\mU(0,1)$, and is $\frac 23 \e_0$-near $C_k$ in $\mU_k(0,1)$, but not $\e$-near any translation of $C_k$ in $\mU_k(0,1)$. Then $E$ is not $\d_1$-near any translation of $C_k$ in $\mU_k(0,1)\bs \mU_k(0,\frac 14)$.
\end{cor}

\nd We prove by contradiction. Fix $k>k_0$, and set $\d_1=\frac 23\d(\e)$ ($\d(\e)$ being the one in Proposition \tb{7.6}. Suppose that there exists $p\in \R^n$ such that $E$ is $\d_1$-near $C_k+p$ in $\mU_k(0,1)\bs \mU_k(0,\frac 14)$. Since $E$ is $\frac 23 \e_0$-near $C_k$ in $\mU_k(0,1)$, we know that dist$(0,p)<C(\e_0+\d_1)$, where $C=C(K^1,K^2)$ depends only on the structure of $K^1$ and $K^2$. This means that $\mU_k(p,\frac 23)\bs \mU_k(p, \frac 13)\subset \mU_k(0,1)\bs \mU_k(0,\frac 14)$. Thus in $\mU_k(p,\frac 23)\bs \mU_k(p, \frac 13)$, $E$ is $\d_1$-near $C_k+p$. By scaling, the minimal set $\frac 32 (E-p)$ is $\frac 32\d_1=\d$ near $C_k$ in $\mU_k(0,1)\bs \mU_k(0,\frac 12)$, and is $\e_0$ near $C_k$ in $\mU_k(0,1)$. By the same argument as in Propositions \tb{6.2} and \tb{6.3}, we know that $\frac 32 (E-p)\in \oF_\d(C_k,\mU_k(0,1))$. Therefore we can apply Proposition \tb{7.6}, and get that $\frac 32 (E-p)$ is $\e$-near $C_k$ in $\mU_k(0,1)$. This means that $E$ is $\frac 23\e$ near $C_k+p$ in $\mU_k(p,\frac 32)$. Since $\mU_k(0,\frac 12)\subset \mU_k(p,\frac 32)$, hence $E$ is $\frac 23\e$ near $C_k+p$ in $\mU_k(0,\frac 12)$. By hypothesis, $E$ is already $\d_1$ near $C_k+p$ in $\mU_k(0,1)\bs \mU_k(0,\frac 14)$, hence we know that $E$ is $\e$ near $C_k+p$ in $\mU_k(0,1)$. This contradicts that fact that $E$ is not $\e$-near any translation of $C_k$ in $\mU_k(0,1)$. \qed

Now for each fixed $k$, and $\e<\frac 13\e_0$, and let $\d_1=\d_1(\e)$ be as in Corollary \tb{7.7}.

since the $\e$-process stops at step $m_k$, we know that in $\Omega_0:=\mU(o_k, r_k)$, $F_k$ is $2\e r_k(<\frac 23\e_0r_k)$ near $C_k+o_k$, but not $\e r_k$ near any translation of $C_k$. By Corollary \tb{7.7},  $F_k$ must be not $\d_1 r_k$ near any translation of $C_k$ in the annulus $\Omega:=\mU(o_k, r_k)\bs \mU(o_k, \frac 14 r_k)$. 

Let $\d_2=\min\{10^{-3}\d_1, \d(\e^8)\}$, where $\d(\e^8)$ is the $\d$ corresponds to $\e^8$ in Proposition \tb{7.6}.

As a result, there are 2 possibilities: 

1) there exists two points $p_1$ and $p_2$, such that each $F^i_{k,m}$ is $\d_2$ near $K^i_k+o_k+p_i$ in $\Omega$; 

2) one of the $F^i_{k,m}$ is not $\d_2$ near any translation of $K^i_k$ in $\Omega$. 

\begin{rem} When both $K^1$ and $K^2$ are planes, possibility 1) does not happen, because planes are invariant under translation by points in themselves. For our case, when $K^1$ or $K^2$ are singular minimal cones, the whole next section will be devoted to fix this problem of non translation invariance.
\end{rem}

Let us first treat the second possibility, and leave the possibility 1) for the next section. 

So from now on, suppose, without loss of generality, that $F^1_{k,m}$ is not $\d_2$ near any translation of $K^1_k$ in $\Omega$. 

Then there are two cases: 

$1^\circ$ $F^1_{k,m}$ is not $10^{-2}\d_2$ near any translation of $K^1_k$ in one of the region $(\mC^1_{jl,k}\times {P^1_k}^\perp+o_k)\cap \Omega$ for some $(j,l)\in J^1$, or in one of the region $(\mS^1_{j,k}\times {P^1_k}^\perp+o_k)\cap \Omega$ for some $1\le j\le \mu_1$. That is, $F^1_{k,m}$ get away from $C_k$ in its regular part.

$2^\circ$ For each $(j,l)\in J^1$, $F^1_{k,m}$ is $10^{-2}\d_2$ near some translation of $K^1_k$ in $(\mC^1_{jl,k}\times {P^1_k}^\perp+o_k)\cap \Omega$; and for each $1\le j\le \mu_1$, $F^1_{k,m}$ is $10^{-2}\d_2$ near some translation of $K^1_k$ in $(\mS^1_{j,k}\times {P^1_k}^\perp+o_k)\cap \Omega$. Then there exists $1\le j_0\le m_1$, such that $F^1_{k,m}$ is not $\d_2$ near any translation of $C_k$ in $(\mA^1_{j_0,k}\times {P^1_k}^\perp+o_k)\cap \Omega$. That is, $F^1_{k,m}$ get away from $C_k$ in its singular part, while it stays very close to $C_k$ in its regular part.

We will prove, in both cases, that the excess of measure will be of order $r_k^2$. We will treat the two cases separatedly in the following two subsections.

\subsection{Excess estimates in regular parts}

Let us first prove following theorem, which essentially gives an estimate for regular parts of $F^1_{k,m}$ outside the critical ball.

\begin{thm} For each $\d>0$, small, there exists a constant $C_1(\d)>0$, such that the following holds: 
Let $\theta_0\in [0,2\pi]$, we define $\a=\a_{\theta_0}\in \N$:  when $\theta_0>\pi$, set $\a=1$; otherwise, let $\a\ge 2$ be such that $\frac{2\pi}{2\a}<\theta_0\le \frac{2\pi}{2(\a-1)}$. Let $s>0$. Let $R=\{x=(r,\theta)\in \R^2: \frac14s\le r\le 2s, \theta\in [0,\theta_0]\}$, and $R'=\{x=(r,\theta)\in \R^2: \frac14s\le r\le s, \theta\in [0,\theta_0]\}$, where $(r,\theta)$ is the polar coordinate in $\R^2$. Let $f:R\to \R^m$ be a $C^1$ map, such that $||\nabla f||_\infty\le 1$. Denote by $G_f$ the graph of $f$: $G_f:=\{(x,f(x)):x\in R\}\subset \R^2\times\R^m$. Suppose that the restriction of the graph on $R'$ $G_f\cap (R'\times \R^m)$ of $f$ is $\d s$-far from any translation of $\R^2$, that is, for any $q\in \R^m$, 
\be\sup_{x\in R'}\mbox{dist}((x,f(x),\R^2\times \{q\})=\sup_{x\in R'}||f(x)-q||_{\R^m}\ge \d s.\ee
Then
\be\H^2(G_f)\ge \H^2(R)+C_1(\d)\a^{-1}s^2.\ee
\end{thm}

\begin{pro}Suppose $0<r_0<\frac 12$. For $\a\in \N$, let $R_\a$ be $\{x=(r,\theta)\in \R^2: r_0\le r\le 1, \theta\in [0,\frac{2\pi}{2\a}]\}$, where $(r,\theta)$ is the polar coordinate in $\R^2$, and let $L_\a$ be the arc $\{x=(r,\theta)\in R_\a: r=r_0\}$, which is  part of the boundary of $R_\a$. Let $u_0\in C^1(L_\a,\R)$. Denote by $m(u_0)=\frac{2\a}{2\pi r_0}\int_{L_\a}u_0$ its average.

Then for all $u\in C^1(R_\a,\R)$ such that
\be u|_{L_\a}=u_0,\ee
we have
\be \int_{R_\a}|\nabla u|^2\ge \frac 14r_0^{-1}\int_{L_\a}|u_0-m(u_0)|^2.\ee
\end{pro}

\nd We define $g: \overline{B(0,1)}\bs B(0,r_0)\cap \R^2\to \R$, 
\be g(r,\theta)=\left\{\begin{array}{rcl}u(r,\theta-\frac{2k}{2\a}2\pi)&,\ if&\ \theta\in [\frac{2k}{2\a}2\pi, \frac{2k+1}{2\a}2\pi], 1\le k\le \a-1;\\
u(r, \frac{2k}{2\a}2\pi-\theta)&,\ if&\ \theta\in [\frac{2k-1}{2\a}2\pi, \frac{2k}{2\a}2\pi]. 1\le k\le \a\end{array}\right.
\ee

Then $g$ is $C^1$ except on segments $\{x=(r,\theta): \theta=\frac{2k\pi}{2\a},0\le k\le 2\a-1\}$.  And by definition of $g$, we know that
\be \int_{\overline{B(0,1)}\bs B(0,r_0)}|\nabla g|^2=2\a\int_{R_\a} |\nabla u|^2.\ee

Set $v_0=g|_{\partial B(0,r_0)}$. Then it is continuous and piecewise $C^1$. By definition, it is easy to see that
\be m(v_0)=\frac{1}{2\pi r_0}\int_{\partial B(0,r_0)} u_0=m(v_0),\ee
and
\be \int_{\partial B(0,r_0)}|v_0-m(v_0)|^2=2\a\int_{L_\a}|u_0-m(u_0)|^2.\ee

Since both $v_0$ and $g$ are continuous and piecewise $C^1$, the conclusion (and the proof) of \cite{2p} Proposition 7.1 works also for $v_0$ and $g$. Thus we get
\be \int_{\overline{B(0,1)}\bs B(0,r_0)}|\nabla g|^2\ge \frac 14 r_0^{-1}\int_{\partial B(0,r_0)}|v_0-m(v_0)|^2.\ee
By \tb{(7.33), (7.35) and (7.36)}, we get
\be\int_{R_\a} |\nabla u|^2\ge \frac 14 r_0^{-1}\int_{L_\a}|u_0-m(u_0)|^2.\ee\qed

By the same reflection and periodic extension as in the proof of Proposition \tb{7.10}, we get, by Corollary 7.4 of \cite{2p}, that
\begin{pro} For all $\e\in (0,1)$, there exists $C_2=C_2(\e)>100$, such that if $0<r_0<1$, $u\in C^1(R_\a,\R)$, such that
\be u|_{L_\a}\ge \d r_0-\frac{\d r_0}{C}\mbox{ and }u|_{K_\a}<\frac{\d r_0}{C},\ee
where $K_\a=\{x=(1,\theta):0\le \theta\le \frac{2\pi}{2\a}\}$, then
\be \int_{R_\a}|\nabla u|^2\ge \e\frac{2\pi \d^2r_0^2}{2\a|\log r_0|}.\ee
\end{pro} 

\noindent\textbf{Proof of Theorem }\tb{7.9}.

Let $f_1,\cdots f_m: R\to \R$ be the coordinate functions of $u$, that is, $f=(f_1,\cdots f_m)$. Then for each $1\le j\le m$, $f_j$ is $C^1$, $||\nabla f_j||_\infty\le 1$, and there exists $j_0$ such that for any $a\in \R$, 
\be \sup_{x\in R'}||f_{j_0}-a||\ge\d=\frac{\e}{m}.\ee 

Without loss of generality, suppose that $j_0=1$. Let $G'=\{(x,f_1(x)):x\in R\}$, this is a subset of $\R^2\times \R$.

For $\frac 14s\le t\le s$, define
\be \Gamma_t=\{(x,f_1(x):x\in R, |x|=t\}\ee
the graph of $f_1$ on $\partial B(0,t)\cap R$.

Then there are two cases with respect to \tb{(7.40)}:

\textbf{Case 1:} There exists $t\in [\frac 14s,s]$ such that
\be\sup_{x,y\in \Gamma_t}|f_1(x)-f_1(y)|\ge\frac{\d s}{C},\ee
where $C=4C_2(\frac 12)$ is the constant of Proposition \tb{7.11}.

Then there exists $a,b\in \Gamma_t$ such that $|f_1(a)-f_1(b)|\ge \frac{\d s}{C}$. Note that the part of $\Gamma_t$ between $a$ and $b$ is a circle arc, let $c$ denote the midpoint of this arc. Then we know that 
\be\max\{|f_1(a)-f_1(c)|,|f_1(c)-f_1(b)|\}\ge \frac{\d s}{2C}.\ee
Suppose, without loss of generality, that $|f_1(a)-f_1(c)|\ge \frac{\d s}{2C}$.

For any $x\in \R^2\bs \{0\}$, let $\theta_x$ denote the angular coordinate. Since $0\le \theta_a,\theta_b\le\theta_0$, $[\theta_a,\theta_c]$ is an interval of length less than $\frac 12\theta_0$. We know that $\frac12\theta_0\le \frac{2\pi}{2\a}$, and hence $[\theta_a,\theta_c]$ is contained in an interval $I$ of length $\frac{2\pi}{2\a}$, and $I\subset [0,\theta_0]$. Set $R_t=\{x\in R:t\le |x|\le 1,\theta_x\in I\}$, $L_t=\Gamma_t\cap \{x:\theta_x\in I\}$, and let $u$ be the restriction of $f_1$ on $R_t$, and let $u_0$ be the restriction of $u$ on $L_t$. We let $r_0=t$, then $R_t$ is isometric to $\R_\a$ in Proposition \tb{7.10}. Hence we apply Proposition \tb{7.10} to $u$, and get that 
\be \int_{R_t}|\nabla u|^2\ge \frac 14 t^{-1}\int_{L_t}|u_0-m(u_0)|^2.\ee

Let us estimate the right-hand-side. We know that $a,c\in L_t$. Since $|u_0(a)-u_0(c)|=|f_1(a)-f_1(c)|\ge \frac{\d s}{2C}$, and $||u_0||_\infty\le ||f_1||_\infty<1$, we know that for any $\beta\in \R$, 
\be \int_{L_t}|u_0(x)-\beta|^2\ge \frac{\d^3s^3}{2^9C^3},\ee
in particular, 
\be  \frac 14 t^{-1}\int_{L_t}|u_0-m(u_0)|^2\ge\frac{\d^3}{2^{11}C^3}s^2.\ee

On the other hand, since $u$ is the restriction of $f_1$ on $R_t$, and by \tb{(7.44) and (7.46),} we have
\be \int_R|\nabla f_1|^2\ge \int_{R_t}|\nabla u|^2\ge\frac 14 t^{-1}\int_{L_t}|u_0-m(u_0)|^2\ge\frac{\d^3}{2^{11}C^3}s^2.\ee

\textbf{Case 2:} If the hypothesis of Case 1 does not hold, that is, for any $t\in [\frac 14s,2s]$,
\be\sup_{x,y\in \Gamma_t}|f_1(x)-f_1(y)|<\frac{\d s}{C}.\ee

Then after \tb{(7.28)}, there exists $\frac 14s\le t<t'\le s$ such that
\be \sup_{x\in \Gamma_t, y\in \Gamma_{t'}}|f_1(x)-f_1(y)|\ge \d s.\ee
Fix such a pair $t$ and $t'$, and without loss of generality, suppose that
\be \sup_{x\in \Gamma_t, y\in \Gamma_{t'}}f_1(x)-f_1(y)\ge \d s.\ee
Let $\beta=\inf_{x\in\Gamma_{t'}}f_1(x)$. Then by \tb{(7.48)-(7.50)}, 
\be f_1(x)-\b\le\frac{\d s}{C}, \forall x\in \Gamma_{t'}, \mbox{ and }f_1(x)-\b\ge(\d-\frac{\d }{C})s, \forall x\in \Gamma_t.\ee

 Let $r_0=\frac{t}{t'}\in[\frac 14,1)$. Let $u: R_\a\to \R: u(x)=\frac{r_0}{s}(f_1(t'x)-\b)$. Then $u$ is a dilatation of the restriction of $f_1$ on the region $R_{t,t'}=\{x:r_x\in [t,t'], \theta_x\in [0,\frac{2\pi}{2\a}]\subset R\}$.

 The function $u$ satisfies all the hypothesese in Proposition \tb{7.11} with $\e=\frac 12$ (recall that here $C=C(\frac 12)$). As a result, we know that
 \be \int_{R_\a}|\nabla u|^2\ge \frac 12\frac{2\pi \d^2r_0^2}{2\a|\log r_0|}\ge \frac{\pi}{2\a\log 4}\d^2r_0^2.\ee
 Then we know that 
 \be\int_R |\nabla f_1|^2\ge \int_{R_{t,t'}}|\nabla f_1|^2=s^2r_0^{-2}\int_{R_\a}|\nabla u|^2\ge \frac{\pi\d^2}{2\a\log 4}s^2.\ee
 
 Combining Case 1 and 2, let $C_1(\d)=3\min\{\frac{\d^3}{2^{11}C^3},\frac{\pi}{2\a\log 4}\d^2\}>0$, which only depends on $\d$. Then we have
 \be \int_R |\nabla f_1|^2\ge 3C_1(\d)\a^{-1}s^2.\ee
 
 Since $||\nabla f_1||_\infty\le 1$, we know that
 \be \begin{split}\H^2(G')&=\int_R\sqrt{1+|\nabla f_1|^2}\ge\int_R(1+\frac{|\nabla f_1|^2}{3})\\
 &=\H^2(R)+\frac13\int_R |\nabla f_1|^2\ge \H^2(R)+C_1(\d)\a^{-1}s^2\end{split}.\ee
 
To estimate the measure of $G_f$, note that $G'$ is the projection of $G_f$ on $\R^2\times\R\times \{(0,\cdots 0\}$. Hence
\be \H^2(G_f)\ge \H^2(G')\ge \H^2(R)+C(\d)s^2,\ee
which gives (7.29). \qed

Recall that $\Omega=\mU(o_k,r_k)\bs\mU(o_k,\frac 14 r_k)$. We have the following corollary:

\begin{cor} If $F^1_{k,m}$ is not $\e^8$ near any translation of $K_k^1$ in $\mC^1_{jl,k}\cap \Omega$ for some $(j_0,l_0)\in J^1$, then
\be \H^2(F^1_{k,m}\cap \Omega)\ge \H^2((C_k+o_k)\cap \Omega)+C_1(\e)r_k^2.\ee
\end{cor}

\nd The proof is similar as the proof of Proposition \tb{7.2}: we decompose $F^1_{k,m}\cap \Omega$ into regular regions $ \Omega\cap (\mC^1_{jl,k}\times {P^1_k}^\perp+o_k)$, $(j,l)\in J^1$, $ \Omega\cap (\mS^1_{j,k}\times {P^1_k}^\perp+o_k)$, $1\le j\le \mu_1$, and singular regions $ \Omega\cap (\mA^1_{j,k}\times {P^1_k}^\perp+o_k)$, $1\le j\le m_1$. 

In regular parts, for $(j,l)\in J^1\bs\{(j_0,l_0)\}$, Proposition \tb{6.1} $3^\circ$ yields that $F^1_{k,m}\cap \Omega\cap(\mC^1_{jl,k}\times {P_k^1}^\perp+o_k)$ is a graph on $(\mc^1_{jl,k}+o_k)\cap \Omega\cap(\mC^1_{jl,k}\times {P_k^1}^\perp+o_k)$, which is just $(C_k+o_k)\cap \Omega\cap(\mC^1_{jl,k}\times {P_k^1}^\perp+o_k)$. Hence
\be 
\H^2[F^1_{k,m}\cap\Omega\cap(\mC^1_{jl,k}\times {P_k^1}^\perp+o_k)]\ge\H^2[C_k+o_k)\cap \Omega\cap(\mC^1_{jl,k}\times {P_k^1}^\perp+o_k)];
\ee

Similarly, by Proposition \tb{6.1} $4^\circ$, we have, for $1\le j\le \mu_1$,
 \be 
\H^2[F^1_{k,m}\cap\Omega\cap(\mS^1_{j,k}\times {P_k^1}^\perp+o_k)]\ge\H^2[C_k+o_k)\cap \Omega\cap(\mS^1_{j,k}\times {P_k^1}^\perp+o_k)];
\ee

For $(j_0,l_0)$, we apply Theorem \tb{7.9}: let $\theta_0$ be the angle of the sector $\mC^1_{j_0l_0,k}\bs\mA^{1,k}$, $s=o_k$, then the region $\Omega':=(\mc^1_{j_0l_0,k}\bs\mA^{1,k})\cap \Omega$ is isometric to $R'$. Let $f$ be the restriction of $h^1_{j_0l_0,k}$ on $\Omega'$. Set $\d=\e^8$, then by Theorem \tb{7.9}, we know that 
\be \begin{split}
\H^2[F^1_{k,m}\cap &\Omega\cap(\Omega'\times {P_k^1}^\perp+o_k)]=\H^2(G_f)\\
&\ge \H^2(\Omega')+ C_1(\e^8)\a^{-1}r_k^2=\H^2[\mC^1_{j_0l_0,k}\cap \Omega]+ C_1(\e^8)\a^{-1}r_k^2\\
&=\H^2[(C_k+o_k)\cap \Omega\cap(\mC^1_{j_0l_0,k}\times {P_k^1}^\perp+o_k)]+ C_1(\e^8)\a^{-1}r_k^2
\end{split}\ee

Note that by the structure Theorem \tb{2.21} for 2-dimensional minimal cones, and the definition of $\mU_k$, we know that $\theta_0>\eta_0-10\eta>\frac 12\eta_0$, hence $\a<\frac{4\pi}{\eta_0}$. As a result, by \tb{(7.60)}, 
\be \begin{split}&\H^2[F^1_{k,m}\cap\Omega\cap(\Omega'\times {P_k^1}^\perp+o_k)]\\
&\ge \H^2[(C_k+o_k)\cap \Omega\cap(\mC^1_{j_0l_0,k}\times {P_k^1}^\perp+o_k)]+ C_1(\e^8)(\frac{4\pi}{\eta_0})^{-1}r_k^2\\
&=\H^2[(C_k+o_k)\cap \Omega\cap(\mC^1_{j_0l_0,k}\times {P_k^1}^\perp+o_k)]+ C_3(\e)r_k^2.
\end{split}\ee

In singular parts, for each $1\le j\le m_1$, we have again as in the proof of Proposition \tb{7.2}, by Corollary \tb{7.4} and a slicing as \tb{(7.19)} yields that 
\be\begin{split}\H^2[F^1_{k,m}\cap \Omega\cap ( (\mA^1_{j,k}\times {P^1_k}^\perp)+o_k)]\ge \H^2[(C_k+o_k)\cap \Omega\cap ( (\mA^1_{j,k}\times {P^1_k}^\perp)+o_k)].
\end{split}\ee

We sum up over all $1\le j\le m_1$, all $1\le j\le \mu_1$, and all $(j,l)\in J^1$, and get the conclusion.\qed

\subsection{Excess estimates in singular zones}

In this section we give an estimate for the case $2^\circ$ at the end of Subsection \tb{7.2}, that is: for each $(j,l)\in J^1$, $F^1_{k,m}$ is $10^{-2}\d_2$ near some translation of $K^1_k$ in $(\mC^1_{jl,k}\times {P^1_k}^\perp+o_k)\cap \Omega$, and for each $1\le j\le \mu_1$, $F^1_{k,m}$ is $10^{-2}\d_2$ near some translation of $K^1_k$ in $(\mS^1_{j,k}\times {P^1_k}^\perp+o_k)\cap \Omega$. Then there exists $1\le j_0\le m_1$, such that $F^1_{k,m}$ is not $\d_2$ near any translation of $C_k$ in $(\mA^1_{j_0,k}\times {P^1_k}^\perp+o_k)\cap \Omega$. That is, $F^1_{k,m}$ get away from $C_k$ in its singular part, while it stays very close to $C_k$ in its regular part.

\begin{pro}Under the condition of $2^\circ$ at the end of Subsection \tb{7.2}, 
\be \H^2(F^1_{k,m}\cap (\mA^1_{j_0,k}\times {P_k^1}^\perp+o_k)\cap \Omega)\ge \H^2((C_k+o_k)\cap (\mA^1_{j_0,k}\times {P_k^1}^\perp+o_k)\cap \Omega)+C_4(\e)r_k^2.\ee
\end{pro}

\nd Let $Y\subset \R^n=P_k^1\times {P^1_k}^\perp$ be the 2-dimensional $\Y$ cone such that $Y+o_k$ contains $(C_k+o_k)\cap (\mA^1_{j_0,k}\times {P_k^1}^\perp+o_k)\cap \Omega$. Then $Y\subset P_k^1$. Let $R_i,1\le i\le 3$ be the three  planes that contain the three branched half planes of $Y$. Without loss of generality, suppose that the spine of $Y$ is the first coordinate line $X_1=\{(x_1,0,\cdots, 0)\in P_k^1,x_1\in \R\}$.

We know that in $(\mA^1_{j_0,k}\times {P^1_k}^\perp+o_k)\cap \Omega$, $F^1_{k,m}$ is $2\e r_k$ near $Y+o_k$, but not $\d_2$ near any translation of $Y$. 

Recall that $\mA^1_{j_0,k}$ is the cone over the $n_1-1$-dimensional region $A^1_{j_0,k}$, which is contained in the $n_1-1$ dimensional subspace $Q_1$ of $P_k^1$ orthogonal to $X_1$. Note that $P_k^1=X_1\times Q_1$. Denote by $\pi$ the orthogonal projection from $P_k^1$ to $Q_1$. Set $A(t)=\pi(tA^1_{j_0,k})$. Then we have
\be A(t)=tA(1),\mbox{ and }tA^1_{j_0,k}=\{x_0t\}\times A(t),\ee
where $x_0$, slightly smaller than 1 (in fact $x_0=1-2\eta$ as in \tb {(5.2)}), is the first coordinate of all the points of the set $A^1_{j_0,k}$.

Also, recall that the boundary of $A^1_{j_0,k}$ contains three flat part $I^1_{j_0l,k},(j_0,l)\in J^1$. Denote by $I_j,1\le j\le 3$ the three sets $\pi(I^1_{j_0l,k},(j_0,l)\in J^1)\subset Q_1$. They are exactly the three flat parts of $\partial A(1)$. Therefore for any $t\in\R$, $tI_j,1\le i\le 3$ are the three flat part of $\partial A(t)$.

Let $V=A(\frac14 r_k)\bs A(\frac 18 r_k)\subset Q_1$, and let $I=[\frac 14 r_k x_0, r_kx_0]$. Then $I\times V\subset P_k^1$ is a sub-domain of $\mA^1_{j_0,k}$. And by a same compactness argument as in Proposition 9.6, we can prove that there exists $\d_3=\d_3(\e)>0$, that depends only on $\e$, such that $F^1_{k,m}$ is not $\d_3 r_k$ near any translation of $Y$ in $\Omega_1:=(I\times V\times {P^1_k}^\perp+o_k)\cap \Omega$.

For the rest part of $(\mA^1_{j_0,k}\times {P^1_k}^\perp+o_k)\cap \Omega$, we set
\be \Omega_2:=((\mA^1_{j_0,k}\bs (I\times A(\frac14 r_k)))\times {P^1_k}^\perp+o_k)\cap \Omega,\ee
and 
\be \Omega_3=(I\times A(\frac18 r_k))\times {P^1_k}^\perp+o_k)\cap \Omega.\ee

Then $(\mA^1_{j_0,k}\times {P^1_k}^\perp+o_k)\cap \Omega$ is the disjoint union of $\Omega_i,i=1,2,3$.

For the behavior of $F^1_{k,m}$ in each of the three $\Omega_i,1\le i\le 3$, similar to Proposition \tb{6.1}, we have the following:

\begin{lem} There exists $\e_0$ small, such that for each $\e<\e_0$:

$1^\circ$ For $j=1,2$, $F^1_{k,m}\cap \Omega_j$ coincides with the union of the graphs from $(R_i+o_k)\cap\Omega_j$ to ${R_i}^\perp$, $1\le i\le 3$, with gradiant less than 1, and infinity norm less than $2\e r_k$;

$2^\circ$ $\H^2(F^1_{k,m}\cap \Omega_j)\ge \H^2((Y+o_k)\cap \Omega_j)$, $j=1,2$;

$3^\circ$ For any $t\in I$, the intersection of $F^1_{k,m}\cap\Omega_3$ with the $n_1-1$-dimensional planar part $(\{t\}\times A(\frac 18 r_k)\times{P^1_k}^\perp)+o_k$ is a set that connects the three flat parts $(\{t\}\times\frac18 r_kI_j\times {P^1_k}^\perp+o_k, (j_0,l)\in J^1$, $1\le j\le 3$, of the boundary of $(\{t\}\times A(\frac 18 r_k)\times{P^1_k}^\perp)+o_k$.

$4^\circ$ $\H^2(F^1_{k,m}\cap \Omega_3)\ge \H^2((Y+o_k)\cap \Omega_3)$.
\end{lem}

\nd The proof of $1^\circ$ and $3^\circ$ are the same as in the proof of Proposition \tb{6.1} $3^\circ$ and $4^\circ$. $2^\circ$ follows directly from $1^\circ$. $4^\circ$ follows from $2^\circ$, Lemma \tb{7.3}, and the same argument as \tb{(7.20)}.\qed

Let us return to the proof of Proposition \tb{7.13}.

Note that $F^1_{k,m}$ is $2\e r_k$ near $Y+o_k$ in $\Omega_1$, in which $Y+o_k$ coincides with the union of the three planes $R_i,1\le i\le 3$. Note that $R_i\cap \Omega_1$ is a rectangle. Hence when $\e$ is small enough, by Lemma \tb{7.14} $1^\circ$, $F^1_{k,m}\cap \Omega_1$ coincides with three graphs $G_i,1\le i\le 3$, such that each $G_i$ is the graph of a map $g_i: R_i\cap\Omega_1\to R_i^\perp$, with $\nabla g_i\le 1$, and $||g_i||_\infty\le 2\e r_k$.

Since $F^1_{k,m}$ is not $\d_3 r_k$ near any translation of $Y$ in $\Omega_1$, we have two cases: 

\textbf{Case 1:} There exists an $1\le i\le 3$, such that the graph $G_i$ is not $10^{-2} \d_3$ near any translation of $R_i$ in $\Omega_1$. Without loss of generality, suppose that $i=1$. In this case, we can use the same argument as in the proof of Theorem \tb{7.9} and Corollary \tb{7.12} (the case for regular regions) to obtain that 
\be \H^2(G_1)\ge \H^2((R_1+o_k)\cap\Omega_1)+C_5(\e)r_k^2.\ee

While for $G_2$ and $G_3$, since they are graphs on $(R_2+o_k)\cap \Omega_1$ and $(R_3+o_k)\cap \Omega_1$, we know that 
\be \H^2(G_i)\ge \H^2((R_i+o_k)\cap\Omega_1),i=2,3.\ee
In all we have
\be \begin{split}\H^2(F^1_{k,m}\cap\Omega_1)&=\sum_{i=1}^3\H^2(G_i)\ge \sum_{i=1}^3\H^2((R_i+o_k)\cap \Omega_1)+C_5(\e)r_k^2\\
&=\H^2((Y+o_k)\cap \Omega_1)+C_5(\e)r_k^2.\end{split}\ee

Now by Lemma \tb{7.14}, we get
\be\begin{split} &\H^2(F^1_{k,m}\cap(\mA^1_{j_0,k}\times {P^1_k}^\perp+o_k)\cap \Omega)\\
=&\H^2(F^1_{k,m}\cap\Omega_1)+\H^2(F^1_{k,m}\cap\Omega_2)+\H^2(F^1_{k,m}\cap\Omega_3)\\
\ge&\H^2((Y+o_k)\cap \Omega_1)+C_5(\e)r_k^2+ \H^2((Y+o_k)\cap \Omega_2)+ \H^2((Y+o_k)\cap \Omega_3)\\
=&\H^2((Y+o_k)\cap(\mA^1_{j_0,k}\times {P^1_k}^\perp+o_k)\cap \Omega)+C_5(\e)r_k^2\\
=&\H^2((C_k+o_k)\cap(\mA^1_{j_0,k}\times {P^1_k}^\perp+o_k)\cap \Omega)+C_5(\e)r_k^2.
\end{split}\ee

\textbf{Case 2:} For each $1\le i\le 3$, the graph of $G_i$ is $10^{-2}\d_3$ near some translation of $R_i$, but there union is not $\d_3$ near any translation of $Y$ in $\Omega_1$.

This time we will do slicing in $\Omega_1$ and $\Omega_3$. Recall that $\Omega_1=(I\times V\times {P^1_k}^\perp+o_k)\cap \Omega$, and $\Omega_3=(I\times A(\frac18 r_k)\times  {P^1_k}^\perp+o_k)\cap \Omega$. We slice by $t\in I$: let $\Omega_1^t=(\{t\}\times V\times {P^1_k}^\perp+o_k)\cap \Omega$, and $\Omega_3^t=(A(\frac18 r_k)\times \{t\}\times {P^1_k}^\perp+o_k)\cap \Omega$.

Let $Y_t=(Y+o_k)\cap \partial \Omega_3^t$, which is the union of three points. The three points are just the centers of the three $\{t\}\times \frac 18 r_kI_j\times{P^1_k}^\perp+o_k,1\le j\le 3$.

We claim that 
\begin{lem} For each $t\in I$, set $F_t=F^1_{k,m}\cap \partial \Omega_3^t$. Then $F_t$ is not $\frac 12 \d_3r_k$ near any translation of $Y_t$.
\end{lem}

\nd We prove by contradiction. 

So suppose there exists $t\in I$ and some point $o\in \R^n$ such that $F_t$ is $\frac 12 \d_3r_k$ near $Y_t+o$. Since $Y_t$ is in the slice $\{t\}\times Q_1\times{P^1_k}^\perp+o_k$, hence we can suppose that $o\in \{t\}\times Q_1\times{P^1_k}^\perp+o_k$ as well.


Note that 
\be F_t= F^1_{k,m}\cap\partial \Omega_1^t=\cup_{i=1}^3 G_i\cap\partial \Omega_1^t,\ee
and each $G_i$ intersects $\cap \partial \Omega_3^t$ at exactly on  point $y_i$, we know that $F_t= \{y_i,1\le i\le 3\}$. 

Recall that $G_i$ is the graph of $g_i$ over $(R_i+o_k)\cap \Omega_1$, suppose that $y_i=(z_i, g_i(z_i)), z_i\in (R_i+o_k)\cap \Omega_1$. Then the fact that $F_t= \{y_i,1\le i\le 3\}$ is $\frac 12 \d_3r_k$ near $Y_t+o$ means that $|g_i(z_i)-\pi_{R_i^\perp}(o)|<\frac 12 \d_3r_k$, where $\pi_{R_i^\perp}$ is the orthogonal projection to $R_i^\perp$. 

Recall that we are in Case 2, the graph of $G_i$ is $10^{-2}\d_3$ near some translation of $R_i$. Therefore, for each $z\in R_i\cap \Omega_1$, $|g_i(z)-\pi_{R_i^\perp}(o)|<\frac 34 \d_3r_k$. That is, for each $z\in R_i\cap \Omega_1$, $\pi_{R_i^\perp}((z,g_i(z))-o)<\frac 34 \d_3r_k$, which implies that dist$((z,g_i(z)),R_i+o)<\frac 34 \d_3r_k$. As a result, $G_i\subset B((R_i+o)\cap\Omega_1,\frac 34 \d_3r_k)$, $1\le i\le 3$. 

On the other hand, since each $G_i$ is a almost flat graph on $R_i$, hence the above inclusion results in that each $G_i$ is $\d_3 r_k$ near $R_i+o$ in $\Omega_1$, $1\le i\le 3$. As their union, we know that $ F^1_{k,m}$ is $\d_3 r_k$ near $Y+o_k+o$ in $\Omega_1$, this contradicts our hypothesis.\qed

Now we are going to estimate the length of each slice $F^1_{k,m}\cap\Omega_3^t$ for each $t\in I$.

Fix any $t$. Denote by $I_1,I_2,I_3$ the three flat parts $(\frac18 r_kI_j\times \{t\}\times{P^1_k}^\perp+o_k$, $1\le j\le 3$, of the boundary of $\Omega_3^t$, for short. The center of each $I_j$ being the intersection of $(R_j+o_k)\cap \partial \Omega_3^t$.

 Let $H=F^1_{k,m}\cap\Omega_3^t$. By Lemma \tb{7.14}, $H$ connects the three flat parts $I_j$, $1\le j\le 3$, of the boundary of $\Omega_3^t$. On the other hand, denote by $a_j$ the center of $I_j$, $1\le j\le 3$.

Note that $H\cap \partial\Omega_3^t$ is just the $F_t=\{y_j,1\le j\le 3\}$ in Lemma \tb{7.14}, each $y_j\in G_j\cap I_j$. Hence $H$ connects the three points $y_j,1\le j\le 3$, while $F_t=\{y_j,1\le j\le 3\}$ is not $\d_3$ near any translation of $Y_t=\{a_1,a_2,a_3\}$. Denote by $o_t$ is the center of $\Omega_3^t$. Then
$(Y+o_k)\cap\Omega_3^t$ is the almost disjoint union of the three segments $[o_t,a_j],1\le j\le 3$.

we are going to prove that 
\be \H^1(H)\ge \H^1((Y+o_k)\cap \Omega_3^t.\ee

Let $q\in \Omega_3^t$ be the point such that the angle between the segments $[q,y_j],1\le j\le 3$ are $120^\circ$, and this implies that 
\be\H^1(H)\ge\sum_{j=1}^3 |q-y_j|.\ee

Such a point $q$ exists and belongs to $B(o_t, 4\e r_k)$, because $\{y_j,1\le j\le 3\}$ are $2\e r_k$ near $Y_t$. Moreover by the same reason, we know that for each $j$, there exists $z_j\in I_j$ such that the segment $[q,z_j]$ is orthogonal to the planar part $I_j$. 

As a result, by the same argument as in Lemma \tb{7.3}, we know that
\be \sum_{j=1}^3 |q-z_j|=\sum_{j=1}^3 
|o_t-a_j|=\H^1((Y+o_k)\cap\Omega_3^t).\ee

Note that the set $\cup_{j=1}^3 \{z_j\}=Y_t-o_k+q$ is a translation of $Y_t$, hence it is not $\d_3 r_k$ near $\{y_j,1\le j\le 3\}$. As a result, there exists $1\le j\le 3$ such that $d(y_j,z_j)>\d_3 r_k$. Suppose, for example, that $d(y_1,z_1)>\d_3 r_k$. Then we have
\be \begin{split}|q-y_1|&=\sqrt{|q-z_1|^2+|y_1-z_1|^2}=|q-z_1|\sqrt{1+(\frac{|y_1-z_1|}{|q-z_1|})^2}\\
&\ge |q-z_1|(1+\frac 14(\frac{|y_1-z_1|}{|q-z_1|})^2)=|q-z_1|+\frac 14\frac{|y_1-z_1|^2}{|q-z_1|}
\end{split}\ee
The inequality is because $|y_1-z_1|\le $diam$I_j<10^{-2}|o_t-a_1|$, and $|q-z_1|\ge |o_t-a_1|-|o_t-q|-|a_1-z_1|\ge |o_t-a_1|-4\e r_k-$diam$I_j\ge \frac 12|o_t-a_1|$. Hence $\frac{|y_1-z_1|}{|q-z_1|}<1$.

Again we know that $|q-z_1|<$diam$\Omega_3^t\le \frac14 r_k$, and $|y_1-z_1|\ge \d_3 r_k$, hence
\be |q-y_1|\ge |q-z_1|+\d_3^2r_k.\ee

For $j=2,3$, since $[q,z_j]$ is perpendicular to $I_j$, we know that $|q-y_j|\ge |q-z_j|$. As a result, we have
\be \H^1(H)\ge\sum_{j=1}^3 |q-y_j|\ge|q-z_1|+\d_3^2r_k+\sum_{j=2}^3|q-z_j|=\H^1((Y+o_k)\cap\Omega_3^t)+\d_3^2r_k.\ee

This gives, that for each $t\in I$, we have
\be \H^1(F^1_{k,m}\cap\Omega_3^t)\ge \H^1((Y+o_k)\cap\Omega_3^t)+\d_3^2r_k.\ee
Therefore
\be \begin{split}\H^2(F^1_{k,m}\cap\Omega_3)&=\int_{t\in I}\H^1(F^1_{k,m}\cap\Omega_3^t)\ge\int_{t\in I}\H^1((Y+o_k)\cap\Omega_3^t)+\d_3^2r_k\\
&=|I|\times (\H^1((Y+o_k)\cap\Omega_3^t)+\d_3^2r_k)
=\H^2((Y+o_k)\cap\Omega_3)+|I|\times \d_3^2r_k\\
&=\H^2((Y+o_k)\cap\Omega_3)+\frac34 \d_3^2r_k^2.
\end{split}\ee

Together with Lemma \tb{7.14} $2^\circ$, we get
\be\begin{split} &\H^2(F^1_{k,m}\cap(\mA^1_{j_0,k}\times {P^1_k}^\perp+o_k)\cap \Omega)\\
=&\H^2(F^1_{k,m}\cap\Omega_1)+\H^2(F^1_{k,m}\cap\Omega_2)+\H^2(F^1_{k,m}\cap\Omega_3)\\
\ge&\H^2((Y+o_k)\cap \Omega_1)+ \H^2((Y+o_k)\cap \Omega_2)+ \H^2((Y+o_k)\cap \Omega_3)+\frac34 \d_3^2r_k^2\\
=&\H^2((Y+o_k)\cap(\mA^1_{j_0,k}\times {P^1_k}^\perp+o_k)\cap \Omega)+\frac34 \d_3^2r_k^2\\
=&\H^2((C_k+o_k)\cap(\mA^1_{j_0,k}\times {P^1_k}^\perp+o_k)\cap \Omega)+\frac34 \d_3^2r_k^2.
\end{split}\ee

Now combining \textbf{Case 1} and \textbf{Case 2}, we set $C_4(\e)=\min\{C_5(\e),\frac34 \d_3^2r_k^2\}$ (recall that $\d_3$ depends only on $\e$), and get Proposition \tb{7.13}.\qed

\begin{cor}If $F^1_{k,m}$ is not $\d_2$ near any translation of $K_k^1$ in $\mA^1_{j_0,k}\cap \Omega$ for some $1\le j_0\le m_1$, then
\be \H^2(F^1_{k,m}\cap \Omega)\ge \H^2((C_k+o_k)\cap \Omega)+C_4(\e)r_k^2.\ee
\end{cor}

\nd The result comes directly from Proposition \tb{7.13}, and the same argument as Corollary \tb{7.12}.\qed

To summerize, by Corollary \tb{7.12} and Corollary \tb{7.16} (for case $1^\circ$ and $2^\circ$ at the end of Subsection 7.2 respectively), when the possibility (2) at the end of Subsection 7.2 happens, we have that
\be \H^2(F^1_{k,m}\cap \Omega)\ge \H^2((C_k+o_k)\cap \Omega)+C_6(\e)r_k^2,\ee
where $C_6(\e)=\min\{C_3(\e), C_4(\e)\}$ is a constant that only depends on $\e$.

\section{An argument of reduction}

Now let us turn to Possibility 1) stated after Corollary \tb{7.7}. That is, in $\Omega_0=\mU(o_k,r_k)$, $F_k$ is not near any translation of the union $K^1_k\cup K^2_k$, but it is near the union of a translated $K^1_k$, and a translated $K^2_k$. In other words, we cannot approximate $F_k$ by the union of translated $K^1_k$ and $K^2_k$ with the same center, but we can approximate $F_k$ if we translate $K^1_k$ and $K^2_k$ separately to different centers.

To treat this case, we are going to do the following: suppose we can approximate $F_k$ by the union of a translated $K^1_k+p_1$ and a translated $K^2_k+p_2$ with different centers $p_1$ and $p_2$. Then the point $o$ of "intersection" of these two cones cannot be near both $p_1$ and $p_2$. Suppose for example $o$ is far from $p_1$. As a result, near the point $o$ there is no point of type $K^1$ for $K^1_k+p_1$. 

By the structure Theorem \tb{2.21} for 2-dimensional minimal cones, in a neighborhood of $o$, all points of $K^1_k+p_1$ must be $\Y$ or $\P$ points. As a result, in the region near this intersection point $o$, $F_k$ will be near a minimal cone $C_k'$ which is the union of a $\Y$ cone or a plane, with $K_k^2$ (when $o$ is near $p_2$) or a $\Y$ cone or a plane (when $o$ is far from $p_2$ as well). Moreover, we will prove that $F_k$ is a limit of deformation of this cone $C'_k$. Note that both $\Y$ cones and planes are Almgren unique (cf. \cite{uniquePYT}) and sliding stable (cf. \cite{stablePYT}). Hence we can redo the argument that starts from Sections 6, to get an estimate. Of course, for the rest part of $F_k$ far from $o$, we can also prove that the measure is no less than the measure of $C_k'$.

Note that the cone $C_k'$ is ''simpler'' than $C_k$, because its first component is strictly simpler than $K^1$, and the second component is simpler than (i.e. strictly simpler than or equal to) $K^2$. Here the word ''simpler'' means the following:

A plane is strictly simpler than any singular 2-dimensional minimal cone, and $\Y$ set is strictly simpler than any singular 2-dimensional minimal cone other than $\Y$. This is natural, because any singular 2-dimensional minimal cone admit a plane as a blow-up limit at some point, but the inverse does not hold; and by the structure Theorem \tb{2.21},  any singular 2-dimensional minimal cone other than $\Y$ admit $\Y$ as a blow-up limit at some point, but the inverse does not hold.

Next we will do the whole discuss (Sections \tb{5 to 7}) again, but with respect to a minimal cone $C_k'$ which is strictly simpler than $C_k$. And if Possibility 1) happens again, we can reduce to an even simpler minimal cone, in the sense that it is a union of two minimal cones as well, each component of the union is simpler than that of $C_k'$, and at least one of the component is strictly simpler than that of $C_k'$. We can continue this reduction, and will either ends by possibility 2) at some step, or ends by possibility 1), and reduced finally to a union of two planes. When we come to the union of two planes, since for planes every point has the same type of blow-up limit, possibility 1) cannot happen, and again our argument ends by possibility 2).

So let us do it more precisely.

\subsection{The new critical region}

Recall that $F_k$ is $2\e r_k$ near $C_k+o_k$, but not $\e r_k$ near any translation of $C_k$ in $\Omega_0=\mU(o_k,r_k)$; not $\d_1r_k$ near any translation of $C_k$ in the annulus $\Omega=\mU(o_k,r_k)\bs \mU(o_k, \frac 14 r_k)$. But there exists two points $p_1$ and $p_2$ such that $F^i_{k,m}$ is $\d_2r_k$ near $K^i_k+o_k+p_i$ in the annulus $\Omega$. Since $F_k$ is $2\e r_k$ near $C_k+o_k$ in $\mU(o_k, r_k)$, we know that $p_1,p_2$ must belong to $\bar B(o_k, \frac 14 r_k)$, and $d(p_1,p_2)<3\e r_k$.

Now since $F^i_{k,m}$ is $\d_2r_k$ near $K^i_k+o_k+p_i$ in $\Omega$, we know that $F_k$ is $\d_2 r_k$ near $C_{k,p_1,p_2}+o_k$ in $\Omega$. Recall that $\d_2<\d(\e^8)$, hence by Proposition \tb{7.6}, we know that $F_k$ is $\e^8 r_k$ near $C_{k,p_1,p_2}+o_k$ in $\Omega_0$. 

Recall that $F_k$ is not $\e r_k$ near any translation of $C_k$ in $\Omega_0$, hence we have that $ d(p_1,p_2)>\frac12 \e r_k.$

If the distance between the two $K^i_k+o_k+p_i,i=1,2$ is larger than $2\times \e^8 r_k$, then since $F_k$ is $10^{-3}\e r_k$ near $C_{k,p_1,p_2}+o_k$ in $\Omega_0$, we can decompose $F_k\cap \Omega_0$ into a disjoint union $E_1$ and $E_2$, each $E_i$ is $\e^8 r_k$ near $K^i_k+o_k+p_i$ in $\Omega_0$. But by Proposition \tb{6.3}, $F_k\cap \Omega_0\in\oF_{4\e}((C_k+o_k)\cap \bar\Omega_0,\bar\Omega_0)$, it is a limit of $4\e$-sliding deformation of $C_k+o_k$ in $\Omega_0$. In particular, $F_k\cap \Omega_0$ is connected, and hence cannot be decomposed into disjoint union of two closed subset.

Hence the distance between $K^i_k+o_k+p_i,i=1,2$ is smaller than $2\times \e^8 r_k$. Take $x_i\in K^i_k+o_k+p_i$, $i=1,2$, such that $d(x_1,x_2)<3\times \e^8 r_k$. Set $p_2'=p_2+x_1-x_2$, then we know that $(K^1_k+o_k+p_1)\cap (K^2_k+o_k+p_2')=x_1$, and $F_k$ is $4\e^8 r_k<\e^7 r_k$ near $C_{k,p_1,p_2'}+o_k$ in $\Omega_0$.

So in all the following context, to save notations, we write $p_2$ for $p_2'$ instead. Then we have:

\be F_k\mbox{ is }2\e r_k\mbox{ near }C_k+o_k\mbox{ in }\Omega_0;\ee
\be F_k\mbox{ is not }\e r_k\mbox{ near any translation of }C_k\mbox{ in }\Omega_0;\ee
\be F_k\mbox{ is }\e^8 r_k\mbox{ near }C_{k,p_1,p_2}+o_k\mbox{ in }\Omega_0.\ee

As a result of the above three equations, we have
\be \frac 12 \e r_k<d(o_k+p_1,o_k+p_2)<3\e r_k.\ee

Now let us study the situation near the intersection point $x_1$. Recall that by structure of 2-dimensional minimal cones, we have, for each $i=1,2$: except for $o_k+p_i$, all other points of $K^i_k+o_k+p_i$ are of type $\Y$ of $\P$. Obviously $x_1$ cannot be very close to both $o_k+p_i$. Then have the following cases:

\medskip

\noindent \textbf{Case 1:} Suppose that $x_1$ is $\e^3 r_k$ near one of the $o_k+p_i$, say, $o_k+p_1$. Then by \tb{(8.4)}, $d(x_1,o_k+p_2)>\frac 13 \e r_k$. Let $y\in K^2_k+o_k+p_2$ be the nearest $\Y$ point near $x_0$, then we have the following cases:

\textbf{Case 1-1:} if $d(x_1,y)>\e^2 r_k$, then $x_1$ must be a $\P$ point of $K^2_k+o_k+p_2$. Let $P\subset P^2_k$ be the linear tangent plane of $K^2_k+o_k+p_2$ at $x_1$. Set $C_k'=K^1_k\cup P-$, set $t_k=\e^2 r_k$. Then in $\mU(x_1,t_k)$, we know that $F_k$ is $\e t_k$ near $C_k'+x_1$;

\textbf{Case 1-2:} if $d(x_1,y)\le \e^2 r_k$, then let $Y\subset P^2_k$ be the $\Y$ type tangent cone of $K^2_k+o_k+p_2$ at $y$. Set $C_k'=K^1_k\cup Y$, set $t_k=\e r_k$. Then in $\mU(x_1, t_k)$, we know that $F_k$ is $\e t_k$ near $C_k'+x_1$;

\bigskip

\noindent \textbf{Case 2:} Suppose that $x_1$ is $\e^3 r_k$ far from both $o_k+p_i,i=1,2$. For each $i=1,2$, let $y_i\in K^i_k+o_k+p_i$ be the nearest $\Y$ point near $x_0$, and let $d_i=d(x_1,y_i)$. Without loss of generality, suppose that $d_1>d_2$. Then we have the following cases:

\textbf{Case 2-1:} if $d_2>\e^6 r_k$. Then we must have $d_1>\e^6$. For $i=1,2$, let $P_i\subset P^i_k$ be the tangent plane of $K^i_k+o_k+p_i$ at $x_1$. Set $C_k'=P_1\cup P_2$. Set $t_k=\e^6 r_k$. Then $F_k$ is $\e t_k$ near $C_k'+x_1$ in $\mU(x_1,t_k)$;

\textbf{Case 2-2:} if $d_2\le \e^6 r_k$, and $d_1\le \e^4 r_k$. For $i=1,2$, let $Y_i\subset P^i_k$ be the $\Y$ type tangent cone of $K^i_k+o_k+p_i$ at $y_i$. Set $C_k'=Y_1\cup Y_2$. Set $t_k=\e^6 r_k$. Then $F_k$ is $\e t_k$ near $C_k'+x_1$ in $\mU(x_1,t_k)$;

\textbf{Case 2-3:} if neither Case 2-1 or Case 2-2 happen. Then we must have $d_1>\e^4 r_k$ and $d_2\le \e^6 r_k$. Let $P\subset P^1_k$ be the tangent plane of $K^1_k+o_k+p_1$ at $x_1$, and let $Y\subset P^2_k$ be the $\Y$ type tangent cone of $K^2_k+o_k+p_2$ at $y_2$. Set $C_k'=P\cup Y$. Set $t_k=\e^4 r_k$, then $F_k$ is $\e t_k$ near $C_k'+x_1$ in $\mU(x_1,t_k)$.

Note that in all the above cases, $C_k'$ is a union of two minimal cones $A_1$ and $A_2$, the angle between which is the same as the angle between $K^1_k$ and $K^2_k$. Thus $C_k'$ is a minimal cone, and is strictly simpler than $C_k$, in the sense that each component $A_i$ is simpler than $K^i$, and at least one of the $A_i$ is not $K^i$. Also, by Theorem \tb{2.21}, we know that $A_i$ is one of $K^i$, $\P$ and $\Y$, all of which are sliding stable, and Almgren unique.

After the above argument, we have found $x_1\in \mU(o_k,\frac 19 r_k)$, and $t_k\in (\e^6 r_k,\e r_k)$, such that $F_k$ is $\e$ near a minimal cone $C_k'$ in $\mU(x_1,t_k)$, where $C_k'$ is a union of two sliding stable Almgren unique minimal cones $A_1$ and $A_2$, with $A_1$ being $K^1$, a $\Y$ cone or a $\P$ cones, and $A_2$ being a $\Y$ cone or a $\P$ cone. 

Recall that our goal is to compare the measure of $F_k$ and $C_k+o_k$ inside $\mU(o_k,r_k)$. But by stability of measure \cite{stablePYT} for each $K^i$, we know that $\H^2((C_k+o_k)\cap \mU(o_k,r_k))=\H^2((C_{k,p_1,p_2}+o_k)\cap\mU(o_k,r_k))$, hence we only have to compare the measure of $F_k$ with $C_{k,p_1,p_2}+o_k$ in $\mU(o_k,r_k)$.

Let $\mU'_k$ denote the convex domain according to the cone $C_k'$. Then to restart our argument as before for $C_k$, we have to prove two things:

$1^\circ$ The measure of $F_k\cap \mU(o_k,r_k)\bs \mU'(x_1,t_k)$ is at least that of $C_{k,p_1,p_2}+o_k\cap\mU(o_k,r_k)\bs \mU'(x_1,t_k)$;

$2^\circ$ The part $F_k\cap \bar\mU'(x_1,t_k)$ contains an element in $\oF_{\e t_k}((C_k'+x_1)\cap \bar\mU'(x_1,t_k),\bar\mU'(x_1,t_k))$, which is $10^{-2}\e t_k$ near $(C_k'+x_1)$ in  $\mU'(x_1,t_k)$.

The proof of $1^\circ$ is similar to the proof of Proposition \tb{7.2} and part of the argument in Proposition \tb{6.1}: since we are outside the intersection area of the two minimal cones $K^1_k+o_k+p_1$ and $K^2_k+o_k+p_2$, we can decompose $F_k\cap  \mU(o_k,r_k)\bs \mU'(x_1,t_k)$ into two disjoint parts $F'_1$ and $F'_2$ that are close to $K^1_k+o_k+p_1$ and $K^2_k+o_k+p_2$ respectively. Take $F'_1$ for example, if Case 1 happens, then we use the exactly same argument as Proposition \tb{6.1} and Proposition \tb{7.2}, and get that 
\be \H^2(F_1')\ge \H^2((K^1_k+o_k+p_1)\cap \mU(o_k,r_k)\bs \mU'(x_1,t_k));\ee
otherwise, if Case 2 happen, then we have to decompose again $F_1'$ into two parts: a region with the shape $\mU$ centered at $o_k+p_1$, and the rest. Again the same argument as in Proposition \tb{6.1} and Proposition \tb{7.2}, gives \tb{(8.5)}.

Similarly we get
\be \H^2(F_2')\ge \H^2((K^2_k+o_k+p_2)\cap \mU(o_k,r_k)\bs \mU'(x_1,t_k)).\ee

Altogether we have $1^\circ$. 

We are then going to prove $2^\circ$ in the next subsection.

\subsection{Deformation}

In this section we would like to prove that $F_k\cap \bar\mU'(x_1,t_k)\in\oF_{\e t_k}((C_k'+x_1)\cap \bar\mU'(x_1,t_k),\bar\mU'(x_1,t_k))$. 

We know that $F_k\cap \bar\mU(o_k,r_k)\in \oF_{4\e r_k}((C_k+o_k)\cap \bar\mU(o_k,r_k), \bar\mU(o_k,r_k))$. It is easy to see that $(C_k+o_k)\cap \bar\mU(o_k,r_k)\in \oF_{4\e r_k}((C_{k,p_1,p_2}+o_k)\cap \bar\mU(o_k,r_k), \bar\mU(o_k,r_k))$. Hence we have
\be F_k\cap \bar\mU(o_k,r_k)\in \oF_{4\e r_k}((C_{k,p_1,p_2}+o_k)\cap \bar\mU(o_k,r_k), \bar\mU(o_k,r_k)).\ee

On the other hand, since $F_k$ is $\e^8 r_k$ near $(C_{k,p_1,p_2}+o_k)$ in $\bar\mU(o_k,r_k)$, modulo taking the limit of some reparametrization of $F_k$ by $(C_{k,p_1,p_2}+o_k)$ near $\partial\mU(o_k,r_k)$, we know that in fact, 
\be F_k\cap \bar\mU(o_k,r_k)\in \oF_{\e^8 r_k}((C_{k,p_1,p_2}+o_k)\cap \bar\mU(o_k,r_k), \bar\mU(o_k,r_k)).\ee

Note that $(C_k'+x_1)\cap \bar\mU'(x_1,t_k)\subset (C_{k,p_1,p_2}+o_k)\cap \bar\mU(o_k,r_k)$.

\begin{pro}Let $K\subset \R^d$ be a 2-dimensional minimal cone, and let $x\in K$. Suppose that $r>0$ is such that $\bar B(x,r)\cap K$ coincides with a minimal cone $A+x$ centered at $x$. Then there exists a Lipschitz map $\varphi: \bar B(x,r)\cap (A+x)\to K\cap \bar B(0,1)$, such that 

$1^\circ$ $\varphi(x)=x$, and $\varphi((A+x)\cap \partial B(x,r))\subset K\cap \partial B(0,1)$;

$2^\circ$ $\varphi: \bar B(x,r)\cap (A+x)\to \varphi(\bar B(x,r)\cap (A+x))$ is bi-Lipschitz;

$3^\circ$ There exists a Lipschitz deformation retract $\psi: K\cap \bar B(0,1)\to \varphi(\bar B(x,r)\cap (A+x))$, such that $\varphi^{-1}\circ\psi\circ \varphi ((A+x)\cap \partial B(x,r))=(A+x)\cap \partial B(x,r)$.
\end{pro}

\nd 
We have three cases: 

\textbf{Case 1:} $x=0$. In this case, $A$ and $K$ are the same cone, and we can just set $\varphi$ to be the dilatation $\varphi(y)=\frac{1}{r}x$, and set $\psi$ to be the identity map;

\textbf{Case 2:} $x\ne 0$ is a $\Y$ type point of $K$. In this case, we know that the $\Y$ line $[0,x]\cap \bar B(x,r))\subset (A+x)_Y$, and since $A+x$ is a $\Y$ cone centered at $x$, we know that the whole half line issued from $x$ and containing $[0,x]\cap \bar B(x,r)$ is contained in the spine $(A+x)_Y$ of $A+x$. Since a 2-dimensional $\Y$ set is invariant under translation along its spine, we know that 
\be A+x=A.\ee

Set $x'=\frac{x}{|x|}\in K\cap \partial B(0,1)$, then it is a $\Y$ point of $K\cap \partial B(0,1)$. By the structure Theorem \tb{2.21} of 2-dimensional minimal cones, $X=K\cap \partial B(0,1)$ is a union of finitely many great circles $S_1,\cdots, S_m$, and a net composed of finitely many arcs of great circles $L_j,j\in J$,  that can only meet at their endpoints, and each endpoint is a common endpoint of exactly three of these arcs, which meet by $120^\circ$ angles. 

Let $X_0$ be the union of $L_j,j\in J$, and let $K_0$ be the cone over $X_0$. Then $K_0\subset K$. 

Let $L_1,L_2,L_3$ be the three arcs of great circles that meet at $x'$, and let $z_i$ be the other end point of $L_i,i=1,2,3$ (it may happen that $z_1=z_2=z_3$, and in this case $K_0=A$, and we turn to the above \textbf{Case 1}). Then in $\bar B(x,r)$, $K_0$ coincides with the cone over $\cup_{i=1}^3L_i$.

When the three $z_1,z_2,z_3$ are distinct, let $L_4,L_5$ be the two other arcs of great circles with endpoint $z_1$. Let $z_4,z_5$ be the other endpoints of $L_4,L_5$ respectively. Note that $z_i,1\le i\le 3$ are distinct, but $z_4,z_5$ may coincide with $z_i,1\le i\le 3$.

\centerline{\includegraphics[width=0.5\textwidth]{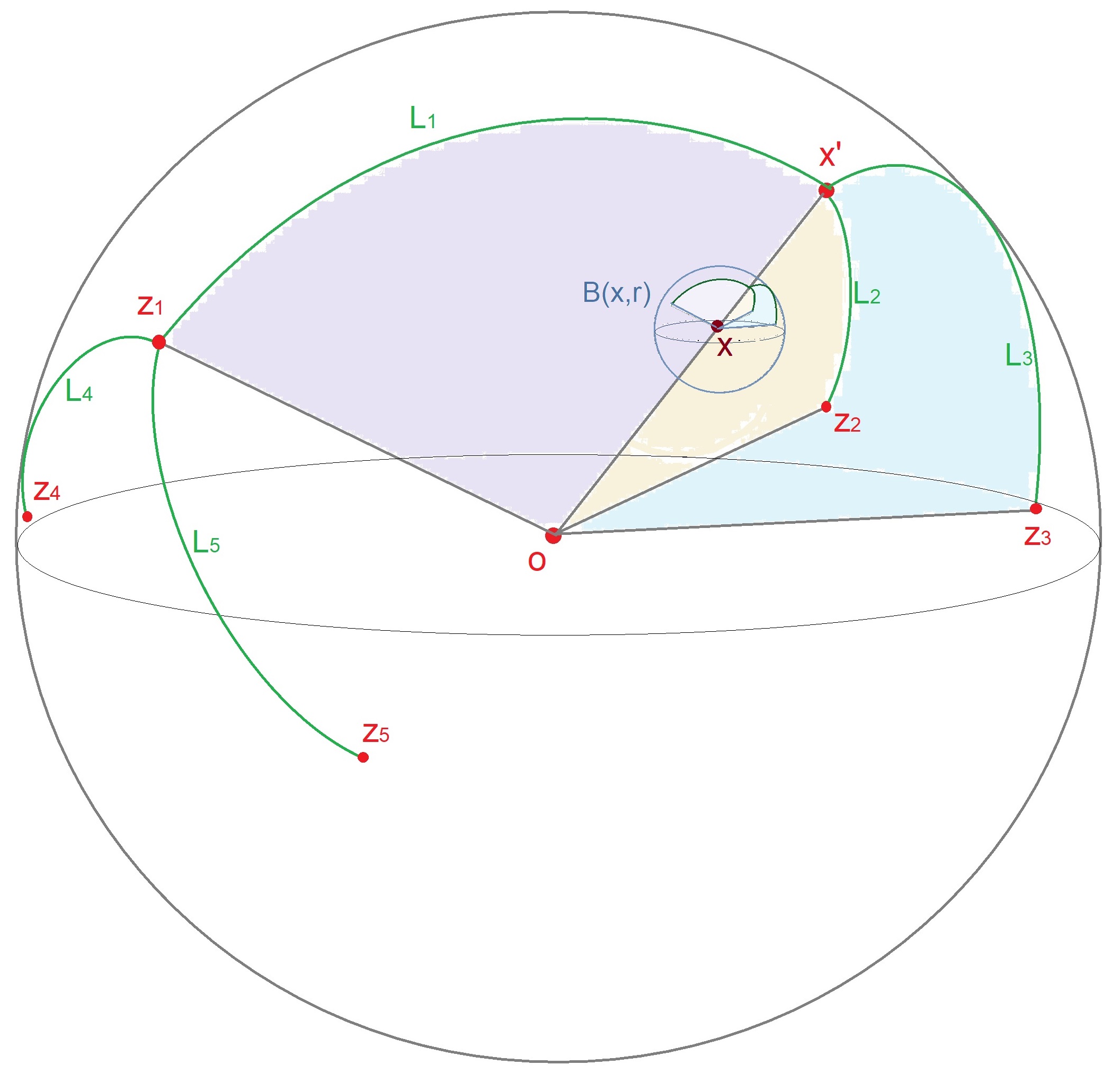} }

Now we are going to prove that 

\begin{lem} There exists a subset $J'\in J$, such that $\{1,2,3\}\subset J'$, and except for $x'$ and $z_1$, each of the other endpoints of $L_j,j\in J'$ is the common endpoint of exactly 2 of them. In other word, if we denote by $A_0$ the intersection of the $\Y$ cone $A$ with the unit sphere, then the union $X_1:=\cup_{j\in J'}L_j$ is a bi-Lipschitz version of $A_0$.\end{lem}

\nd Let $z_j,1\le j\le l$ be the set of $\Y$ points of $X_0$. Then $l$ must be even. 

\textbf{Case 1:} $l=2$, then $X_0$ must be a $\Y$ cone, and the conclusion is trivial, with $J'=J$;

\textbf{Case 2:} $l=4$, then $X_0$ is topologically the complete graph $K_4$, and we set $J'=\{1,2,3,4,5\}$ (this is equivalent to deleting an edge of $K_4$, which gives a topological version of $A_0$.

In case $l\ge 6$, we can do the following: since $l\ge 6$, take a $\Y$ point $z_6$ of $X_0$ other than $z_i,1\le i\le 5$. Let $j_1,j_2,j_3\in J$ be such that $z_6$ is an endpoint of $L_{j_1}, L_{j_2},L_{j_3}$. We delete $L_{j_1}, L_{j_2},L_{j_3}$, and regard the other three endpoints of these three $L_j$ as regular points, because now these points are the endpoints of 2 elements in $\{L_j:j\in J\bs \{j_1,j_2,j_3\}\}$. Now the new net is topologically a graph whose edges can only meet by group of three, and with $l-4$ endpoints.

\centerline{\includegraphics[width=0.8\textwidth]{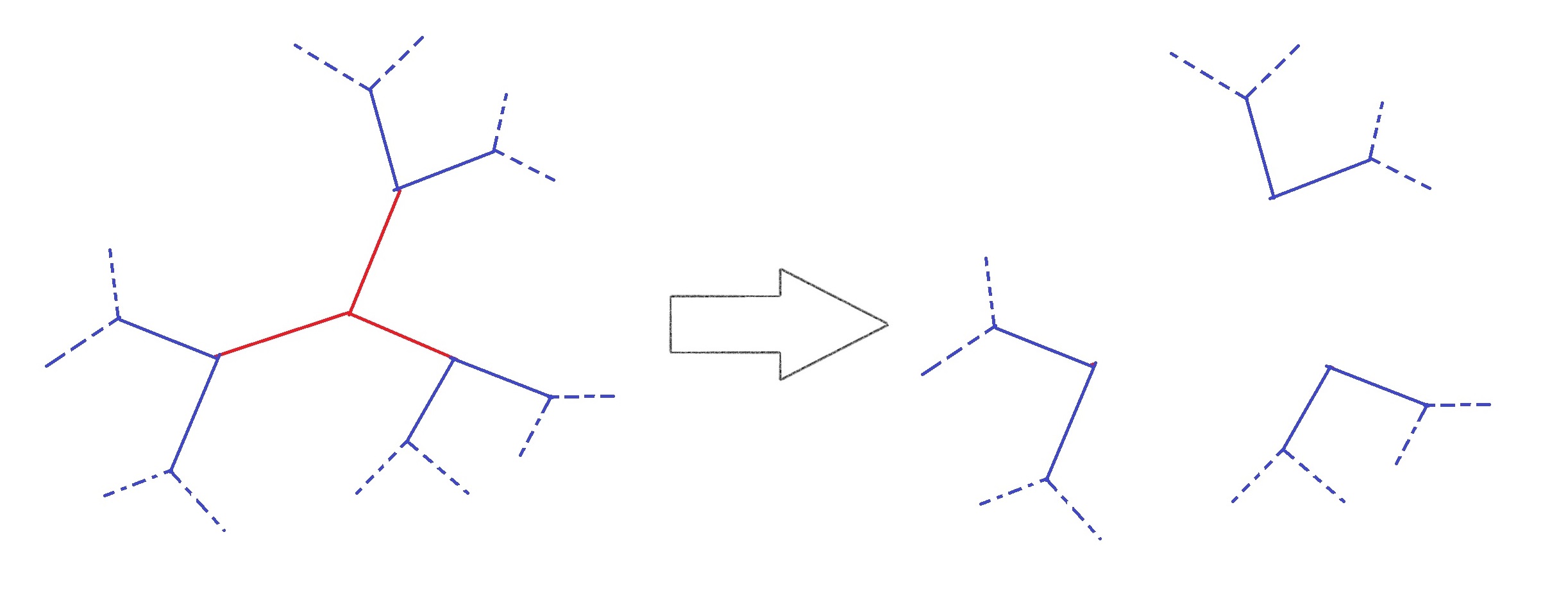} }

If $l-4\le 4$ we turn into case 1 or 2, and finish; otherwise, we continue the above reduction. This reduction should end at a finite step, and we get the conclusion.
\qed

Now let us return to the proof of Proposition \tb{8.1}. By Lemma \tb{8.2}, $X_1$ is a bi-Lipschitz version of $A_0$, with $x'$ and $z_1$ being the two $\Y$ type points. Moreover, $X_1\cap A_0=\cup_{j=1}^3 L_j$. Let $A^r$ be the cone over $A\cap \bar B(x,r)$, and let $A^r_0=A^r\cap \partial B(0,1)$. Then $A^r_0\subset \cup_{j=1}^3 L_j\bs\{z_1,z_2,z_3\}$. This means that $A^r_0\subset X_1\cap A_0$, and it only contains the common $\Y$ points $x'$ of each set, but does not touch the other $\Y$ point of neither $X_1$ nor $A_0$.

As a result, we can find a bi-Lipschitz map $f_0: A_0\to X_1$, such that $f_0=id$ on $A^r_0$, $f_0(-x')=z_1$ ($-x'$ and $z_1$ being the $\Y$ points of $A_0$ and $X_1$ respectively). Then we have $f_0(x')=x'$.

Let $K_1$ denote the cone over $X_1$. Then we can extend $f_0$ radially to a bi-Lipschitz map $f$ between the cones $A\cap \bar B(0,1)$ and $K_1\cap\bar B(0,1)$. That is, for each $z\in A\cap \bar B(0,1)$, set $f(z)=|z|f_0(\frac{z}{|z|})$. Note that $f(x)=x$.

Next, define $g: A\cap \bar B(x,r)\to A\cap \bar B(0,1)$ to be the bi-Lipschitz map as following: for each $\theta\in S^{n-1}$, let $z_\theta$ be the intersection of the ray $\{x+t\theta,t\ge 0\}$ with $\partial B(0,1)$, and let $t_\theta\in (r, 2-r )$ be such that $z_\theta=x+t_\theta \theta$; then for each $r'\in [0,r]$, set $g(x+r'\theta)=x+\frac{r't_\theta}{r}\theta$. Then we have $g(x)=x$.

 Then set $\varphi=f\circ g:  A\cap \bar B(x,r)\to K_1\cap\bar B(0,1)$. Then $\varphi$ is a bi-Lipschitz map, such that $\varphi(x)=x$, and $\varphi(A\cap \partial B(x,r))=K_1\cap\partial B(0,1)\subset K\cap \partial B(0,1)$. Note that $A=A+x$, hence the map $\varphi$ satisfies both $1^\circ$ and $2^\circ$.

Now let us prove $3^\circ$. We first define $\psi_0$ on $X=K\cap \partial B(0,1)$. Note that 
\be X=(\cup_{i=1}^m S_i)\cup (\cup_{j\in J}L_j)=\cup (\cup_{i=1}^m S_i)\cup (\cup_{j\in J'}L_j)\cup(\cup_{j\in J\bs J'}L_j).\ee

For each $z\in \cup_{i=1}^m S_i$, set $\psi_0(x)=0$; for $z\in \cup_{j\in J'}L_j=X_1$, set $\psi_0(z)=z$. For $j\in J\bs J'$, we have two cases: if $L_j\cap X_1=\emptyset$, then set $\psi_0(z)=0$; if $L_j$ intersects $X_1$ at one endpoint $y$, then let $\psi_0$ be a bi-Lipschitz map that maps $L_j$ to the segment $[0,y]$ such that $\psi_0(y)=y$; if $L_j$ intersects $X_1$ at both of its endpoints $y$ and $j$, then let $\psi_0$ be a bi-Lipschitz map from $L_j$ to $[y,0]\cup [0,z]$, such that $\psi_0(y)=y$ and $\psi_0(z)=z$. 

After the definition, $\psi_0$ is a map from $X$ to $K_1\cap \bar B(0,1)$, with $\psi_0(z)=z$ for $z\in X_1$.

Next we extend $\psi_0$ radially to a bi-Lipschitz map $\psi$ from $K\cap \bar B(0,1)$ to $K_1\cap \bar B(0,1)$: for each $z\in K\cap \bar B(0,1)$, set $\psi(z)=|z| \psi(\frac{z}{|z|}$. Then $\psi$ is a Lipschitz deformation retract from $K\cap \bar B(0,1)$ to $K_1=\varphi(\bar B(x,r)\cap (A+x))$, $\varphi=id$ on $K_1$, and $\psi\circ \varphi ((A+x)\cap \partial B(x,r))=(A+x)\cap \partial B(x,r)$.

This finishes the proof of \textbf{Case 2}. The last case is when $x$ is a regular point of $K$. The proof is similar to that of \textbf{Case 2}. We left the proof to the readers. \qed

\begin{cor}For any $t>0$ such that $(C_{k,p_1,p_2}+o_k)\cap \bar\mU'(x_1,t)=(C_k'+x_1)\cap \bar\mU'(x_1,t)$, there exists a Lipschitz map $\varphi:  (C_k'+x_1)\cap \bar\mU'(x_1,t)\to (C_{k,p_1,p_2}+o_k)\cap \bar\mU(o_k,r_k)$, such that 

$1^\circ$ $\varphi(x_1)=x_1$, and $\varphi((C_k'+x_1)\cap \partial\mU'(x_1,t))\subset (C_{k,p_1,p_2}+o_k)\cap \partial\mU(o_k,r_k)$; 

$2^\circ$ $\varphi:(C_k'+x_1)\cap \bar\mU'(x_1,t)\to \varphi ((C_k'+x_1)\cap \bar\mU'(x_1,t))$ is bi-Lipschitz;

$3^\circ$ There exists a Lipschitz deformation retract $\psi: (C_{k,p_1,p_2}+o_k)\cap \bar\mU(o_k,r_k)\to \varphi ((C_k'+x_1)\cap \bar\mU'(x_1,t))$, such that $\varphi^{-1}\psi\circ \varphi ((C_k'+x_1)\cap \partial\mU'(x_1,t))=(C_k'+x_1)\cap \partial\mU'(x_1,t)$.
\end{cor}

\nd $1^\circ$ For $i=1,2$ denote by $\t K^i=K^i_k+o_k+p_1$. Then $\t K^i$ is a translation of $K^i_k$ with center $o_k+p_i$. Recall that $\t K^1\cap \t K^2=x_1$, and in $\mU'(x_1,t_k)$, $\t K^1\cup \t K^2=C_{k,p_1,p_2}+o_k$ coincides with the minimal cone $C_k'+x_1$, where $C_k'$ is a union of two minimal cones $A_1$ and $A_2$, each $A_i\subset P_k^i, i=1,2$. As a result, each $(A_i+x_1)\cap \mU'(x_1,t)$ is part of $\t K^i$. Hence $A_i$ is of type $\P$, $\Y$ or $K^i$.

Then it is enough to prove the following: 
\be \begin{split}&\mbox{For each }i=1,2,\ \mbox{there exists a Lipschitz map }\varphi_i\mbox{ and a Lipschitz deformation retract }\psi_i, \\
&\mbox{ such that }1^\circ, 2^\circ\mbox{ and }3^\circ\mbox{ hold when we replace }\varphi\mbox{ by }\varphi_i,\ C_k'\mbox{ by }A_i,C_{k,p_1,p_2}+o_k\mbox{ by }\t K^i,\\
&\mbox{ and }\psi\mbox{ by }\psi_i.
\end{split}\ee

In fact, if \tb{(8.11)} holds, then we can just define $\varphi=\varphi_i(x)$ for $x\in (A_i+x_1)\cap \partial\mU'(x_1,t)$, and $\psi=\psi_i$ for $x\in \t K^i\cap \bar\mU(o_k,r_k)$.

Note that each convex domain $\mU(0,1)$ and $\mU'(0,1)$ is bi-Lipschitz equivalent to the unit ball $B(0,1)$, via a 1-homogeneous bi-Lipschitz map $l$, hence it is enough to prove \tb{(8.10)} while replacing $\mU$ and $\mU'$ by $B$. But this is follows immediately from Proposition \tb{8.2}, after translation and dilatation.\qed

\begin{pro}There exists $F'_k\subset F_k\cap \bar\mU'(x_1,t_k)$, such that
\be F'_k\in \oF_{\e t_k}((C'_k+x_1)\cap\bar \mU'(x_1,t_k), (C'_k+x_1)\cap\bar \mU'(x_1,t_k),\ee
and 
$F'_k$ is $10^{-2}\e t_k$ near $C'_k+x_1$ in $\mU'(x_1,t_k)$.
\end{pro}

\nd Let $p$ be a Lipschitz neighborhood deformation retract from $B((C_{k,p_1,p_2}+o_k)\cap\bar\mU(o_k,r_k), \e r_k)$ to $(C_{k,p_1,p_2}+o_k)\cap\bar\mU(o_k,r_k)$. Let $t$ be slightly larger than $t_k$, such that 
\be(C_k'+x_1)\cap \bar\mU'(x_1,t)=(C_{k,p_1,p_2}+o_k)\cap\bar\mU'(x_1,t).\ee

  Let $f:F_k\cap \bar\mU(o_k,r_k)\to [\bar\mU'(x_1,t_k)\cap F_k]\cup[(\bar\mU'(x_1,t)\bs \bar\mU'(x_1,t_k))\cap B((C_{k,p_1,p_2}+o_k),\e^8)]\cup [(C_{k,p_1,p_2}+o_k)\cap\bar\mU(o_k,r_k)\bs\bar\mU'(x_1,t)]$ be defined as the following :
\be f(x)=\left\{\begin{array}{rlc}x&,\ if\ &x\in \bar\mU'(x_1,t_k);\\
p(x)&,\ if\ &x\in \bar\mU(o_k,r_k)\bs \bar\mU'(x_1,t);\\
t_xp(x)+(1-t_x)x&,\ if\ &x\in \bar\mU'(x_1,t)\bs \bar\mU'(x_1,t_k),\end{array}\right.\ee
where $t_x$ is the value such that $x\in \partial \mU'(x_1, (1-t_x)t_k+t_xt)$.

Now let us construct the deformation from $(C_k'+x_1)\cap \bar\mU'(x_1,t_k)$ to $F_k\cap \bar\mU'(x_1,t_k)$:

First, take $\varphi$ as in Corollary \tb{8.3}. It maps $(C_k'+x_1)\cap \bar\mU'(x_1,t)$ to a subset of $(C_{k,p_1,p_2}+o_k)\cap\bar\mU(o_k,r_k)$, and for all $z\in (C_k'+x_1)\cap \partial\mU'(x_1,t_k)$, $\varphi(z)\in (C_{k,p_1,p_2}+o_k)\cap\partial\mU(o_k,r_k)$.

Let $\theta:\bar\mU'(x_1,t_k)\to \bar\mU'(x_1,t)$  be the dilatation centered at $x_1$, that is:
$\theta(x)=x_1+\frac{t}{t_k}(x-x_1)$

Since $F_k\cap \bar\mU(o_k,r_k)\in \oF_{\e^8}((C_{k,p_1,p_2}+o_k)\cap \bar\mU(o_k,r_k),\bar\mU(o_k,r_k))$, there exists a deformation $\varphi_1:\bar\mU(o_k,r_k)\to \bar\mU(o_k,r_k)$, such that $\varphi_1((C_{k,p_1,p_2}+o_k)\cap \bar\mU(o_k,r_k)=F_k\cap \bar\mU(o_k,r_k)$, $\varphi_1(\partial\mU(o_k,r_k))\subset \mU(o_k,r_k)$, and for $z\in \partial\mU(o_k,r_k)$, $|\varphi_1(z)-z|<\e^8 r_k$. In particular, $\varphi_1(F_k\cap \partial\mU(o_k,r_k))\subset B((C_{k,p_1,p_2}+o_k)\cap \partial\mU(o_k,r_k), \e^8)$.

Then $\varphi_1\circ\varphi\circ\theta[(C_k'+x_1)\cap \bar\mU'(x_1,t_k)]$ is a subset of $F_k$, such that 
\be \varphi_1\circ\varphi\circ\theta[(C_k'+x_1)\cap \partial\mU'(x_1,t_k)]\subset \partial \mU(o_k,r_k)\cap B((C_{k,p_1,p_2}+o_k)\cap \partial\mU(o_k,r_k), \e^8).\ee

Next we apply $f$ to $\varphi_1\circ\varphi\circ\theta[(C_k'+x_1)\cap \bar\mU'(x_1,t_k)]$. Then the image 
\be \begin{split}f\circ\varphi_1&\circ\varphi\circ\theta[(C_k'+x_1)\cap \bar\mU'(x_1,t_k)]\subset \\
&[\bar\mU'(x_1,t_k)\cap F_k]
\bigcup[(\bar\mU'(x_1,t)\bs \bar\mU'(x_1,t_k))\cap B((C_{k,p_1,p_2}+o_k),\e^8)]\\
&\bigcup [(C_{k,p_1,p_2}+o_k)\cap\bar\mU(o_k,r_k)\bs\bar\mU'(x_1,t)],\end{split}\ee
and $f\circ\varphi_1\circ\varphi\circ\theta[(C_k'+x_1)\cap \partial\mU'(x_1,t_k)]\subset \partial \mU(o_k,r_k)\cap (C_{k,p_1,p_2}+o_k)$. Moreover, for $z\in \varphi[(C_k'+x_1)\cap \bar\mU'(x_1,t_k)]$, we know that $|f\circ\varphi_1(z)-z|<2\e^8$.

Now we define $g: [\bar\mU'(x_1,t_k)\cap F_k]
\cup[(\bar\mU'(x_1,t)\bs \bar\mU'(x_1,t_k))\cap B((C_{k,p_1,p_2}+o_k),\e^8)]
\cup [(C_{k,p_1,p_2}+o_k)\cap\bar\mU(o_k,r_k)\bs\bar\mU'(x_1,t)]\to  [\bar\mU'(x_1,t_k)\cap F_k]
\cup[(\bar\mU'(x_1,t)\bs \bar\mU'(x_1,t_k))\cap B((C_{k,p_1,p_2}+o_k),\e^8)]$, such that $g|_{[(C_{k,p_1,p_2}+o_k)\cap\bar\mU(o_k,r_k)\bs\bar\mU'(x_1,t)]}$ is just the deformation retract $\varphi^{-1}\circ\psi $ from $(C_{k,p_1,p_2}+o_k)\cap\bar\mU(o_k,r_k)$ to $C_k'\cap \bar \mU'(x_1,t)$, as in Corollary \tb{8.3}, and $g|_{\bar \mU'(x_1,t)}=id$.

Then the set 
\be\begin{split}g \circ f&\circ\varphi_1\circ\varphi\circ\theta[(C_k'+x_1)\cap \bar\mU'(x_1,t_k)]\\
&\subset[\bar\mU'(x_1,t_k)\cap F_k] \cup[(\bar\mU'(x_1,t)\bs \bar\mU'(x_1,t_k))\cap B((C_{k,p_1,p_2}+o_k),\e^8)]\\
&=[\bar\mU'(x_1,t_k)\cap F_k] \cup[(\bar\mU'(x_1,t)\bs \bar\mU'(x_1,t_k))\cap B((C'_k+x_1),\e^8)]\\
&\subset \bar\mU'(x_1,t),\end{split}\ee 
and 
\be g \circ f\circ\varphi_1\circ\varphi\circ\theta[(C_k'+x_1)\cap \partial\mU'(x_1,t_k)]= (C_k'+x_1)\cap\partial\mU'(x_1,t).\ee

Let \be\pi: \bar\mU'(x_1,t)\to \bar\mU'(x_1,t_k): 
\pi (z)=\left\{\begin{array}{rcl}\frac{t_k}{|z-x_1|}(z-x_1)+x_1&,\ if\ &z\in \bar\mU'(x_1,t)\bs \bar\mU'(x_1,t_k)\\
z&,\ if\ &z\in \bar\mU'(x_1,t_k). \end{array}\right.\ee

Then $\pi$ is a "radial projection" centered at $x_1$ from $\bar\mU'(x_1,t)$ to $\bar\mU'(x_1,t_k)$. By \tb{(8.16)}, we know that 
\be\begin{split}\pi\circ g&\circ f\circ\varphi_1\circ\varphi\circ\theta[(C_k'+x_1)\cap \bar\mU'(x_1,t_k)]\\
&\subset [\bar\mU'(x_1,t_k)\cap F_k]\cup [\partial \mU'(x_1,t_k)\cap B((C_k'+x_1), \e^8 r_k))]\\
&\subset [\bar\mU'(x_1,t_k)\cap F_k]\cup [\partial \mU'(x_1,t_k)\cap B(F_k, 2\e^8 r_k))],
\end{split}\ee
and for any $z\in \bar\mU'(x_1,t_k)\cap (C_k'+x_1)$, we have
\be |\pi\circ g\circ f\circ\varphi_1\circ\varphi\circ\theta(z)-z|<3\e^8 r_k\le 3 \e^2 t_k.\ee

Last, let $h_0$ be a 2-Lipschitz neighborhood deformation retract from $B(F_k, 4\e^2 t_k)\cap \partial\mU'(x_1,t_k)\to F_k\cap \partial\mU'(x_1,t_k)$. The existence of such an $h_0$ is by a same argument as in the proof of Proposition \tb{6.2}. Define $h:\pi\circ g\circ f\circ\varphi_1\circ\varphi\circ\theta[(C_k'+x_1)\cap \bar\mU'(x_1,t_k)]\to F_k\cap \bar\mU'(x_1,t_k)$ by $h=h_0$ on $B(F_k, 4\e^2 t_k)\cap \partial\mU'(x_1,t_k)$, and $h=id$ on $F_k\cap \mU'(x_1,t_k)$. Then $h\circ\pi\circ g\circ f\circ\varphi_1\circ\varphi\circ\theta[(C_k'+x_1)\cap \bar\mU'(x_1,t_k)]\subset F_k\cap\bar\mU'(x_1,t_k)$. Moreover, for $z\in (C_k'+x_1)\cap \pa\mU'(x_1,t_k)$, we know that $h\circ\pi\circ g\circ f\circ\varphi_1\circ\varphi\circ\theta(z)\in \partial \mU'(x_1,t_k)$, and $|h\circ\pi\circ g\circ f\circ\varphi_1\circ\varphi\circ\theta(z)-z|\le 8\e^2 t_k$. 

Now we extend $h\circ\pi\circ g\circ f\circ\varphi_1\circ\varphi\circ\theta$ to a Lipschitz deformation $\xi$ from $\bar \mU'(x_1,t_k)$ to $\bar \mU'(x_1,t_k)$, such that $\xi(\partial\mU'(x_1,t_k))\subset \partial\mU'(x_1,t_k)$ . Then $|\xi(z)-z|\le 8\e^2 t_k$ for $z\in (C_k'+x_1)\cap \partial\mU'(x_1,t_k)$. Therefore 
\be F'_k:=\xi((C_k'+x_1)\cap \bar\mU'(x_1,t_k))\in \oF_{\e t_k}((C_k'+x_1)\cap \bar\mU'(x_1,t_k),\bar\mU'(x_1,t_k)),\ee and 
\be F'_k\subset F_k\cap \bar\mU'(x_1,t_k).\ee

The rest is to prove that 
\be d_{\mU'(x_1,t_k),t_k}(F'_k,(C_k'+x_1)\cap \bar\mU'(x_1,t_k))<10^{-2}\e. \ee

We already know that $F_k$ is $\e^8r_k\le \e^2 t_k$ near $C_k'+x_1$ in $\mU(x_1,t_k)$, hence 
\be F'_k\subset F_k\cap \bar\mU'(x_1,t_k)\subset B(C_k'+x_1, \e^2 t_k).\ee
So we only have to prove that 
\be (C_k'+x_1)\cap \mU'(x_1,t_k)\subset B(F'_k, 10^{-2}\e t_k).\ee

 Suppose not. Then there exists $z_0\in (C_k'+x_1)\cap \mU'(x_1,t_k)$, such that $B(z_0,10^{-2}\e t_k)\cap F'_k=\emptyset$. Set $\Omega=\bar\mU'(x_1,t_k)\bs \bar B(z_0,\frac 12\times 10^{-2}\e t_k)$. Since $10^{-2}\e$ is much larger than $\e^2$, by regularity of minimal cones, there exists a Lipschitz neighborhood deformation retract $p'$ from $\Omega\cap B((C_k'+x_1), \e^2)$ to $\Omega\cap (C_k'+x_1)$. Then $p'(F'_k)\subset (C_k'+x_1)\cap \mU'(x_1,t_k)\bs \bar B(z_0,\frac 12\times 10^{-2}\e t_k)$. Note that $F'_k\in 
\oF_{\e t_k}((C_k'+x_1)\cap \bar\mU'(x_1,t_k),\bar\mU'(x_1,t_k))$, hence $p'(F'_k)\in\oF(C_k'+x_1,\mU'(x_1,t_k))$. Since the density of $C'_k+x_1$ at $z_0$ is at least 1, we have

\be \begin{split}\H^2(p'(F'_k))&\le \H^2((C_k'+x_1)\cap \bar\mU'(x_1,t_k)\bs \bar B(z_0,\frac 12\times 10^{-2}\e t_k))\\
&= \H^2((C_k'+x_1)\cap \bar\mU'(x_1,t_k))-H^2((C_k'+x_1)\cap\bar B(z_0,\frac 12\times 10^{-2}\e t_k))\\
&\le \H^2((C_k'+x_1)\cap \bar\mU'(x_1,t_k))-\frac 14 10^{-4}\e^2 t_k^2.\end{split}\ee

Now let us consider the restriction $\xi_1$ of $\xi$ to $(C_k'+x_1)\cap \partial\mU'(x_1,t_k)$. Since $\xi$ is a deformation in $\bar \mU'(x_1,t_k)$, we know that $p'\circ \xi_1$ is a deformation from $(C_k'+x_1)\cap \partial\mU'(x_1,t_k)$ to $(C_k'+x_1)\cap \partial\mU'(x_1,t_k)$. So let $H:(C_k'+x_1)\cap \partial\mU'(x_1,t_k)\times [0,1]$ be continuous, such that $H(\cdot, 0)=id$, and $H(\cdot, 1)=p'\circ \xi_1$.

Now we define $G: (C_k'+x_1)\cap\bar\mU'(x_1, 2t_k)\bs \bar\mU'(x_1, t_k)\to (C_k'+x_1)\times [0,1]$: for any $z\in (C_k'+x_1)\cap\bar\mU'(x_1, 2t_k)\bs \bar\mU'(x_1, t_k)$, set $t_z$ be such that $z\in \partial \mU'(x_1, t_z)$, and set $G(z)=(z, \frac{2t_k-t_z}{t_k})$. Then we set $\xi_2:(C_k'+x_1)\cap \bar\mU'(x_1, 2t_k)\to (C_k'+x_1)\cap\bar\mU'(x_1, 2t_k)$: 
\be\xi_2(z)=\left\{\begin{array}{rcl}\xi(z)&,\ if\ &z\in (C_k'+x_1)\cap \bar\mU'(x_1, t_k);\\
\frac{t_z}{t_k}(H\circ G(z)-x_1)+x_1\end{array}\right.\ee

We extend $\xi_2$ to a Lipschitz deformation from $\mU'(x_1, 2t_k)$ to itself, with $\xi_2|_{\partial \mU'(x_1, 2t_k)}=id$. Then $F_k'':=\xi_2((C_k'+x_1)\cap\bar\mU'(x_1, 2t_k))\in \oF(C_k'+x_1, \mU'(x_1,2t_k)$, and
 
\be \begin{split}F_k''&=\xi_2((C_k'+x_1)\cap\bar\mU'(x_1, 2t_k))\\
&=\xi_2((C_k'+x_1)\cap\bar\mU'(x_1, t_k))\cup \xi_2((C_k'+x_1)\cap\bar\mU'(x_1, 2t_k)\bs \bar\mU'(x_1, t_k))\\
&\subset p'(F'_k)\cup[(C_k'+x_1)\cap\bar\mU'(x_1, 2t_k)\bs \bar\mU'(x_1, t_k)].\end{split}\ee

The last union is disjoint, hence by \tb{(8.27)} we have
\be \begin{split}\H^2(F_k'')&=\H^2(\xi_2((C_k'+x_1)\cap\bar\mU'(x_1, 2t_k))\\
&\le \H^2(p'(F'_k))+\H^2((C_k'+x_1)\cap\bar\mU'(x_1, 2t_k)\bs \bar\mU'(x_1, t_k))\\
&\le\H^2((C_k'+x_1)\cap \bar\mU'(x_1,t_k))-\frac 14 10^{-4}\e^2 t_k^2+\H^2((C_k'+x_1)\cap\bar\mU'(x_1, 2t_k)\bs \bar\mU'(x_1, t_k))\\
&=\H^2((C_k'+x_1)\cap \bar\mU'(x_1,2t_k))-\frac 14 10^{-4}\e^2 t_k^2<\H^2((C_k'+x_1)\cap \bar\mU'(x_1,2t_k)).
\end{split}\ee

But $F_k'':=\xi_2((C_k'+x_1)\cap\bar\mU'(x_1, 2t_k))\in \oF(C_k'+x_1, \mU'(x_1,2t_k)$, and $C_k'+x_1$ is a minimal cone, hence \tb{(8.30)} contradicts Corollary \tb{4.7} $2^\circ$ of \cite{uniquePYT}.

Hence we get \tb{(8.24)}. Thus we complete the proof of Proposition \tb{8.4}.\qed

\subsection{The end of the reduction}

Now let us summerize what happened:

Recall that we have the following:

$1^\circ$ $\H^2(F_k\bs \mU(o_k,r_k))\ge \H^2((C_k+o_k)\cap B(0,1)\bs \mU(o_k,r_k))$;

$2^\circ$ $F_k\cap \bar\mU(o_k,r_k)\in \oF_{4\e r_k}((C_k+o_k)\cap\bar\mU(o_k,r_k), \bar\mU(o_k,r_k))$;

$3^\circ$ In the critical region $\mU(o_k,r_k)$, $F_k$ is $2\e r_k$ near $C_k+o_k$, but is not $\e r_k$ near any translation of $C_k$. In the annulus $\mU(o_k,r_k)\bs \mU(o_k,\frac 14 r_k)$, $F_k$ is decomposed into two disjoint part $F_k^1$ and $F_k^2$, and each $F_k^i$ is $2\e r_k$ near $K^i_k$ in $\mU(o_k,r_k)\bs \mU(o_k,\frac 14 r_k)$. Then we have two situations :

\textbf{Case 1:} At least one of the $F_k^i$ is $\d(\e^8) r_k$ far from $K^i_k$ in $\mU(o_k,r_k)\bs \mU(o_k,\frac 14 r_k)$; 

\textbf{Case 2:} Each $F_k^i$ is $\d(\e^8) r_k$ near some translation $K^i_k+o_k+p_i$ of $K^i_k$ in $\mU(o_k,r_k)\bs \mU(o_k,\frac 14 r_k)$, $i=1,2$.

If we goes to \textbf{Case 1}, the reduction ended, and we can do all the estimates in Section 9;

If we goes to \textbf{Case 2}, then we can find a point $x_1\in B(o_k,\frac 19 r_k)$, a radius $t_k\in (\e^6 r_k,\e r_k)$, and a minimal cone $C_k'$ which is a union of two Almgren unique minimal cones $A_1\subset P_k^1$ and $A_2\subset P_k^2$, such that 

$1^\circ$ $A_i$ is simpler than $K^i$, and at least one of the $A_i$ is strictly simpler than $K^i$.

$2^\circ$ $F_k\cap \bar \mU'(x_1,t_k)$ contains a part $F'_k$ that belongs to the class $\oF_{\e t_k}((C'_k+x_1)\cap \bar\mU'(x_1,t_k),\bar\mU'(x_1,t_k))$;

$3^\circ$ $F'_k$ is $10^{-2}\e t_k$ near $C'_k+x_1$ in $\mU'(x_1,t_k)$; 

$4^\circ$ $\H^2(F_k\cap\mU(o_k,r_k)\bs \mU'(x_1,t_k))\ge \H^2[((K^1_k+o_k+p_1)\cup (K^2_k+o_k+p_2))\cap \mU(o_k,r_k)\bs \mU'(x_1,t_k)]$.

Here $\mU'$ denote the convex domain corresponding to $C_k'$.
Due to $1^\circ$-$4^\circ$, we can restart the $\e$-process with respect to $F'_k$, $C_k'$ in $\mU'(x_1,t_k)$, and arrive to a new critical region $\mU(o_k', r_k')$, where we fall again onto the two cases above. We replace $o_k$ and $r_k$ by $o_k'$ and $r_k'$. If \textbf{Case 1} happens, the reduction ended , and get the same estimates after Section \tb{7}; If \textbf{Case 2} happens, we continue the reduction in this section, and get a new minimal cone $C_k^{(1)}$, which is strictly simpler than $C'_k$, and new critical region $\mU^{(1)}(o_k',r_k')$.......

This reduction procedure will stop at a finite time, because there are only finitely many choices of pairs of minimal cones that are strictly simpler than the pair $(K^1,K^2)$. We stop by falling into \textbf{Case 1} in some intermediate step, or we get finally to the cone $Q^1\cup Q^2$, where $Q^1$ and $Q^2$ are planes. For the union of two planes, when we do the $\e$-process, we will finally goes to \textbf{Case 1}, because \textbf{Case 2} cannot happen since planes are translation invariant.

In any case, when the reduction procedure stops, we get some $o_k\in B(0,1)$, and $r_k>0$, such that we have the following estimates:

\begin{pro}There exists $0<\e_0<\min\{10^{-5}, 10^{-3}\}$, such that for each $\e<\e_0$, for $k$ large, there exists $o_k\in B(0,10^{-2})$, $r_k<10^{-2}$, a minimal cone $C_k'$ which is a union of two sliding stable Almgren unique minimal cones $A^1_k\subset P_k^1$, and $A^2_k\subset P^2_k$, each $A^i_k$ being simpler than $K^i_k$, such that the following holds:

$1^\circ$ $\H^2(F_k\bs \mU'(o_k,r_k))+\H^2((C_k'+o_k) \cap\mU'(o_k,r_k))\ge \H^2(C_k)$, where $\mU'=\mU_k'$ is the convex domain corresponds to $C_k'$;

$2^\circ$ $\H^2(F_k\cap \mU'(o_k,r_k)\bs \mU'(o_k,\frac 14 r_k))\ge \H^2((C_k'+o_k)\cap \mU'(o_k,r_k)\bs \mU'(o_k,\frac 14 r_k))+C_4(\e)r_k^2$, where $C_4(\e)$ depends only on $\e$ and $n$; 

$3^\circ$ $\H^2(F_k\bs \mU'(o_k,\frac 14r_k))+\H^2((C_k'+o_k)\cap \mU'(o_k,\frac 14r_k))\ge \H^2(C_k)+C_4(\e)r_k^2$;

$4^\circ$ $F_k\cap \bar\mU'(o_k,\frac 14 r_k)\subset \oF_{4\e r_k}((C_k'+o_k)\cap \bar\mU'(o_k,\frac 14r_k),\bar\mU'(o_k,\frac 14r_k)$.
\end{pro}

\nd $1^\circ$ We prove by recurrence. When we finish the initial $\e$-process using $C_k$, then by Corollary \tb{7.5}, we have
\be \begin{split}&\H^2(F_k\bs \mU(o_k,r_k))+\H^2((C_k+o_k)\cap \mU(o_k,r_k))\\
&\ge \H^2(C_k\bs \mU(o,r_k))+\H^2((C_k+o_k)\cap \mU(o_k,r_k))\\
&=\H^2(C_k\bs \mU(o,r_k))+\H^2(C_k\cap \mU(o,r_k))\\
&=\H^2(C_k).
\end{split}\ee
If at this step, we come to Possibility (2) at the end of Subsection \tb{7.2}, then we get $1^\circ$ directly. Otherwise, we start the reduction argument. In $\mU(o_k,r_k)$, $F_k$ is $\e^8r_k$ near some $C_{k,p_1,p_2}+o_k$ and is also $2\e r_k$ near $C_k+o_k$, hence the two points $p_1,p_2$ are in the ball $B(o_k,3\e r_k)$. By stability of measure \cite{stablePYT} of each $K^i_k$, we know that
\be \H^2((K^i_k+o_k)\cap \mU(o_k,r_k))=\H^2((K^i_k+o_k+p_i)\cap \mU(o_k,r_k)),\ee
which gives
\be \H^2(C_k\cap \mU(o_k,r_k))=\H^2((C_{k,p_1,p_2}+o_k)\cap \mU(o_k,r_k)).\ee
On the other hand, by \tb{(8.5) and (8.6)}, we know that
\be \H^2(F_k\cap \mU(o_k,r_k)\bs \mU'(x_1,t_k))\ge \H^2((C_{k,p_1,p_2}+o_k)\cap \mU(o_k,r_k)\bs \mU'(x_1,t_k)).\ee
Then we restart the $\e$-process for $F_k\cap \mU'(x_1,t_k)$, with the cone $(C_k'+x_1)\cap \mU'(x_1,t_k)$ which coincides with $C_{k,p_1,p_2}+o_k$ in $\mU'(x_1,t_k)$ and find the new critical region $\mU'(o_k',r_k')$. Then we get similarly as in Corollary \tb{7.5}, that
\be \begin{split}\H^2&(F_k\cap \mU'(x_1,t_k)\bs \mU'(o_k',r_k'))
\ge \H^2((C_k'+x_1)\cap \mU'(x_1,t_k)\bs \mU'(x_1,r_k'))\\
&=\H^2((C_k'+x_1)\cap \mU'(x_1,t_k))-\H^2((C_k'+x_1)\cap\mU'(x_1,r_k'))\\
&=\H^2((C_{k,p_1,p_2}+o_k)\cap \mU'(x_1,t_k))-\H^2((C_k'+o_k')\cap\mU'(o_k',r_k')).\end{split}\ee
Combining \tb{(8.31), (8.33), (8.34) and (8.35)}, we get
\be \H^2(F_k\bs \mU'(o_k',r_k'))+\H^2((C_k'+o_k')\cap\mU'(o_k',r_k'))\ge \H^2(C_k).\ee

By recurrence, we get $1^\circ$;

$2^\circ$ This follow directly from the reduction procedure, and \tb{(7.81)};

$3^\circ$ Follows directly from $1^\circ$ and $2^\circ$;

$4^\circ$ Follows from Proposition \tb{8.4}, and Proposition \tb{6.3}.\qed

\section{Conclusion}

We will finish our proof in this section. A final step is to control the measure of $F_k$ in $\mU'(o_k,r_k)$. 

\begin{pro}Take all the notations as before. Then for each $k$, we have
\be \H^2(F_k\cap \mU'(o_k,\frac 14 r_k))\ge (1-C(\theta_k))\H^2((C_k'+o_k)\cap \mU'(o_k,\frac 14r_k)),\ee
with $C(\theta)\to 0$ as $\theta\to \frac \pi 2$.
\end{pro}

\nd By $4^\circ$ of Proposition \tb{8.5}, we know that for each $k$ large,
\be  F_k\cap \bar\mU'(o_k,\frac 14 r_k)\subset \oF_{4\e r_k}((C_k'+o_k)\cap \bar\mU'(o_k,\frac 14r_k),\bar\mU'(o_k,\frac 14r_k)).\ee

Fix any $k$ large. Let $C'_k=A_1\cup A_2$, where $A_i\subset P^i_k,i=1,2$ are measure and sliding stable Almgren unique minimal cones. Let $U_i\subset P^k_i,i=1,2$ be the corresponding convex domain of $A_i$. Then $\mU'_k=(U_1\times {P_k^1}^\perp)\cap ({P_k^2}^\perp\times U_2)$.

Now we do projections to $P^k_1$ and $P^k_2$, and get that for $i=1,2$, 
\be p_k^i(F_k\cap \bar\mU'(o_k,\frac 14 r_k))\in \oF_{4\e r_k}((A_i+o_k)\cap \bar U_i(o_k,\frac 14r_k),\bar U_i(o_k,\frac 14r_k),\ee
where $p^k_i$ denote the orthogonal projection from $\R^n$ to $P_k^i$.

By stability of $A_i$, we know that
\be \H^2(p_k^i(F_k\cap \bar\mU'(o_k,\frac 14 r_k)))\ge \H^2((A_i+o_k)\cap \bar U_i(o_k,\frac 14r_k)).\ee

Now by Lemmas \tb{2.34 and 2.35}, we know that 
\be \begin{split}\H^2(F_k\cap \bar\mU'(o_k,\frac 14 r_k))&\ge \frac{1}{1-d\cos\theta_k}\sum_{i=1}^2 \H^2(p_k^i(F_k\cap \bar\mU'(o_k,\frac 14 r_k)))\\
&\ge \frac{1}{1-d\cos\theta_k}\sum_{i=1}^2\H^2((A_i+o_k)\cap \bar U_i(o_k,\frac 14r_k))\\
&=\frac{1}{1-d\cos\theta_k}\H^2((C_k'+o_k)\cap \bar\mU'(o_k,\frac 14r_k)).
\end{split}\ee\qed

Set $C(\theta)=1-\frac{1}{1+d\cos\theta}=\frac{d\cos\theta}{1+d\cos\theta}$. Then \tb{(9.1)} holds, and $C(\theta)\to 0$ as $\theta\to \frac \pi 2$. \qed

To finish the proof of Theorem \ref{main}, we combine Proposition \tb{8.5} $3^\circ$ and Proposition \tb{9.1}, and get
\be \begin{split}\H^2(F_k)=&\H^2(F_k\cap \bar\mU'(o_k,\frac 14 r_k))+\H^2(F_k\bs\bar\mU'(o_k,\frac 14 r_k))\\
=&\H^2(F_k\cap \bar\mU'(o_k,\frac 14 r_k))+\H^2(F_k\bs\mU'(o_k,\frac 14 r_k))+\H^2(F_k\cap \partial \mU'(o_k,\frac 14 r_k))\\
\ge &(1-C(\theta_k))\H^2((C_k'+o_k)\cap \bar\mU'(o_k,\frac 14r_k))+\H^2(C_k)-\H^2((C_k'+o_k)\cap \mU'(o_k,\frac 14 r_k))\\
&+C_5(\e)r_k^2+\H^2(F_k\cap \partial \mU'(o_k,\frac 14 r_k))\\
=& \H^2(C_k)-C(\theta_k)\H^2((C_k'+o_k)\cap \bar\mU'(o_k,\frac 14r_k))+C_5(\e)r_k^2+\H^2(F_k\cap \partial \mU'(o_k,\frac 14 r_k)).
\end{split}
\ee
Since $F_k$ is $2\e r_k$ near $C_k'+o_k$ in $\mU'(o_k,\frac 14 r_k)$, by Corollary \tb{2.24}, we know that 
\be \H^2(F_k\cap \partial \mU'(o_k,\frac 14 r_k))=0;\ee
On the other hand, let $d$ denote the density of $C_k$ at $0$, then since $C_k'$ is simpler than $C_k$, we know that the density of $C_k'$ at zero is no more than $d$, therefore
\be \H^2((C_k'+o_k)\cap \bar\mU'(o_k,\frac 14r_k))\le d(\frac 14 r_k)^2=2^{-4}d r_k^2.\ee
As a result, 
\be \H^2(F_k)\ge \H^2(C_k)-C(\theta_k)2^{-4}d r_k^2+C_5(\e)r_k^2=\H^2(C_k)+(C_5(\e)-C(\theta_k)2^{-4}d)r_k^2.\ee

But we know that $\lim_{k\to\infty}C(\theta_k)=0$, therefore for $k$ is large enough, we have
\be C_5(\e)-C(\theta_k)2^{-4}d>\frac 12C_5(\e),\ee
which gives
\be \H^2(F_k)\ge \H^2(C_k)+\frac 12C_5(\e)r_k^2>\H^2(C_k).\ee
This contradicts Proposition \tb{4.1 $1^\circ$}. Thus, the proof of Theorem \ref{main} is completed.\qed

As a result, combining with the discuss in \cite{stablePYT, uniquePYT, stableYXY}, we have

\begin{cor}For all the known 2-dimensional Almgren minimal cones, their almost orthogonal union is minimal.
\end{cor}

\section{Perspectives}

\subsection{The Almgren uniqueness, sliding and measure stability for the union}

Theorem \ref{main} gives us families of new 2-dimensional minimal cones. As stated before, all other known 2-dimensional minimal cones are Almgren unique and sliding stable, so the almost orthogonal union of any two of them is still minimal. But we wonder whether this union is Almgren unique and sliding stable, so that we can continue to do almost orthogonal unions and get new types of 2-dimensional minimal cones.

The answer is partially affirmative in the following sense:

\begin{thm}Let $K_0^1\subset \R^{n_1}$ and $K_0^2\subset \R^{n_2}$ be two 2-dimensional Almgren unique minimal cones. We suppose in addition that both cones are sliding stable. Then there exists $0<\theta_1<\frac\pi2$, such that if $K^1$ and $K^2$ are two copies (centered at 0) of $K_0^1$ and $K_0^2$ in $\R^{n_1+n_2}$, with $\a_{K^1,K^2}\ge\theta_1$, 
then the union $K^1\cup K^2$ is a sliding stable Almgren unique minimal cone.
\end{thm}

Note that here the angle $\theta_1$ might be different from the angle $\theta_0$ in Theorem \ref{main}. We do not guarantee that $K^1\cup K^2$ is Almgren unique or sliding stable once the union is minimal.

\noindent{\bf Proof of Theorem \tb{10.1}.} The proof is almost the same as that of Theorem \ref{main}. Let us sketch it: 

First suppose the conclusion of Almgren uniqueness in Theorem \tb{10.1} does not hold for two 2-dimensional sliding stable Almgren unique minimal cones $K_0^1$ and $K_0^2$. Then there exists a sequence of angles $\theta_k\to \frac \pi2$, such that the cone $C_k=K_k^1\cap_{\theta_k} K_k^2$ is minimal but not Almgren unique, $k\to\infty$. Then for each $k$, there exists a set $F_k\in \oF(C_k,\mU_k(0,1))$, $F_k\ne C_k$, with $\H^2(F_k)=\H^2(C_k)$. By the upper semi continuity Theorem \tb{4.1} of \cite{uniquePYT}, we know that $F_k$ is also Almgren minimal in $\mU_k(0,1)$. Then the limit of $F_k$ is $C_0$. For the sequence $F_k$, we employ the same argument from Section \tb{5}: first we start a stopping time argument. But this time, to prove that the $\e-$process has to stop at a finite step, we use Proposition \tb{5.3}: if the $\e$-process never stops, then by the argument in Proposition \tb{5.4}, we can decompose $F_k$ into the essentially disjoint union of two parts $F^1$ and $F^2$, $F^i$ contains an element $F^i_0\in \oF(K^i_k, B(0,1))$, $i=1,2$. Therefore we have
\be \H^2(F_k)\ge \H^2(F^1_0)+\H^2(F^2_0)\ge \H^2(K^1_k)+\H^2(K^2_k)=\H^2(C_k)=\H^2(F_k).\ee

As consequence, $\H^2(F^i)=\H^2(F^i_0)=\H^2(K^i_k)$, $i=1,2$. By Almgren uniqueness of $K^i_k$, we know that $F^i=F^i_0=K^i_k$ (modulo $\H^2$-null sets), and hence $F_k=K^1_k\cup K^2_k=C_k$, which contradicts the fact that $F_k\ne C_k$.

Hence the $\e$-process has to stop at a finite step. 

Then we continue all the same arguments until the end of Section \tb{9}, and get that for $k$ large enough, $\H^2(F_k)>\H^2(C_k)$, which leads to a contradiction.

Next for the sliding stability: suppose there exists a sequence of angles $\theta_k\to \frac \pi2$, such that the cone $C_k=K_k^1\cap_{\theta_k} K_k^2$ is sliding stable, $k\to\infty$. Then for each $k$, there exists a set $F_k\in \oF_{\d_k}(C_k,\mU_k(0,1))$, $F_k\ne C_k$, with $\H^2(F_k)<\H^2(C_k)$, where $\d_k\to 0$, and $\d_k$ is less than the $\d$ for which $K^i$ is $(\eta,\d)$-sliding stable, $i=1,2$.

Then the limit of $F_k$ is again $C_0$, and we employ exactly same argument as from Section \tb{5 to 9}, and get that $\H^2(F_k)>\H^2(C_k)$ for $k$ large, contradiction.\qed

Using recurrence and Theorem \tb{10.1}, we have

\begin{cor}The almost orthogonal unions of several sliding stable Almgren unique 2-dimensional minimal cones are sliding stable Almgren unique 2-dimensional minimal cones. More precisely, for any group of $K_1, \cdots K_n$ of 2-dimensional sliding stable Almgren unique 2-dimensional minimal cones, there exists $\theta>0$, such that if $\a_{K_i,K_j}\ge \theta$ for $i\ne j$, then their union
\be \cup_{i=1}^n K_j\ee
is an sliding stable Almgren unique 2-dimensional minimal cone.
\end{cor}

\subsection{Topological minimality}

It is natural to ask whether the almost orthogonal union of two 2-dimensional $G$-topological minimal cones is still $G$-topologically minimal. The answer is yes with the same conditions ($G$-topological unique, and some condition similar to sliding minimality). The proof is similar, and simpler at some places. 

The main differences are places where we proved that some set is either a deformation or a limit of deformations, there we have to prove that those sets are $G$-topological competitors instead. But fortunately, since deformations are $G$-topological competitors (Proposition \tb{2.7}), and limits of $G$-topological competitors are still $G$-topological competitors, we do not need any additional work, and at some places we can even save some time. An example is that we do not need the upper semi continuity Theorem (Theorem 4.1 of \cite{uniquePYT}). 

Otherwise, all regularity theorems for Almgren minimal sets hold for $G$-topological minimal ones, by Proposition \tb{2.7}.

So let us state the result for almost orthogonal unions of $G$-topological minimal cones:

\begin{thm}[minimality of the union of two almost orthogonal 2-dimensional $G$-topological minimal cones]Let $K_0^1\subset \R^{n_1}$ and $K_0^2\subset \R^{n_2}$ be two 2-dimensional $G$-topologically unique minimal cones. We suppose in addition that both cones are $G$-topological sliding stable. Then there exists $0<\theta_0<\frac\pi2$, such that if $K^1$ and $K^2$ are two copies (centered at 0) of $K_0^1$ and $K_0^2$ in $\R^{n_1+n_2}$, with $\a_{K^1,K^2}\ge\theta_0$, 
then the union $K^1\cup K^2$ is a $G$-topological unique $G$-topological sliding stable minimal cone.
\end{thm}

See the beginning of Section \tb{3} for the definition of $G$-topological uniqueness. The $G$-topological sliding stability means the following:

\begin{defn}Let $K$ be a 2-dimensional $G$-topological minimal cone in $\R^n$. Let $0<\d<\eta<10^{-2}\eta_0$. 
$1^\circ$ We say that a closed set $F$ is an $(\eta,\d)$-sliding $G$-topological competitor for $K$, if there exists a 2-dimensional $G$-topological competitor $E$ of $K$ in $\cU(K,\eta)$, such that $F$ is a $\d$-sliding deformation of $E$ in $\bar\cU(K,\eta)$.

$2^\circ$ We say that $K$ is $(\eta,\d)$-$G$-topological sliding stable, if for all $(\eta,\d)$-sliding $G$-topological competitor $F$ of $K$, we have
\be \H^d(F\cap \bar\cU(K,\eta))\ge \H^d(K\cap \bar\cU(K,\eta)).\ee
\end{defn}

%
%

\begin{rem}$1^\circ$ Similarly results for $G$-topological minimal cones as in Theorem \tb{10.1} and Corollary \tb{10.2} also hold.

$2^\circ$ The $G$-topological uniqueness and sliding stability also hold for all known $G$-topological minimal cones. See \cite{stablePYT} and \cite{stableYXY}.
\end{rem}
%
%
%
%
%
%
%
%
%

\renewcommand\refname{References}
\bibliographystyle{plain}
\bibliography{reference}

\end{document}